\documentclass[letterpaper,11pt,oneside]{amsart}

\usepackage{amsmath,amssymb,amsthm}
\usepackage{graphicx}
\usepackage{subfig}
\usepackage{bbm}

\usepackage{algorithm,algpseudocode}

\usepackage{hyperref}
\usepackage[dvipsnames]{xcolor}
\usepackage[mathscr]{euscript}

\theoremstyle{definition}
\theoremstyle{definition}
\theoremstyle{definition}\newtheorem{theorem}{Theorem}
\theoremstyle{definition}
\theoremstyle{definition}\newtheorem{remark}{Remark}
\theoremstyle{definition}\newtheorem{corollary}{Corollary}
\numberwithin{corollary}{section}
\numberwithin{equation}{section}
\numberwithin{theorem}{section}
\numberwithin{definition}{section}
\numberwithin{lemma}{section}
\numberwithin{figure}{section}

\newcommand{\rd}{\mathrm{d}}

\newcommand{\rma}{\mathrm{a}}
\newcommand{\rmb}{\mathrm{b}}

\newcommand{\rmi}{\mathrm{i}}
\newcommand{\rmo}{\mathrm{o}}
\newcommand{\rmr}{\mathrm{r}}
\newcommand{\rmL}{\mathrm{L}}

\newcommand{\eps}{\epsilon}

\newcommand{\Res}{\mathrm{res}}

\newcommand{\rmNN}{\mathrm{NN}}
\newcommand{\rminput}{\mathrm{input}}

\newcommand{\Rb}{\mathbb{R}}

\newcommand{\Lc}{\mathcal{L}}

\newcommand{\algrule}[1][.2pt]{\hspace*{-.6in}\hrulefill}
\algdef{SE}[SUBALG]{Indent}{EndIndent}{}{\algorithmicend\ }%
\algtext*{Indent}
\algtext*{EndIndent}

\newcommand{\ol}[1]{\overline{#1}}
\newcommand{\tabincell}[2]{\begin{tabular}{@{}#1@{}}#2\end{tabular}}

\graphicspath{{./figures/}}

\title[NN-Reduced Order Schwarz for Nonlinear Elliptic Equations]{A reduced order Schwarz method for nonlinear multiscale elliptic equations based on two-layer neural networks}

\author{Shi Chen}
\address{Mathematics Department, University of Wisconsin-Madison, Madison, WI 53706.}
\email{schen636@wisc.edu}

\author{Zhiyan Ding}
\address{Mathematics Department, University of Wisconsin-Madison, Madison, WI 53706.}
\email{zding49@math.wisc.edu}

\author{Qin Li}
\address{Mathematics Department and Discovery Institute, University of Wisconsin-Madison, Madison, WI 53706.}
\email{qinli@math.wisc.edu}

\author{Stephen J. Wright}
\address{Computer Sciences Department, University of Wisconsin, Madison, WI 53706.}
\email{swright@cs.wisc.edu}

\date{\today}

\begin{document}

% \subjscheme[2020]{65N55, 35J66, 41A46, 68T07.}
% \keywords{Nonlinear homogenization, multiscale elliptic problem, neural networks, domain decomposition.}

\begin{abstract}
  Neural networks are powerful tools for approximating high dimensional data that have been used in many contexts, including solution of partial differential equations (PDEs).
  We describe a solver for multiscale fully nonlinear elliptic equations that makes use of domain decomposition, an accelerated Schwarz framework, and two-layer neural networks to approximate the boundary-to-boundary map for the subdomains, which is the key step in the Schwarz procedure.
  Conventionally, the boundary-to-boundary map requires solution of boundary-value elliptic problems on each subdomain.
  By leveraging the compressibility of multiscale problems, our approach trains the neural network offline to serve as a surrogate for the usual implementation of the boundary-to-boundary map.
  Our method is applied to a multiscale semilinear elliptic equation and a multiscale $p$-Laplace equation.
  In both cases we demonstrate significant improvement in efficiency as well as good accuracy and generalization performance.

  \medskip


\end{abstract}

\maketitle

\section{Introduction} \label{sec:intro}
Approximation theory plays a key role in  scientific computing, including in the design of numerical PDE solvers.
This theory  prescribes a certain form of ansatz to approximate a solution to the PDE, allowing derivation of an algebra problem whose solution yields the
coefficients in the ansatz.
Various methods are used to fine-tune the process of translation to an algebraic problem, but the accuracy of the computed solution is essentially determined by the the underlying approximation theory.
New approximation methods have the potential to produce new strategies for numerical solution of PDEs.

During the past decade, driven by some remarkable successes in machine learning, neural networks (NNs) have become popular in many contexts.
They are extremely powerful in such areas as computer vision, natural language processing, and games~\cite{LeBeHi:2015deep,GoBeCo:2016deep}.
What kinds of functions are well approximated by NNs, and what are the advantages of using NNs in the place of  more traditional approximation methods?
Some studies~\cite{Ba:1993universal,KlBa:2016risk,EMaWu:2019barron} have revealed that NNs can represent functions in high dimensional spaces very well.
For Barron functions, in particular, unlike traditional approximation techniques that require a large number of parameters (exponential on the dimension), the number of parameter required for a NN to achieve a prescribed accuracy is rather limited. In this sense, NN approximation overcomes the ``curse of dimensionality.'' This fact opens up many possibilities in scientific computing, where the discretization of high dimensional problems often plays a crucial role.
One example is problems from uncertainty quantification, where many random variables are needed to represent a random field, with each random variable essentially adding an extra dimension to the  PDE~\cite{Xi:2010numerical,XiKa:2002wiener,GhSp:2003stochastic,BaTeZo:2004galerkin}.
Techniques that exploit intrinsic low-dimensional structures can be deployed on the resulting high-dimensional problem~\cite{FrScToRa:2005finite,BiSc:2009sparse,BaScZo:2011multi,CoDeSc:2011analytic,HoLiZh:2017exploring}.
Another example comes from  PDE problems in which the medium contains structures at multiple scales or is highly oscillatory, so that traditional discretization techniques require a large number of grid points to achieve a prescribed error tolerance.
Efficient algorithms must then find ways to handle or compress the many degrees of freedom.

Despite the high dimensionality in these examples, successful algorithms have been developed, albeit specific to certain classes of  problems.
With the rise of NN approximations, with their advantages in high-dimensional regimes, it is reasonable to investigate whether strategies based on NNs can be developed that may even outperform  classical strategies.
In this paper, we develop an approach that utilizes a two-layer NN to solve multiscale elliptic PDEs. We test our strategy on two nonlinear problems of this type.

The use of NN in numerical PDE solvers is no longer a new idea.
Two approaches that have been developed are to use NN to approximate the {\em solutions} (\cite{EHaJe:2017deep,EYu:2018deep,SiSp:2018dgm,RaPeKa:2019physics,ZaBaYeZh:2020weak,LiCaXu:2020multi,LiXuZh:2020multi,BeDaGr:2020numerically, KhLuYi:2019solving}) or the {\em solution map} (\cite{FaLiYiZe:2019multiscale,FaBoYi:2019bcr,LuJiKa:2019deeponet,FeFaYi:2020meta,LiKoAzLiBhStAn:2020fourier,WaChLeChEfWh:2020reduced, KhLuYi:2021solving, KhYi:2019switchnet, LiChLu:2018butterfly, XuLiCh:2020butterfly, WuXi:2020data, LoLuMaDo:2018pde, BaHoHiBr:2019learning, DiHeRa:2020controlling,LiLuMa:2020variational}).
Due to the complicated and unconventional nature of approximation theory for NN, it is challenging to perform rigorous numerical analysis, though solid evidence has been presented of the computational efficacy of these  approaches.

The remainder of our paper is organized as follows.
In Section~\ref{sec:elliptic_Schwarz} we formulate the multiscale PDE problem to be studied.
We give an overview of our domain decomposition strategy and the general specification of the Schwarz algorithm.
In Section~\ref{sec:NN}, we discuss our NN-based approach in detail and justify its use in this setting.
We then present our reduced-order Schwarz method based on two-layer neural networks.
Numerical evidence is reported in Section~\ref{sec:numerics}.
Two comprehensive numerical experiments for the semilinear elliptic equation and the $p$-Laplace equation are discussed, and efficiency of the methods is evaluated.
We make some concluding remarks in Section~\ref{sec:conclusion}.

\section{Domain Decomposition and the Schwarz method for multiscale elliptic PDEs} \label{sec:elliptic_Schwarz}

We start by reviewing some key concepts.
Section~\ref{sec:elliptic} describes nonlinear multiscale elliptic PDEs and discussed the homogenization limit for highly oscillatory medium.
Section~\ref{sec:Schwarz} outlines the domain decomposition framework and the Schwarz iteration strategy.

\subsection{Nonlinear elliptic equation with multiscale medium} \label{sec:elliptic}
Consider the following general class of nonlinear elliptic PDEs with Dirichlet boundary conditions:
\begin{equation}\label{eqn:elliptic}
\begin{cases}
F^\eps\left(D^2u^\eps(x),Du^\eps(x),u^\eps(x),x\right) = 0,\quad& x\in\Omega\,, \\
u^\eps(x) = \phi(x),\quad& x\in\partial\Omega\,,
\end{cases}
\end{equation}
where $\Omega\subset\mathbb{R}^d$ is a domain in $d$-dimensional space, $\eps>0$ represents the small scale, and $F^\eps:S^{d\times d}\times\Rb^d\times\Rb\times\Omega \to \Rb$
(where $S^{d\times d}$ denotes the space of real symmetric $d\times d$ matrices) is a smooth function.
To ensure ellipticity, we require for all $(R,p,u,x) \in  S^{d\times d} \times \Rb^d \times \Rb \times \Omega$ that
\[
F^\eps(R+Q,p,u,x) \leq F^\eps(R,p,u,x) \,,
\]
for all nonnegative semidefinite $Q \in S^{d\times d}$.

This class of problems
has fundamental importance in modern science and engineering, in such areas as synthesis of composite materials, discovery of geological structures, and design of aerospace structures.
The primary computational challenges behind all these problems lie in the complicated interplay between the nonlinearity
and the extremely high number of degrees of freedom necessitated by the smallest scale.
We assume that for an appropriately chosen boundary condition $\phi$, the PDE \eqref{eqn:elliptic} has a unique viscosity solution $u^\eps\in C(\overline{\Omega})$.
For details on the theory of fully nonlinear elliptic equations, see for example, \cite{CaCa:1995fully,IsLi:1990viscosity}.

To achieve a desired level of numerical error, classical numerical methods require refined discretization strategies with a mesh width $\Delta x = o(\eps)$, making the leading to at least $O(\eps^{-d})$ degrees of freedom in the discretized problem.
The resulting numerical cost is prohibitive when $\eps$ is small.
The homogenization limit of \eqref{eqn:elliptic} as $\eps\to 0$ can be specified under additional assumptions, such as when the medium is pseudo-periodic.
Let
\begin{equation}\label{eqn:pseudo_periodic}
F^\eps(R,p,u,x) = F\left(R,p,u,x,\frac{x}{\eps}\right)
\end{equation}
for some $F:S^{d\times d}\times\Rb^d\times\Rb\times\Omega\times\Rb^d \to \Rb$ that is periodic in the last argument with period $Y$.
We have the following theorem.
\begin{theorem}[\cite{Ev:1992periodic}, Theorem 3.3]
Suppose that the nonlinear function $F^\eps$ is uniform elliptic and $u\mapsto F^\eps(\cdot,\cdot,u,\cdot)$ is nondecreasing.
Let $F^\eps$ be pseudo-periodic as defined in \eqref{eqn:pseudo_periodic}.
The solution $u^\eps$ to \eqref{eqn:elliptic} converges uniformly as $\epsilon \to 0$ to the unique viscosity solution $u^\ast$ of the following equation
\begin{equation}\label{eqn:effective}
\begin{cases}
\bar{F}(D^2u^\ast(x),Du^\ast(x),u^\ast(x),x) = 0 , \quad  &x \in \Omega \,, \\
u^\ast(x) = \phi(x),  \quad &x \in \partial\Omega \,,
\end{cases}
\end{equation}
where the homogenized nonlinear function $\bar{F}(R,p,u,x)$ is defined as follows: For a fixed set of $(R,p,u,x)\in S^{d\times d}\times\Rb^d\times\Rb\times\Omega$, there exists a unique real number $\lambda$ for which the following cell problem has a unique viscosity solution $v\in C^{1,\gamma}(\Rb^d)$ for some $\gamma>0$:
\begin{equation}\label{eqn:homogenized_F}
\begin{cases}
F(D_y^2v(y)+R,p,u,x,y) = \lambda, \quad &y \in \Rb^d \,, \\
v(y+Y) = v(y), \quad &y\in\Rb^d \,,
\end{cases}
\end{equation}
(where $Y$ is the period in the last argument of $F$).
We set $\bar{F}(R,p,u,x)=\lambda$.
\end{theorem}

This result
can be viewed as the extension of a linear homogenization result~\cite{Al:1992}.
Although the medium is highly oscillatory for small $\epsilon$,  the solution $u^\epsilon$ approaches that of a certain limiting equation with a one-scale structure, as $\epsilon \to 0$.
In practice, the form of the limit $\bar{F}$ is typically unknown, but this observation has led to an exploration of numerical homogenization algorithms, in which one seeks to capture the limit numerically without resolving the fine scale $\epsilon$.
We view this problem as one of manifold reduction.
The solution $u^\epsilon$ can be ``compressed'' significantly;  its ``information'' is stored mostly in $u^\ast$, which can be computed from \eqref{eqn:homogenized_F} using mesh width $\Delta x = O(1)$, in contrast to the $\Delta x = o(\eps)$ required to solve  \eqref{eqn:elliptic}.
In other words, the $O(\eps^{-d})$-dimensional solution manifold can potentially be compressed into an $O(1)$-dimensional solution manifold, up to small homogenization error that vanishes as $\eps \to 0$.

\begin{remark}
Due to the popularity of the elliptic multiscale problem, the literature is rich.
For {\em linear} elliptic PDEs, many influential methods have been developed, including the multiscale finite element method (MsFEM) \cite{HoWu:1997,EfHoWu:2000,HoWuCa:1999}, the heterogeneous multiscale method (HMM) \cite{EWEn:2003,AbSc:2005,EMiZh:2005}, the generalized finite element method \cite{BaMe:1997,BaLi:2011}, localization methods~\cite{MaPe:2014}, methods based on random SVD \cite{ChLiLuWr:2020randomized,ChLiLuWr:2020random,ChLiLuWr:2021low,ChLiWr:2019schwarz}, and many others \cite{AbBaVi:2015,AbBa:2012,OwZh:2007,OwZhBe:2014polyharmonic,Ow:2017multigrid,Be:2007,Ha:2015}.
Many of these methods adopt an offline-online strategy.
In the offline stage, local bases that encode the small-scale information and approximate the local solution manifold (space) with few degrees of freedom are constructed.
In the online stage, the offline bases are used to compute global solutions on coarse grids, thus reducing online computation requirements drastically over naive approaches.
For {\em nonlinear} problems, there is less prior work, and almost all methods can be seen as extensions of classical methods \cite{EfHoGi:2004,ChSa:2008,EfHo:2009,EMiZh:2005,AbVi:2014,EfGaLiPr:2014,HeMaPe:2014,AbBaVi:2015,LiChZh:2021iterated}.
There is no counterpart on the nonlinear solution manifold for a linear basis, so most classical solvers construct local basis function iteratively, which accounts for a large amount of overhead time.
One strategy that avoids repeated online computation of local bases is to adopt an idea from manifold learning \cite{ChLiLuWr:2020manifold} based on preparing a dictionary for each local patch in the offline stage to approximate the local solution manifold.
The major computational issue for classical multiscale solvers is thus greatly alleviated:  Repeated basis computation is reduced to basis searching on the manifold.
However, since the method is locally linear, its efficacy depends on the amount of nonlinearity of underlying PDE.
Thus far, the approach is difficult to generalize to fully nonlinear elliptic PDEs, and a more universal methodology is needed to approximate the nonlinear solution map.
\end{remark}

\subsection{Domain decomposition and Schwarz iteration} \label{sec:Schwarz}

A popular framework for solving elliptic PDEs is domain decomposition, where the problem is decomposed and  solved separately in different subdomains, with boundary conditions chosen iteratively to ensure regularity of the solution across the full domain.
This approach is naturally parallelizable, with potential savings in memory and computational cost.
It essentially translates the inversion of a large matrix into the a composition of inversions of many smaller matrices.
The many variants of domain decomposition include the Schwarz iteration strategy that we adopt in this paper.
This strategy makes use of a partition-of-unity function that resolves the mismatch between two solutions in adjacent subdomains.
We briefly review the method here.

For simplicity we describe the case of  $d = 2$ and assume throughout the paper that $\Omega = [0,L]^2$ for some $L>0$.
The approach partitions the domain $\Omega$ into multiple overlapping subdomains, also called {\em patches}.
It starts with an initial guess of the solution on the boundaries of all subdomains, and solves the Dirichlet problem on each patch.
The computed solutions then serve as the boundary conditions for neighboring patches, for purposes of computing the next iteration.
The entire process is repeated until convergence.

\begin{figure}[htbp]
  \centering
  \includegraphics[width=0.45\textwidth]{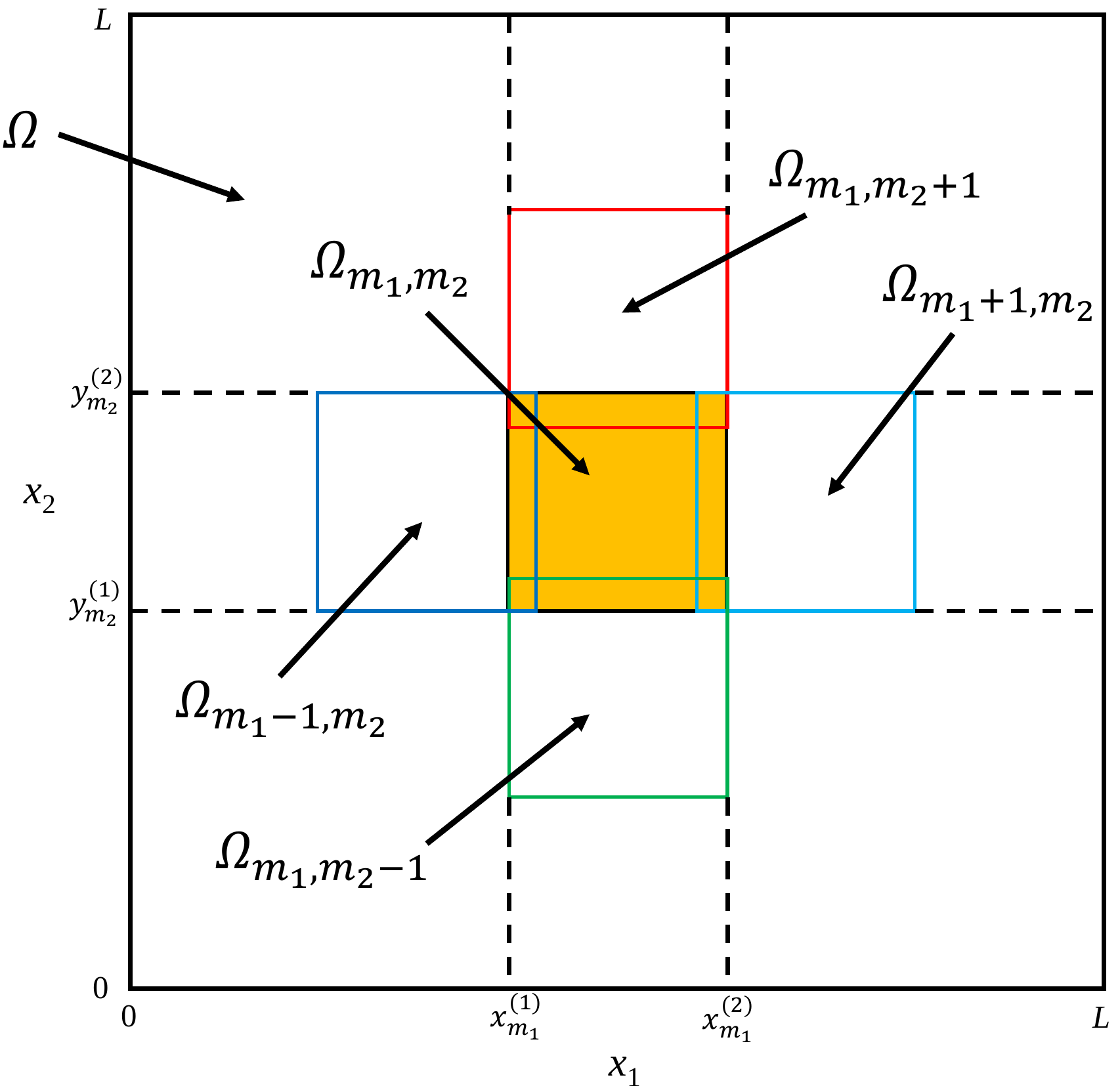}
  \caption{Domain decomposition for a square 2D geometry. Each patch is labeled by a multi-index $m = (m_1,m_2)$. The patches adjacent to $\Omega_m$ are those on its north/south/west/east sides.}
  \label{fig:elliptic_decomp}
\end{figure}

In the current setting, the overlapping rectangular patches are defined as follows:
\begin{equation}
\Omega = \bigcup_{m\in J} \Omega_m,\quad \text{with}\quad \Omega_m = (x_{m_1}^{(1)},x_{m_1}^{(2)})\times(y_{m_2}^{(1)},y_{m_2}^{(2)})\,,
\end{equation}
where $m = (m_1,m_2)$ is a multi-index and $J$ is the collection of the indices
\[
J = \{m = (m_1,m_2):\,m_1 = 1,\dotsc,M_1,\,m_2 = 1,\dotsc,M_2\}\,.
\]
We plot the setup in Figure~\ref{fig:elliptic_decomp}. For each patch we define the associated partition-of-unity function $\chi_m$, which has $\chi_m(x) \geq 0$ and
\begin{equation}\label{eqn:POU_elliptic}
\chi_m(x) = 0\quad \text{on}\,\, x \in \Omega\backslash\Omega_m\,,\quad \sum_m\chi_m(x) = 1, \quad \forall x \in \Omega\,.
\end{equation}
We set $\partial\Omega_m$ to be the boundary of patch $\Omega_m$ and denote by $\mathscr{N}(m)$ the collection of indices of the neighbors of $\Omega_m$. In this 2D case, we have
\begin{equation}\label{eqn:neighbors}
\mathscr{N}(m) = \{(m_1\pm 1, m_2)\} \cup \{(m_1,m_2\pm 1)\} \subset J\,.
\end{equation}
Naturally, indices that are out of range, which correspond to patches adjacent to the boundary $\partial\ \Omega$, are omitted from $\mathscr{N}(m)$.

In the framework of domain decomposition, the full-domain problem is decomposed into multiple smaller problems supported on the subdomains. Define the local Dirichlet problem on patch $\Omega_m$ by:
\begin{equation}\label{eqn:general_local}
  \begin{cases}
  F^\eps\left(D^2u_m^\eps(x),Du_m^\eps(x),u_m^\eps(x),x\right) = 0\,,&\quad x\in \Omega_m\,, \\
  u_m^\eps(x) = \phi_m(x)\,,&\quad x\in\partial\Omega_m\,.
  \end{cases}
\end{equation}
For this local problem, we define the following operators:
\begin{itemize}
  \item $\mathcal{S}_m^\eps$ is the solution operator that maps local boundary condition $\phi_m$ to the local solution $u_m^\eps$:
  \[
  u_m^\eps =  \mathcal{S}_m^\eps \phi_m\,.
  \]
  Denoting by $d_m$ the number of grid points on the boundary $\partial\Omega_m$ and $D_m$ the number of grid points on the subdomain $\Omega_m$, then $\mathcal{S}_m^\eps$ maps  $\Rb^{d_m}$ to $\Rb^{D_m}$.

  \item $\mathcal{I}_{l,m}$ denotes the restriction (or trace-taking) operator that restricts the solution within $\Omega_m$ to its part that overlaps with the boundary of $\Omega_l$, for all $l\in\mathscr{N}(m)$. That is,
  \[
  \mathcal{I}_{l,m} u_m^\eps = u_m^\eps|_{\partial\Omega_l\cap\Omega_m} \,.
  \]
  Denoting by $p_{l,m}$ the number of grid points in $\partial\Omega_l\cap\Omega_m$, then $\mathcal{I}_{l,m}$ maps $\Rb^{D_m}$ to $\Rb^{p_{l,m}}$.

  \item $\mathcal{Q}_{l,m}^\eps$ is the composition of $\mathcal{S}_m^\eps$ and $\mathcal{I}_{l,m}$. It is a boundary-to-boundary operator that maps the local boundary condition $\phi_m$ to the restricted solution $u_m^\eps|_{\partial\Omega_l\cap\Omega_m}$:
  \[
  \mathcal{Q}_{l,m}^\eps \phi_m = \mathcal{I}_{l,m} \mathcal{S}_m^\eps \phi_m = u_m^\eps|_{\partial\Omega_l\cap\Omega_m} \,.
  \]
  $\mathcal{Q}_{l,m}^\eps$ maps $\Rb^{d_m}$ to $\Rb^{p_{l,m}}$.

  \item $\mathcal{Q}_m^\eps$ denotes the collection of all segments of boundary conditions $\psi_{l,m}$ that is computed from the full-domain boundary condition $\phi_m$:
  \begin{equation}\label{eqn:def_Q_m}
    \mathcal{Q}_m^\eps \phi_m = \bigoplus_{l\in\mathscr{N}(m)} \psi_{l,m} = \bigoplus_{l\in\mathscr{N}(m)} \mathcal{Q}_{l,m}^\eps \phi_m
    = \bigoplus_{l\in\mathscr{N}(m)} \mathcal{I}_{l,m} \mathcal{S}_m^\eps \phi_m \,.
  \end{equation}
  Letting $p_m = \sum_{l\in\mathscr{N}(m)}p_{l,m}$,  $\mathcal{Q}_m^\eps$ maps $\Rb^{d_m}$ to $\Rb^{p_m}$.
\end{itemize}

The Schwarz procedure  starts by making a guess of boundary condition on each $\Omega_m$.
At the $n$th iteration, \eqref{eqn:general_local} is solved for each  subdomains $\Omega_m$ (possibly in parallel) and these solutions are used to define new  boundary conditions for the neighboring subdomains $\Omega_l$, $l \in \mathscr{N}(m)$. The boundary conditions for $\Omega_m$ at iteration $n+1$ are thus:
\begin{equation}\label{eqn:patch2_general}
\phi_{m}^{(n+1)} =
\begin{cases}
\psi_{m,l}^{(n)}=\mathcal{I}_{m,l} \mathcal{S}_l^\eps \phi_l^{(n)}\,, & \; \text{on } \partial\Omega_m\cap\Omega_{l}, \quad l\in\mathscr{N}(m)\,,\\
\phi|_{\partial\Omega_m\cap\partial\Omega}\,, &\; \text{on } \partial\Omega_m\cap\partial\Omega\,.
\end{cases}
\end{equation}
Note that  the physical full-domain boundary condition is imposed on the points in $\partial \Omega_m \cap \partial \Omega$.
Each iteration of the Schwarz procedure can be viewed as an application of the map $\mathcal{Q}_{m,l}^\eps$.
The procedure concludes by patching up the local solutions from the subdomains.
The overall algorithm is summarized in Algorithm~\ref{alg:general}.

The convergence of classical Schwarz iteration is guaranteed for fully nonlinear elliptic equations; see, for example \cite{Li:1989schwarz,Li:1988schwarz,Ga:2008schwarz}.
Since the computation of solution $u_m^\eps = \mathcal{S}_m^\eps\phi_m$ can be expensive due to the nonlinearity and oscillation of the medium at small scale $\eps$, the major computational cost for Schwarz iteration comes from the repeated evaluation of the boundary-to-boundary map $\mathcal{Q}_{m,l}^\eps$, which requires solution of an elliptic PDE on each subdomain.

\begin{algorithm}
\caption{The Schwarz iteration for fully nonlinear elliptic equations~\eqref{eqn:elliptic}.}\label{alg:general}

\begin{algorithmic}[1]
\State {\bf Domain Decomposition:}
\State \hspace*{0.2in} Decompose $\Omega$ into overlapping patches: $\Omega=\bigcup_{m\in J}\Omega_m$.
\State Given tolerance $\delta_0$ and initial guesses $\phi_{m}^{(0)}$ of boundary conditions on each patch $m \in J$.
\State {\bf Schwarz iteration: }
\Indent
\State Set $n = 0$ and $\Res = 1$.
\While{$\Res\geq\delta_0$}
    \State For $m \in J$, compute local solutions $u_{m}^{(n)} =\mathcal{S}_m \phi_{m}^{(n)}$;
    \State For $m \in J$ and $l\in\mathscr{N}(m)$, restrict the solutions $\psi_{m,l}^{(n)} = \mathcal{I}_{m,l} u_{m}^{(n)}$;
    \State For $m \in J$, update $\phi_m^{(n+1)}$ by~\eqref{eqn:patch2_general};
    \State Set $\Res = \sum_{m}\|\phi_{m}^{(n+1)}-\phi_{m}^{(n)}\|_{L^2(\partial\Omega_m)}$ and $n \gets n+1$.
\EndWhile
\EndIndent \\
\Return Global solution $u^{(n)} = \sum_{m\in J} \chi_m u_m^{(n)}$.
\end{algorithmic}
\end{algorithm}

\section{Reduced order Schwarz method based on neural networks} \label{sec:NN}
The major numerical expense in the Schwarz iteration comes from the local PDE solves --- one per subdomain per iteration.
However, except at the final step where we assemble the global solution, our interest is not in the local solutions per se: It is in the boundary-to-boundary maps that share information between adjacent subdomains on each Schwarz iteration.
If we can implement  these maps {\em directly}, we can eliminate the need for local PDE solves.
To this end, we propose an offline-online procedure.
In the offline stage, we implement the boundary-to-boundary maps, and in the online stage, we call these maps repeatedly in the Schwarz framework.
This approach is summarized in Algorithm~\ref{alg:NN_online}.
In this description, we replace the boundary-to-boundary map $\mathcal{Q}^\eps_m$ by a surrogate $\mathcal{Q}_m^{\rmNN}(\theta_m)$, which is neural network parametrized by weights $\theta_m$, whose values are found by an offline training process.

\begin{algorithm}
\caption{The NN-Schwarz iteration for nonlinear elliptic equations~\eqref{eqn:elliptic}.}\label{alg:NN_online}

\begin{algorithmic}[1]

\State {\bf Domain Decomposition: }
\State \hspace*{0.2in} Decompose $\Omega$ into overlapping patches: $\Omega=\bigcup_{m\in J}\Omega_m$, and collect the indices for interior patches in $J_\rmi = \{m\in J: \partial\Omega_m\cap\partial\Omega=\varnothing\}$ and boundary patches in $J_\rmb = \{m\in J: \partial\Omega_m\cap\partial\Omega\neq\varnothing\}$.
\State {\bf Offline training:}
\State \hspace*{0.2in}  For each interior patch $\Omega_m$, train the boundary-to-boundary map $\mathcal{Q}_m^\rmNN(\theta_m)$ parametrized by $\theta_m$.

\State {\bf Schwarz iteration (Online):}
\Indent
\State Given the tolerance $\delta_0$ and the initial guess of boundary conditions $\phi_{m}^{(0)}$ on each patch $m \in J$.
\State Set $n = 0$ and $\Res = 1$.
\While{$\Res\geq\delta_0$}
    \State For $m \in J_\rmi$, compute function $(\psi_{m,l}^{(n)})_{l\in\mathscr{N}(m)} = \mathcal{Q}_m^{\rmNN}(\theta_m) \phi_m^{(n)}$;
    \State For $m \in J_\rmb$, compute function $\psi_{l,m}^{(n)} = \mathcal{I}_{m,l}\mathcal{S}_m^\eps \phi_{m}^{(n)}$ for $l\in\mathscr{N}(m)$;
    \State For $m \in J$, update $\phi_m^{(n+1)}$ by~\eqref{eqn:patch2_general};
    \State Set $\Res = \sum_{m}\|\phi_{m}^{(n+1)}-\phi_{m}^{(n)}\|_{L^2(\partial\Omega_m)}$ and $n \gets n+1$.
\EndWhile
\State For $m \in J$, compute function $u_{m}^{(n)} =\mathcal{S}_m^\eps \phi_{m}^{(n)}$;
\EndIndent \\\Return Global solution $u^{(n)} = \sum_{m\in J} \chi_m u_m^{(n)}$.

\end{algorithmic}
\end{algorithm}

Since the online stage is self-explanatory, we focus on the offline stage, and study how to obtain the approximation to $\mathcal{Q}^\eps_m$.

We have two additional comments about our approach.
\begin{itemize}
\item Our algorithm uses the surrogate boundary-to-boundary map only for the interior patches $m\in J_\rmi$. For patches that are adjacent to the physical boundary, we perform the standard Schwarz iteration.
This choice is mainly for convenience of coding.
There are several other options.
For example, one can choose to learn in the offline stage the surrogate boundary-to-boundary map for patches boundary patches $J_\rmb$ as well.
However, in the training stage, one needs to impose a rather general class of functions to serve as the potential physical boundary condition.
Whether this chosen class of functions represents well the boundary condition given in the online phase is a question for approximation theory. We omit a discussion here.
\item If the PDE operator $F^\epsilon$ has no explicit dependence on $x$, then the boundary-to-boundary map is the same across all patches of the same size.
In this case, training can be implemented in parallel, saving computational time.
\end{itemize}

\subsection{Two observations}
A rigorous approach to preparing the boundary-to-boundary map $\mathcal{Q}^\eps_m$ in the offline stage is not straightforward.
In the case of linear PDEs, it amounts to computing all Green's functions in the local subdomains and confining them on the adjacent subdomain boundaries for the map; see~\cite{ChLiWr:2019schwarz}.
When the PDEs are nonlinear,
there would seem to be no alternative to solving the local PDEs with all possible configurations of the boundary conditions, applying the appropriate restrictions, and storing the results.
At the discrete level, $\mathcal{Q}^\eps_m$ would be represented as a high-dimensional function mapping $\Rb^{d_m}$ to  $\Rb^{p_m}$.
To achieve a specified accuracy, both $d_m$ and $p_m$ need to scale as $O(\eps^{-(d-1)})$.
For brute-force training, at least $O(d_m)=O(\eps^{-(d-1)})$ local PDE solves need to be performed to compute the required approximation to $\mathcal{Q}^\eps_m$.
This is a large amount of computation, and it offsets whatever gains accrue in the online stage from efficient deployment of the approximation to $\mathcal{Q}^\eps_m$.

\begin{figure}[htbp]
  \centering
  \includegraphics[width=0.4\textwidth]{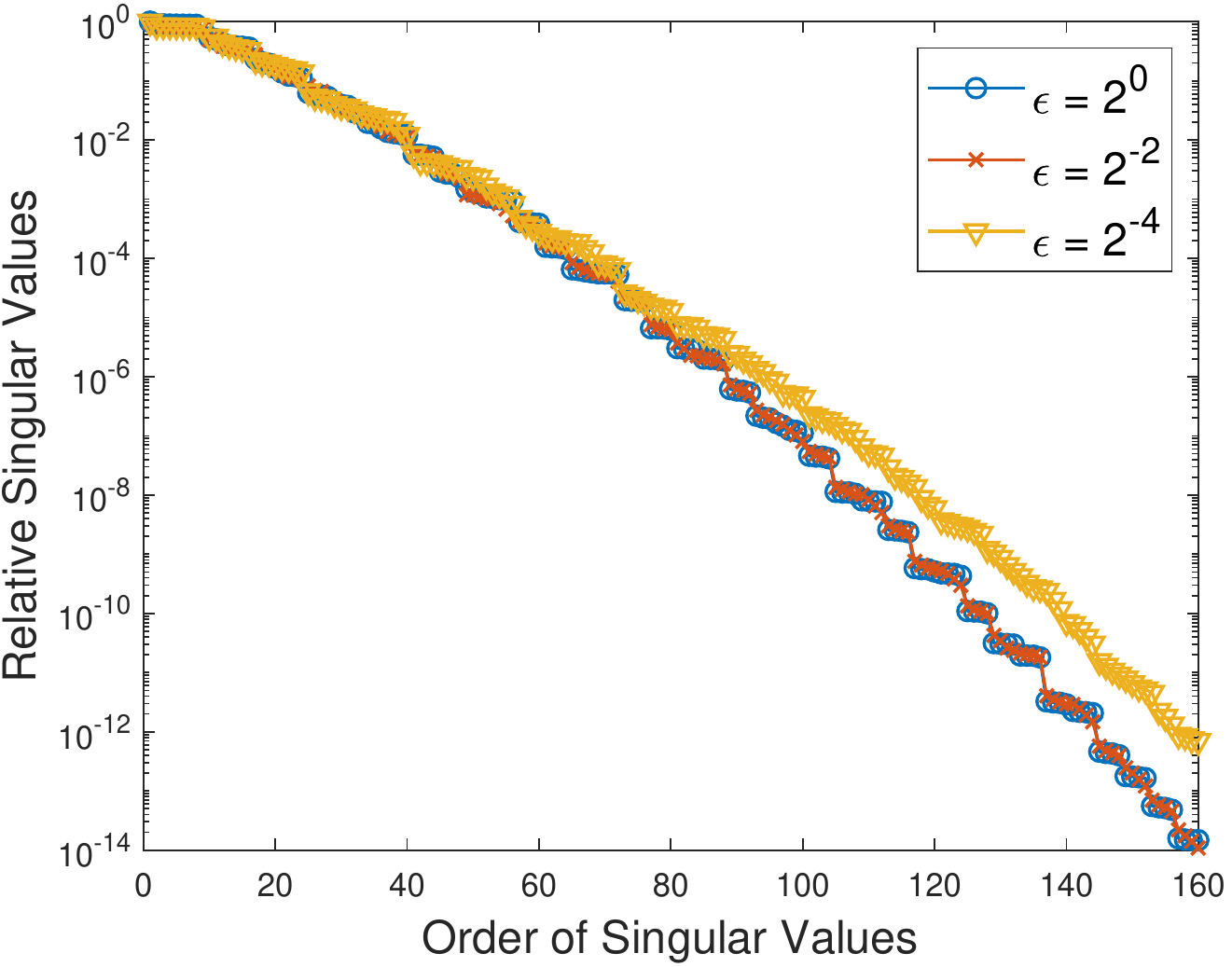}
  \includegraphics[width=0.4\textwidth]{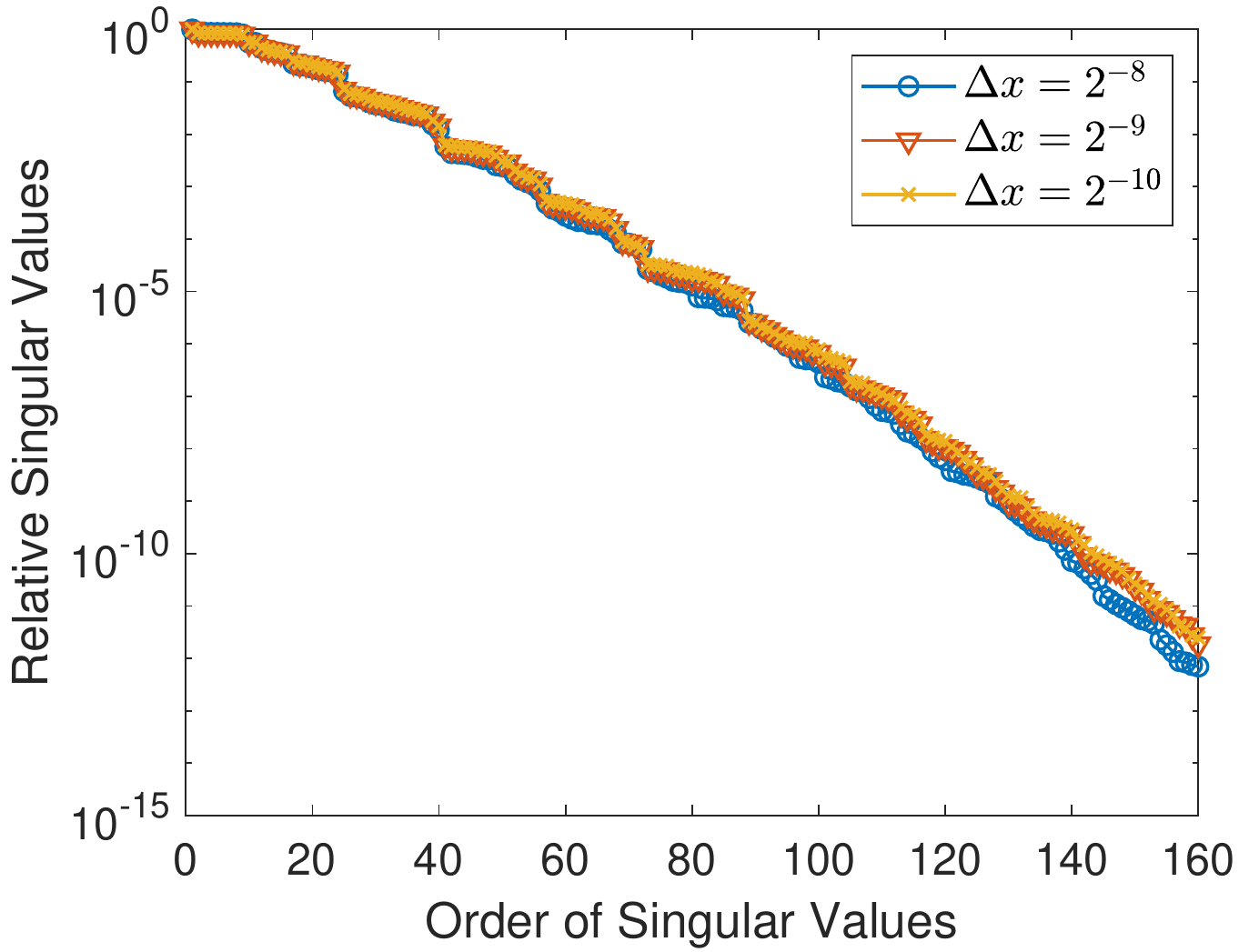}
  \caption{Singular values of the boundary-to-boundary operator $\mathcal{Q}_m^\eps$ for the linear elliptic equation~\eqref{eqn:linear_initial_semi} with medium $\kappa^\eps$ defined in~\eqref{eqn:medium} for different values of $\eps$ and $\Delta x$ on a local patch. Left plot: $\Delta x = 2^{-8}$. Right plot: $\eps = 2^{-4}$. To ensure the regularity of the test function space, the discrete version of the boundary-to-boundary map is represented on basis functions composed of piecewise linear function with fixed step size $2^{-8}$.}
  \label{fig:svd}
\end{figure}

To be cost-effective, a method of the form of Algorithm~\ref{alg:NN_online} must exploit additional properties, intrinsic to $\mathcal{Q}^\eps_m$ and  to the scheme for approximating this mapping.
The first such property is a direct consequence of homogenization.
As argued in Section~\ref{sec:elliptic}, the solution of the effective equation~\eqref{eqn:effective} can preserve the ground truth well, with the effective equation independent of  $\eps$. Therefore, the map $\mathcal{Q}^\eps_m$, though presented as a mapping from $\Rb^{d_m}$ to $\Rb^{p_m}$, is intrinsically of low dimension and can be compressed.
To visualize this relation, we plot the relative singular values of the boundary-to-boundary operator $\mathcal{Q}_m^\eps$ of a linear multiscale elliptic equation (see~\eqref{eqn:linear_initial_semi}) in Figure~\ref{fig:svd}.

With the system being of intrinsically low dimension, we expect that a compression mechanism can be deployed.
Even though the data itself is represented in high dimension, the number of parameters in the compressed representation should not grow too rapidly with the order of discretization. We seek an approximation strategy that can overcome the ``curse of dimensionality.''
These considerations lead us to the use of neural network (NN).
NN, unlike other approximation techniques, is powerful in learning functions supported in high dimensional space; the number of parameters that need to be tuned to fit data in a high dimensional space is typically relaxed from the dimension of the data.

Consider a fully connected feedforward neural network (fully connected NN) representing a function $f:\Rb^{n}\to\Rb^{m}$.
A 2-layer fully connected NN with hidden-layer width $h$ would thus be required to satisfy
\begin{equation} \label{eq:fNN}
f^\rmNN(x) = W_2\sigma(W_1x+b_1)+b_2 \,, \quad x\in\Rb^n\,,
\end{equation}
where $W_1\in\Rb^{h\times n}$, $W_2\in\Rb^{m\times h}$ are weight matrices and $b_1\in\Rb^{h}$, $b_2\in\Rb^{m}$ are biases. The activation function $\sigma:\Rb\to\Rb$ is applied component-wise to its argument.
(The ReLU activation function $\sigma(x) = \max(x,0)$ is especially popular.)
This $2$-layer fully connected NN already can represent high dimensional functions.
A fundamental approximation result \cite{KlBa:2016risk,EMaWoWu:2020towards,Ba:1993universal} is captured in the following theorem.
\begin{theorem}[Barron's Theorem] \label{thm:barron}
Let $D\subset\Rb^n$ be a bounded domain. Suppose a generic function $f\in L^2(D)$ satisfies
\begin{equation}\label{eqn:barron_spectral}
\Delta(f) = \int_{\Rb^n} \|\omega\|_1^2 |\hat{f}(\omega)| \rd \omega < \infty\,,
\end{equation}
where $\hat{f}$ is the Fourier transform of the zero extension of $f$ to $L^2(\Rb^d)$. Then there exists a two-layer ReLU neural network $f^\rmNN$ with $h$ hidden-layer neurons such that
\begin{equation}
\|f-f^\rmNN\|_{L^2(D)} \lesssim \frac{\Delta(f) }{\sqrt{h}} \,.
\end{equation}
\end{theorem}

A natural high dimensional extension of the result is as follows.
\begin{corollary}
Let $D\subset\Rb^n$ be a bounded domain. Suppose a generic function $f=[f_1\,,\cdots,f_m]:\mathbb{R}^n\to\mathbb{R}^m$ so that $f_i\in L^2(D)$ satisfies \eqref{eqn:barron_spectral}, then there exists a two-layer ReLU neural network $f^\rmNN$ with $h$ hidden-layer neurons such that
\begin{equation}
\|f-f^\rmNN\|_{L^2(D)} \lesssim \sqrt{\sum_{i=1}^m\frac{\Delta^2(f_i)}{h/m}}\leq m\frac{\Delta(f)}{\sqrt{h}} \,.
\end{equation}
where $\Delta(f) :=\max_{i=1}^m\Delta(f_i)$.
\end{corollary}
A nice feature of this result is that the approximation error is mostly relaxed from the dimension of the problem, making NN a good fit for our purposes.
In our setting, it is the high-dimensional operator  $\mathcal{Q}_m^\eps$ that needs to be learned.
Theorem~\ref{thm:barron} suggests that if fully connected NN is used as the representation, the number of neurons $h$ required will not depend strongly on this dimension.

\subsection{Offline training and the full algorithm}
The two observations above suggest that using a neural-network approximation for the boundary-to-boundary operator can reduce computation costs and memory significantly.
Following \eqref{eq:fNN}, we define the NN approximation $\mathcal{Q}_m^{\rmNN}$ to $\mathcal{Q}^\eps_m$ as follows:
\begin{equation}\label{eqn:NN_Q}
\mathcal{Q}_m^{\rmNN}(\theta_m) \phi_m = W_{m,2} \sigma(W_{m,1}\phi_m + b_{m,1}) + b_{m,2} \,, \quad \mbox{where $\phi_m \in \Rb^{d_m}$.}
\end{equation}
Here $\theta_m = \{W_{m,1}, W_{m,2}, b_{m,1}, b_{m,2}\}$ denotes all learnable parameters, with weight matrices $W_{m,1} \in \Rb^{h_m\times d_m}, W_{m,2} \in \Rb^{p_m\times h_m}$ and biases $b_{m,1} \in \Rb^{h_m}, b_{m,2}\in \Rb^{p_m}$.
The number of neurons $h_m$ is a tunable parameter that relates to the number of degrees of freedom in  $\mathcal{Q}_m^{\rmNN}(\theta_m)$. Theorem~\ref{thm:barron} and the homogenizability of the elliptic equation suggest that $h_m$ can be chosen to satisfy a prescribed approximation error while being independent of both $d_m$ and $p_m$, and thus of the small scale $\eps$.

Given a fixed NN architecture and a data set, the identification of optimal $\mathcal{Q}_m^{\rmNN}(\theta_m)$ amounts to minimizing a loss function $\Lc(\theta_m)$ that measures the misfit between the data and the prediction.
One needs to prepare a set of data $\mathscr{X}_m =  \left\{ \phi_{m,i} \right\}_{i=1}^N$ and corresponding outputs
\begin{equation}\label{eqn:general_dic1}
\begin{aligned}
\mathscr{Y}_m =  \left\{ \psi_{m,i} = \mathcal{Q}^\eps_m\phi_{m,i}=(\psi_{l,m,i})_{l\in\mathscr{N}(m)} = \left({u}^\eps_{m,i}|_{\partial\Omega_l\cap\Omega_m}\right)_{l\in\mathscr{N}(m)} \right\}_{i=1}^N\,,
\end{aligned}
\end{equation}
where $u^\eps_{m,i}$ solves~\eqref{eqn:general_local}.
The loss function to be minimized is
\begin{equation}\label{eqn:loss}
\begin{aligned}
 \Lc(\theta_m) := \frac{1}{N}\sum_{i=1}^N \ell\left(\mathcal{Q}_m^{\rmNN}(\theta_m)\phi_{m,i}\,,\psi_{m,i}\right)\,,
\end{aligned}
\end{equation}
where $\ell$ evaluates the mismatch between the first and the second arguments. (This measure could be defined using the $L_2$ norm and / or the $H^1$ norm.)
Gradient-based algorithms for minimizing \eqref{eqn:loss} have the general form
\begin{equation}\label{eqn:iteration}
\theta^{(t+1)}_m \leftarrow \theta^{(t)}_m - \eta_t G_t\left(\nabla_{\theta_m}\mathcal{L}\left(\theta_m^{(t)}\right),\dots,\nabla_{\theta_m}\mathcal{L}\left(\theta_m^{(1)}\right)\right) \,,
\end{equation}
where $\eta_t$ is the learning rate
and $G_t$ is based on the all  gradients seen so far.
For example, for the Adam optimizer~\cite{KiBa:2014adam}, the function $G_t$ is a normalized exponentially decaying average of gradients:
\begin{equation}
G_t\left(a_t,\dots,a_1\right)\propto\left(1-\beta_1^t\right)^{-1}\sum_{s=1}^{t}\beta_1^{t-s}\left(1-\beta_1\right)a_{s}\,,
\end{equation}
for some parameter $\beta_1 \in (0,1)$.
The $\propto$ sign means $G_t$ needs to be normalized so that $\|G_t\|_2\sim 1$.

Like many optimization processes, the training and tuning of this NN depends on some prior knowledge.
We propose a mechanism to select training data that represent well the information in $\mathcal{Q}^\eps_m$.
We also initialize the weights $\theta_m$ according to a reduced linear problem.
These mechanisms are described in the following two sections; their effectiveness in numerical testing is demonstrated in Section~\ref{sec:numerics}.

\subsubsection{Generating  training data} \label{sec:training_data}

\begin{figure}[htbp]
  \centering
  \includegraphics[width=0.3\textwidth]{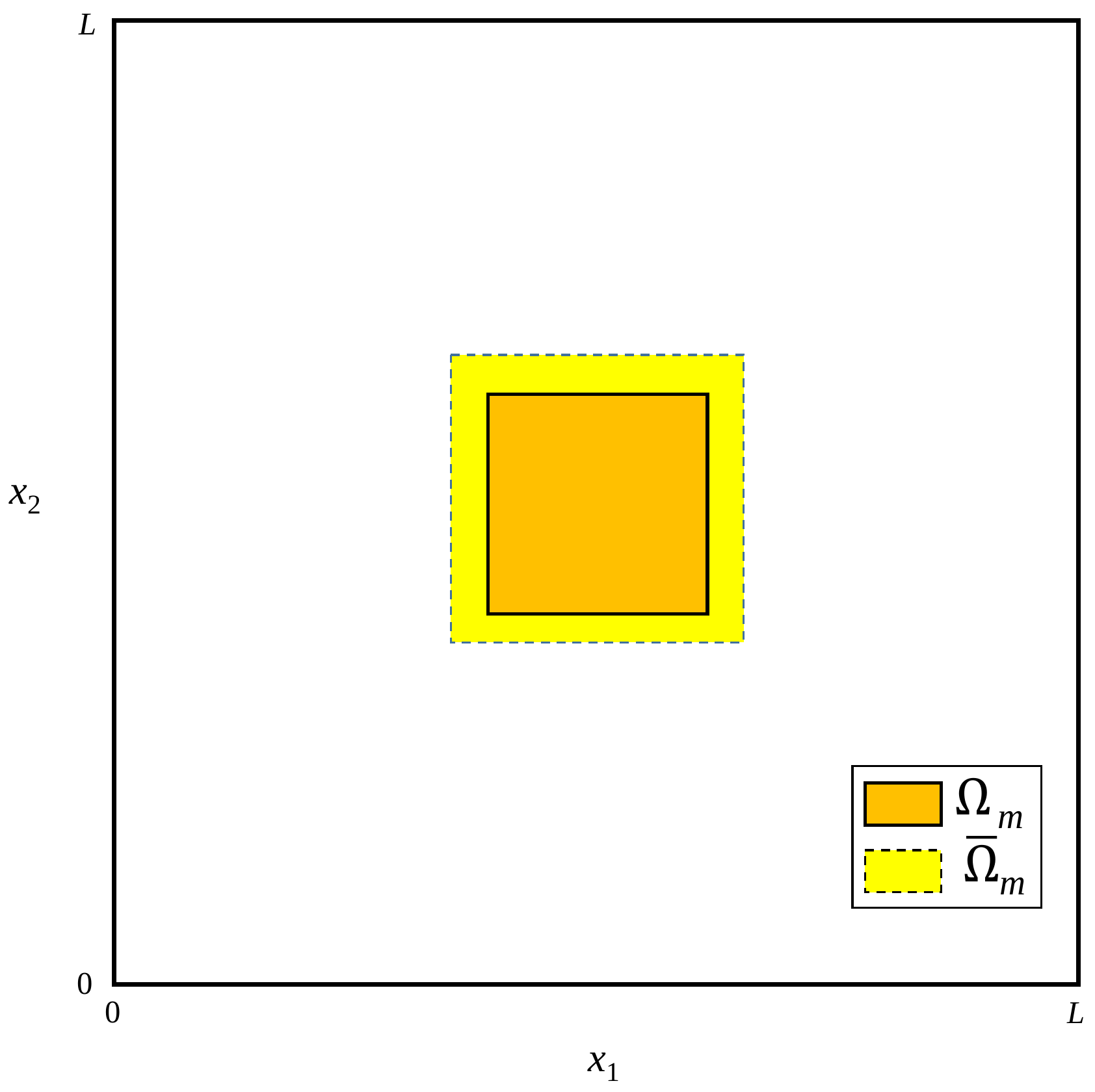}
  \caption{Local enlargement of patches is used to damp boundary effects.}
  \label{fig:elliptic_buffer}
\end{figure}

To learn the parameters in the NN approximation to the boundary-to-boundary map, one needs to provide a training set of examples of the map.
We generate such examples by adding a boundary margin of width $\Delta x_{\text{b}}$ to each interior patch $\Omega_m$  to obtain an enlarged patch $\ol{\Omega}_m$, as shown in Figure~\ref{fig:elliptic_buffer}.
Samples are generated by choosing Dirichlet conditions for the enlarged patch, then solving the equation, and defining the map in terms of restrictions of both input and output conditions to the appropriate boundaries.

Specifically, following~\cite{ChLiLuWr:2020manifold}, we generate $N$ i.i.d. samples of the boundary conditions $\ol{\phi}_m$ for the enlarged patch $\partial\ol{\Omega}_m$ according to $H^{1/2}(\partial\ol{\Omega})$\footnote{The distribution of the sample is uniform in angle and satisfies a power law in the radius. Letting $D>0$, we write $\ol{\phi}_{m} = \ol{\phi}_\rmr\ol{\phi}_\rma$ with $\ol{\phi}_\rma\in\Rb^{\widetilde{d}_m}$ uniformly distributed on the unit sphere $S^{\widetilde{d}_m-1} = \{\phi\in\ol{\mathcal{V}}_m: \|\phi\|_{1/2}=1\}\subset\Rb^{\widetilde{d}_m}$ and $\ol{\phi}_\rmr\in\Rb$ distributed in $[0,R_m]$ according to the density function $f(r) = \frac{D+1}{R_m^D}r^D$. To measure the discrete $H^{1/2}$ norm, we employ the formula $\|\ol{\phi}_m\|_{1/2}^2 = \Delta x \sum_{i=1}^{\ol{d}_m} |(\ol{\phi}_m)_i|^2 + (\Delta x)^2 \sum_{\substack{i,j = 1,\dots,\ol{d}_m\\i\neq j}} \frac{|(\ol{\phi}_m)_i-(\ol{\phi}_m)_j|^2}{|x_i-x_j|^2}$, where we denote denote $\ol{\phi}_m = \left((\ol{\phi}_m)_1,\dots,(\ol{\phi}_m)_{\ol{d}_m}\right)^\top$, and $\Delta$ is the step size.}, and solve the following equations for  $\ol{u}^\eps_{m,i}(x)$:
\begin{equation}\label{eqn:general_local_buffer}
\begin{cases}
F^\eps\left(D^2\ol{u}^\eps_{m,i}(x), D\ol{u}^\eps_{m,i}(x), \ol{u}^\eps_{m,i}(x), x\right) = 0\,,&\quad x\in \ol{\Omega}_m\,, \\
\ol{u}^\eps_{m,i}(x) = \ol{\phi}_{m,i}(x)\,,&\quad x\in\partial \ol{\Omega}_m\,.
\end{cases}
\end{equation}
The boundary-to-boundary map $\mathcal{Q}^\eps_m$ maps each element of  $\mathscr{X}_m=\{\phi_{m,i}\}_{i=1}^N$ to the corresponding element of $\mathscr{Y}_m=\{\psi_{m,i}\}_{i=1}^m$, where
\begin{equation}\label{eqn:general_dic}
\phi_{m,i} = \ol{u}^\eps_{m,i}|_{\partial\Omega_m}\,, \quad \psi_{m,i} = (\psi_{l,m,i})_{l\in\mathscr{N}(m)} = \left(\ol{u}^\eps_{m,i}|_{\partial\Omega_l\cap\Omega_m}\right)_{l\in\mathscr{N}(m)}\,.
\end{equation}
This pair of sets --- input set $\mathscr{X}_m$ and output set $\mathscr{Y}_m$ --- serves as the training data.

We have two  comments regarding the training process.
\begin{itemize}
\item Various  NN architectures could be considered. We use fully connected NN mostly because the initialization procedure requires singular value decomposition of the linearized counterpart of boundary-to-boundary map, and this kind of NN is a natural extension of the linear network that capture the SVD to the nonlinear regime. Another potentially good option is Convolutional Neural Network (CNN) whose structure can potentially alleviate computational difficulty as one refines the discretization.
\item For generating training data, each patch is slightly enlarged before application of the PDE solver.
Most homogenization results need a boundary layer correction adjacent to the physical boundary, so the low-rank property fails to hold near the boundary. The use of a small boundary buffer zone on each patch dampens the boundary layer effect and enables low-rank structure of the boundary-to-boundary map.
\end{itemize}

\subsubsection{Initialization} \label{sec:implement}

The training problem of minimizing $\Lc(\theta_m)$ in \eqref{eqn:loss}  to obtain the NN approximate operator $\mathcal{Q}_m^{\rmNN}(\theta_m)$ is nonconvex, so a  good initialization scheme can improve the performance of a gradient-based optimization scheme significantly.
We can make use of knowledge about the PDE to obtain good starting points.
Our strategy is to assign good initial weights and biases for the neural network using a {\em linearization} of the fully nonlinear elliptic equation \eqref{eqn:elliptic}.
Denoting by $\mathcal{Q}_m^\rmL$ the boundary-to-boundary operator of a linearized version of $\mathcal{Q}^\eps_m$, to be made specific below for the numerical examples in Section~\ref{sec:numerics}, we initialize $\mathcal{Q}_m^\rmNN$ in a way that approximately captures $\mathcal{Q}_m^\rmL$.
The linear boundary-to-boundary operator $\mathcal{Q}_m^\rmL$ has a matrix representation.
Denoting by $r_m$ the approximate rank (up to a preset error tolerance), we can write
\begin{equation}\label{eqn:Q_m^L_rank}
\mathcal{Q}_m^\rmL \approx U_{m,r_m} \Lambda_{r_m} V_{m,r_m}^\top = \left(U_{m,r_m} \sqrt{\Lambda_{r_m}}\right) \left(V_{m,r_m}\sqrt{\Lambda_{r_m}}\right)^\top \,,
\end{equation}
where $U_{m,r_m}\in\Rb^{p_m\times r_m}$ and $V_{m,r_m}\in\Rb^{d_m\times r_m}$ have orthonormal columns while $\Lambda_{r_m}\in\Rb^{r_m\times r_m}$ is diagonal.
As argued in~\cite{ChLiWr:2019schwarz}, due to the fact that the underlying equation is homogenizable, this rank $r_m$ is much less than $\min\{d_m,p_m\}$, and is independent of $p_m$ and $d_m$.

To start the iteration of $\mathcal{Q}_m^\rmNN$, we compare~\eqref{eqn:NN_Q} with the form of~\eqref{eqn:Q_m^L_rank}. This suggests the following settings of parameters in \eqref{eqn:NN_Q}: $b_{m,1} = b_{m,2} = 0$ and
\begin{equation}\label{eqn:initial_weights}
\begin{aligned}
    W_{m,1} & = \left[V_{m,r_m}\sqrt{\Lambda_{r_m}},-V_{m,r_m}\sqrt{\Lambda_{r_m}}\right]^\top, \\
    W_{m,2} & = \left[U_{m,r_m} \sqrt{\Lambda_{r_m}},-U_{m,r_m} \sqrt{\Lambda_{r_m}}\right] \,.
    \end{aligned}
\end{equation}
Note that $h_m = 2r_m$.
These configurations will be used as the initial iteration in~\eqref{eqn:iteration}.

We summarize our offline training method in Algorithm~\ref{alg:NN_offline}.
Integrating into the full algorithm yields the reduced order neural network based Schwarz iteration method.

\begin{algorithm}
\caption{Offline training of $\mathcal{Q}_m^{\rmNN}(\theta_m)$, as a surrogate of $\mathcal{Q}^\eps_m$ on patch $\Omega_m$.}\label{alg:NN_offline}

\begin{algorithmic}[1]
\State Enlarge each interior patch $\Omega_m$ to obtain $\ol{\Omega}_m$;
\State Randomly generate samples $\{\ol{\phi}_{m,i}\}_{i=1}^N$ and solve~\eqref{eqn:general_local_buffer} to obtain $\{\ol{u}^\eps_{m,i}\}_{i=1}^N$.
\State Compute~\eqref{eqn:general_dic} to define $\{\mathscr{X}_m\,,\mathscr{Y}_m\} = \{\{\phi_{m,i}\}_{i=1}^N\,, \{(\psi_{l,m,i})_{l\in\mathscr{N}(m)}\}_{i=1}^N\}$;
\State Initialize $\theta_m$ in $\mathcal{Q}_m^{\rmNN}(\theta_m)$ by using the linearized boundary-to-boundary operator $\mathcal{Q}_m^{\rmL}$, as defined in \eqref{eqn:initial_weights};
\State Find the optimal coefficient $\theta_m^\ast$ in the neural network $\mathcal{Q}_m^{\rmNN}(\theta_m)$ by  applying the gradient descent method~\eqref{eqn:iteration} until convergence.
\end{algorithmic}
\end{algorithm}

\section{Numerical results} \label{sec:numerics}
We present numerical examples using our proposed method to solve a multiscale semilinear elliptic equation and a multiscale $p$-Laplace equation.
In both examples, we use domain $\Omega = [0,1]^2$.
To form the partitioning, $\Omega$ is divided into $M_1\times M_2$ equal non-overlapping rectangles, then each rectangle is enlarged by $\Delta x_{\rmo}$ on the sides that do not intersect with $\partial\Omega$, to create overlap.
We thus have
\begin{align*}
\Omega_m & =\left[\max\left(\tfrac{m_1-1}{M_1}-\Delta x_{\text{o}},0\right),\min\left(\tfrac{m_1}{M_1}+\Delta x_{\text{o}},1\right)\right] \\
& \quad \times \left[\max\left(\tfrac{m_2-1}{M_2}-\Delta x_{\text{o}},0\right),\min\left(\tfrac{m_2}{M_2}+\Delta x_{\text{o}},1\right)\right]\,,\quad m = (m_1,m_2) \in J\,.
\end{align*}
The loss function is defined as in~\eqref{eqn:loss}, with parameter $\mu=10^{-3}$.
For training to obtain  $\mathcal{Q}_m^{\rmNN}(\theta_m)$, we use PyTorch~\cite{PaGrMa:2019pytorch}. For both examples, each neural network is trained for 5,000 epochs using shuffled mini-batch gradient descent with a batch-size of 5\% of the training set size. The Adam optimizer is used with default settings, and the learning rate decays with a decay-rate of 0.9 every 200 epochs. The codes accompanying this manuscript are publicly available~\cite{ChDiLiWr:2021reduced}.

\begin{figure}[b!]
  \centering
  \includegraphics[width=0.35\textwidth]{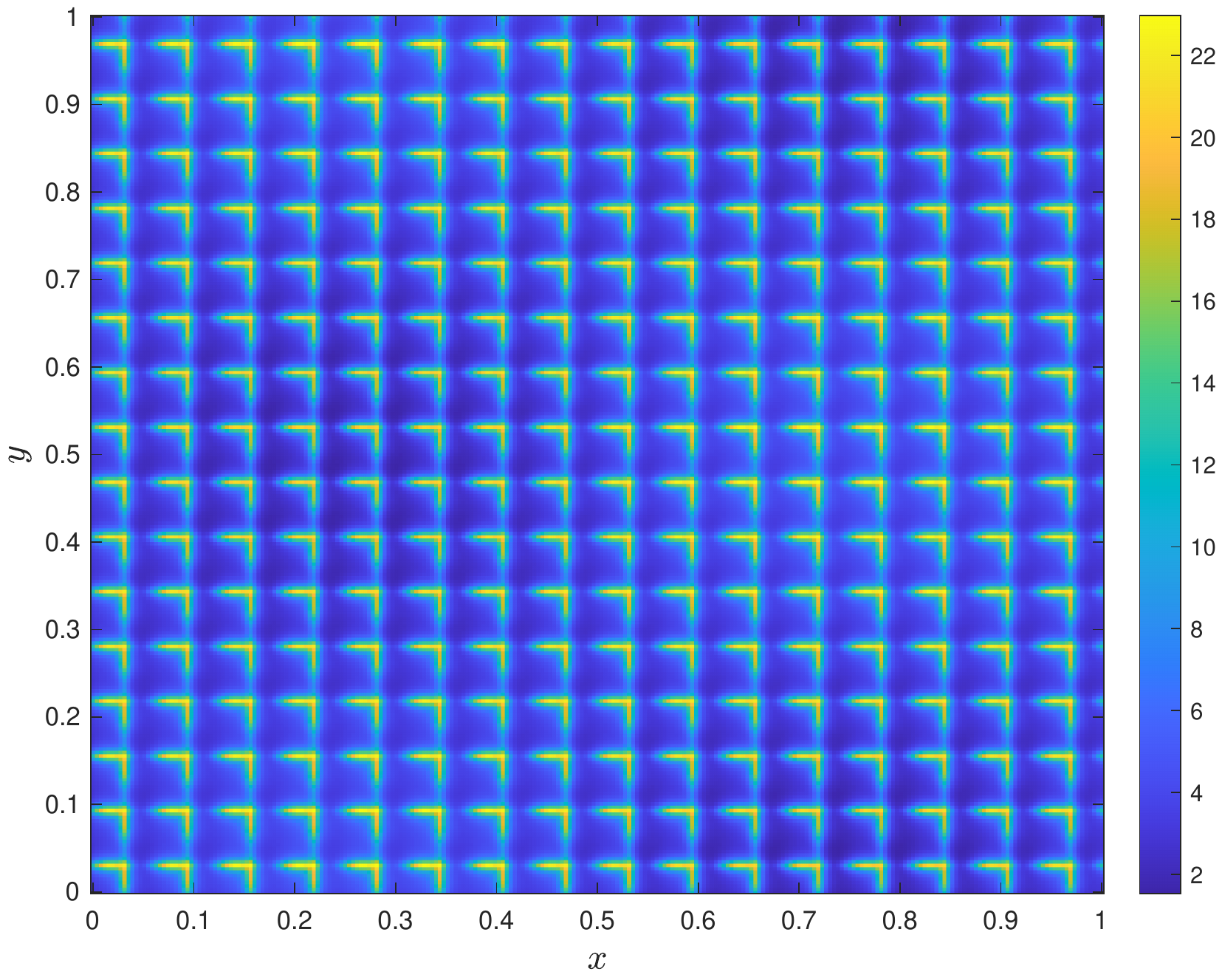}
  \caption{Medium $\kappa$ for semilinear elliptic equation.}
  \label{fig:semilinear_medium}
\end{figure}

\subsection{Semilinear elliptic equations} \label{sec:numerics_semi}
The first example is the semilinear elliptic equation
\begin{equation}\label{eqn:semi_elliptic}
\begin{cases}
-\nabla \cdot (\kappa^\eps(x) \nabla u^\eps(x)) + u^\eps(x)^3 = 0,\quad& x\in\Omega\,, \\
u^\eps(x) = \phi(x),\quad& x\in\partial\Omega\,,
\end{cases}
\end{equation}
with oscillatory medium $\kappa^\eps(x)=\kappa^\eps(x_1,x_2)$ defined by
\begin{equation} \label{eqn:medium}
\kappa^\eps(x_1,x_2) = 2+\sin(2\pi x_1)\cos(2\pi x_2)+\frac{2+1.8\sin(2\pi
  x_1/\eps)}{2+1.8\cos(2\pi x_2 / \eps)} + \frac{2+\sin(2\pi x_2
  /\eps)}{2+1.8\cos(2\pi x_1/\eps)}\,.
\end{equation}
with $\eps = 2^{-4}$. The medium is plotted in Figure~\ref{fig:semilinear_medium}.

The reference solution and the local PDE solves are computed using the standard finite-volume scheme with uniform grid with mesh size $\Delta x = 2^{-8} = \tfrac{1}{256}$ and Newton's method is used to solve the resulting algebraic problem.
For our domain decomposition approach, we set $M_1=M_2 = 4$ to define the patches $\Omega_m$,  with boundary margins $\Delta x_{\rmo} = 2^{-4} = \tfrac{1}{16}$ to form $\Omega_m$.
The input and output dimensions of $Q^\eps_m$ are thus $(d_m,p_m) = (388, 388)$.

To obtain the training data, each patch $\Omega_m$ is further enlarged to a buffered patch $\ol{\Omega}_m$ by adding a margin of $\Delta x_\rmb = 2^{-4} = \tfrac{1}{16}$ to $\Omega_m$. On each patch $\ol{\Omega}_m$, $10,000$ samples are generated with random boundary conditions defined by $R_m = 1000$ and $D = 3$.
To train the NN, we use the loss function~\eqref{eqn:loss} with
\begin{equation}\label{eqn:loss_detail}
\begin{aligned}
\ell(\mathcal{Q}_m^{\rmNN}(\theta_m)\phi_{m}-\psi_{m})=& \|\mathcal{Q}_m^{\rmNN}(\theta_m)\phi_{m}-\psi_{m}\|^2 \\
& + \mu \sum_{l\in\mathscr{N}(m)}\|D_h\mathcal{Q}_{l,m}^{\rmNN}(\theta_m)\phi_{m}-D_h\psi_{l,m}\|^2\,,
\end{aligned}
\end{equation}
where $\mathcal{Q}_{l,m}^{\rmNN}(\theta_m)$ and $\mathcal{Q}_m^{\rmNN}(\theta_m)$ are NN approximation of $\mathcal{Q}_{l,m}(\theta_m)$ and $\mathcal{Q}_{l,m}(\theta_m)$, as defined in~\eqref{eqn:neighbors} and~\eqref{eqn:def_Q_m}, respectively; and $D_h$ is the discrete version of the derivative operator with step size $h$. The second term measures mismatch in the derivative so as to enforce the regularity.

To initialize the neural networks, we take $\mathcal{Q}_m^\rmL$ to be the boundary-to-boundary operator of the following linear elliptic equation
\begin{equation}\label{eqn:linear_initial_semi}
-\nabla \cdot (\kappa^\eps(x) \nabla u^\eps(x)) = 0,\quad x\in\Omega\,.
\end{equation}
We truncate the rank representation of $\mathcal{Q}_m^\rmL$ at rank $r_m=40$ to preserve all singular values bigger than a tolerance $\delta_1=10^{-2}$ so that the width of the hidden layer is $h_m = 80$.

\subsubsection{Offline training}
\begin{figure}[t!]
  \centering
  \includegraphics[width=0.6\textwidth]{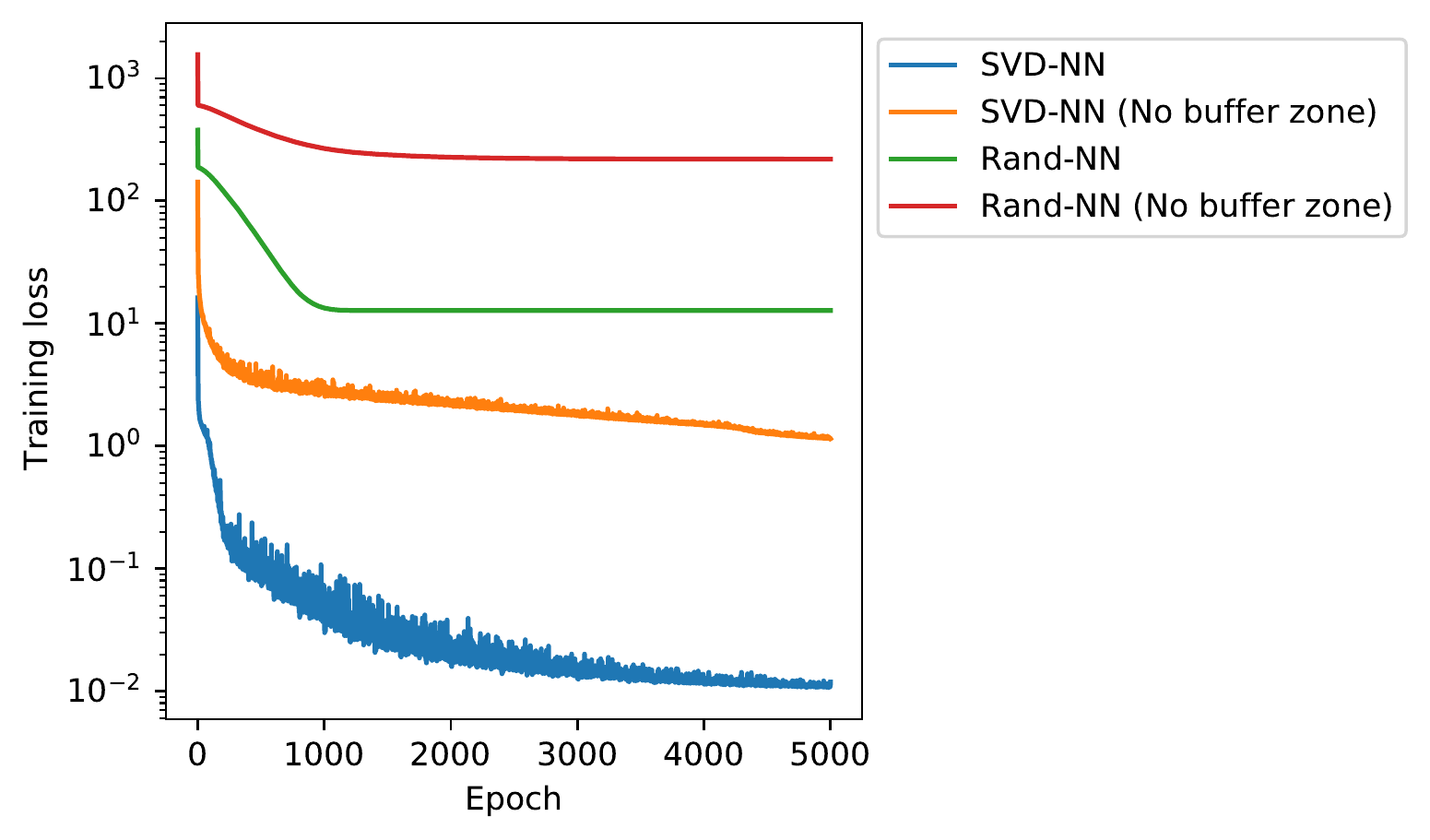}
  \caption{Training loss for loss function $\Lc$ \eqref{eqn:loss_detail} for patch (2,2). For the variants that use random initializations, we use the  PyTorch default, which generate the weights and biases in each layer uniformly from $(-\sqrt{d_\rminput},\sqrt{d_\rminput})$, where  $d_\rminput$ is the input dimension of the layer.}
  \label{fig:Lc_learning_curve}
\end{figure}

We show the improvements in the offline process for training $\mathcal{Q}_m^\rmNN$ due to the two strategies described in Subsection~\ref{sec:implement}: the use of enlarged patches, and initialization using SVD of a matrix representation of a linearized equation. Figure~\ref{fig:Lc_learning_curve} plots number of epochs in the offline training vs the training loss function $\Lc$ \eqref{eqn:loss_detail} associated with $\mathcal{Q}_{m}^{\rmNN}$ for the patch $m=(2,2)$ in four different settings: SVD-initialization on training data with buffer zone, SVD-initialization on training data without buffer zone, and the counterpart without SVD-initialization.
The same NN model is used in all four settings.
It is immediate that the training process has a much faster decay in error if buffer zone is adopted, and that the SVD initialization gives a much smaller error than random initialization.

\begin{figure}[t!]
  \centering
  \includegraphics[width=0.45\textwidth]{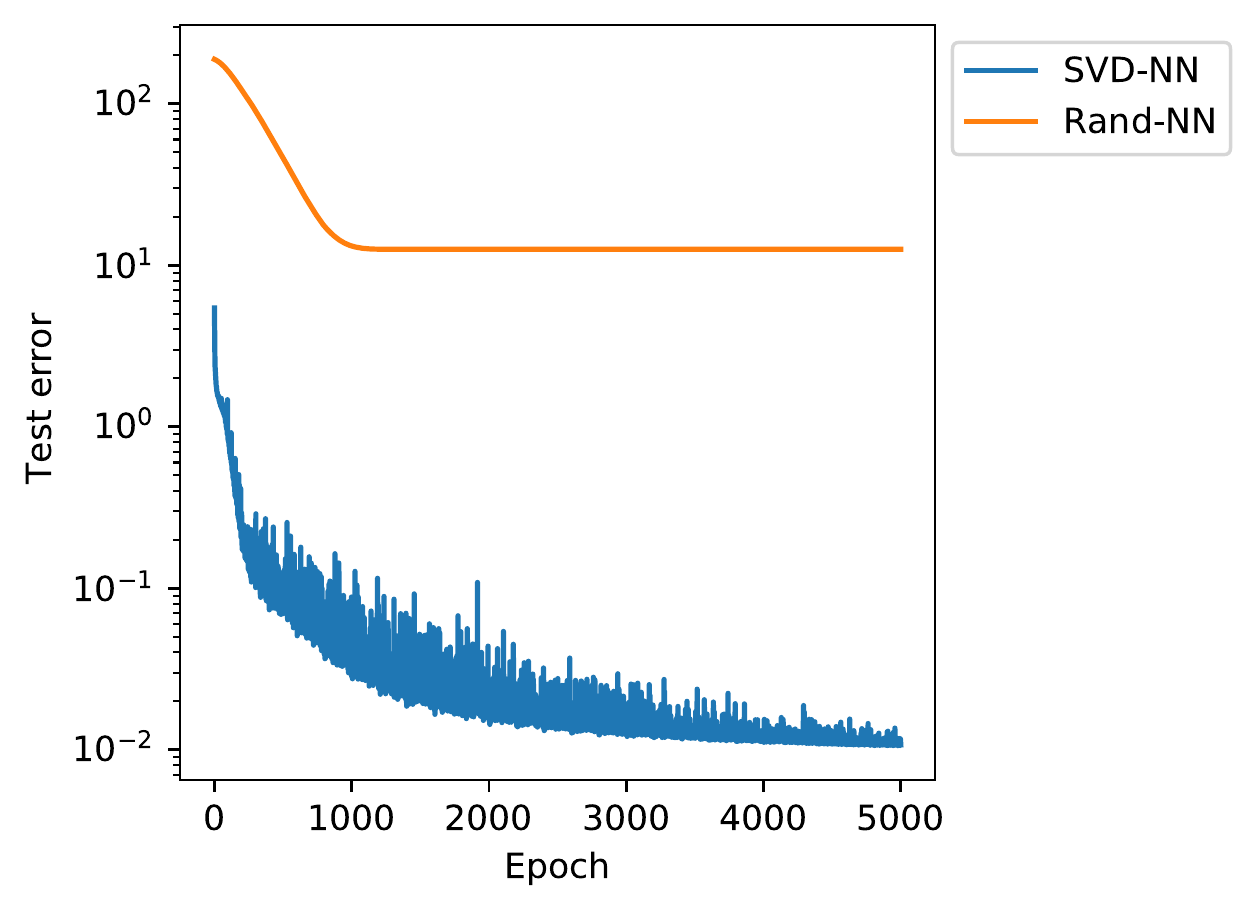}
  \caption{Testing error during the training for patch (2,2).}
  \label{fig:Lc_learning_curve_test}
\end{figure}

To show the generalization performance of the resulting trained NN,
we generate a test data set from the same distribution as the buffered training data set with 1,000 samples, for the same patch $m=(2,2)$.
Since the NNs trained using non-buffered data produce larger error, we only test the NNs trained with buffered data.
The test errors~\eqref{eqn:loss_detail} in the training process for different models are plotted in Figure~\ref{fig:Lc_learning_curve_test}.
Again, the use of buffered data along with SVD-initialization yields the best performance.

\begin{figure}[t!]
  \centering
  \includegraphics[width=0.35\textwidth]{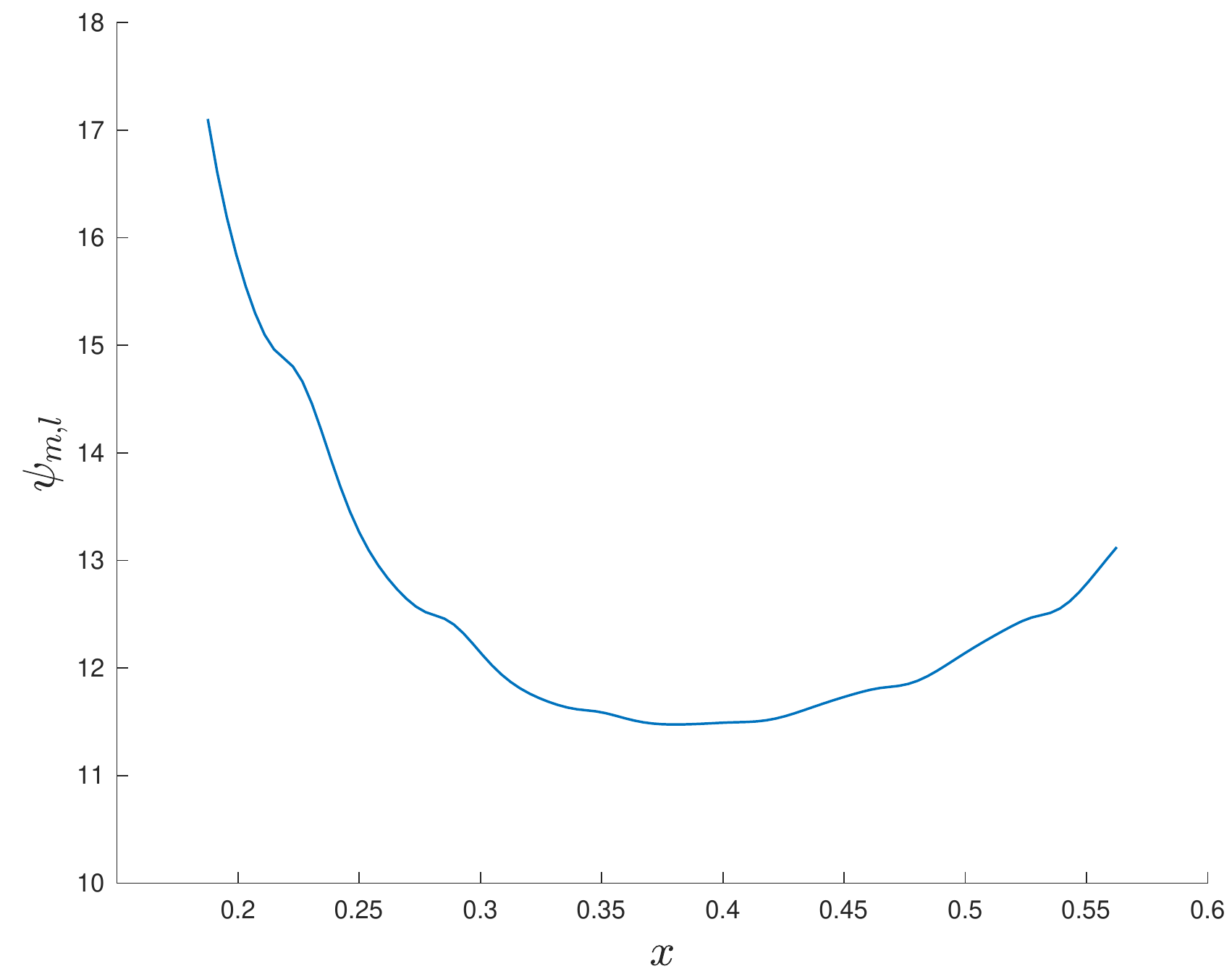}
  \includegraphics[width=0.35\textwidth]{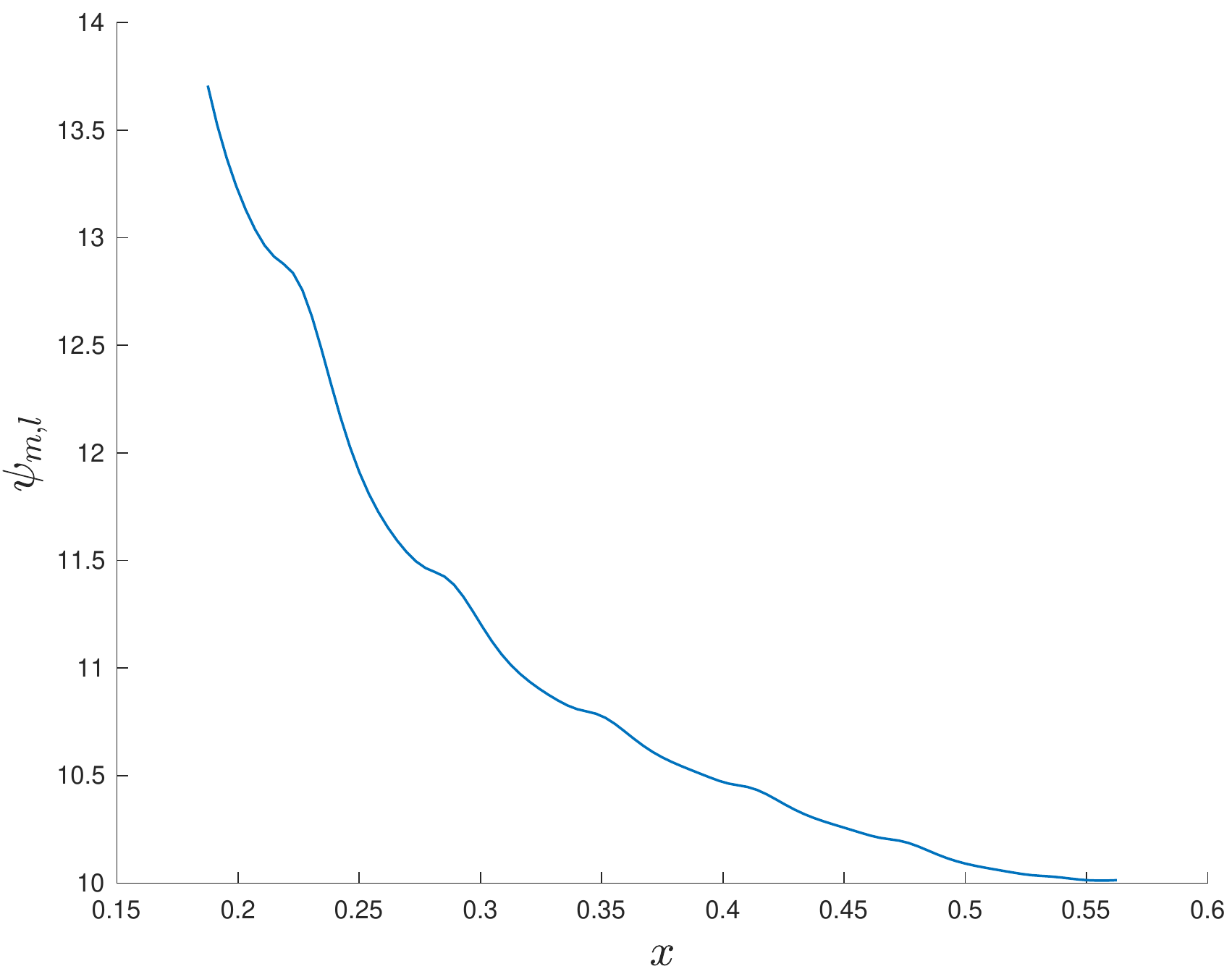}
  \includegraphics[width=0.35\textwidth]{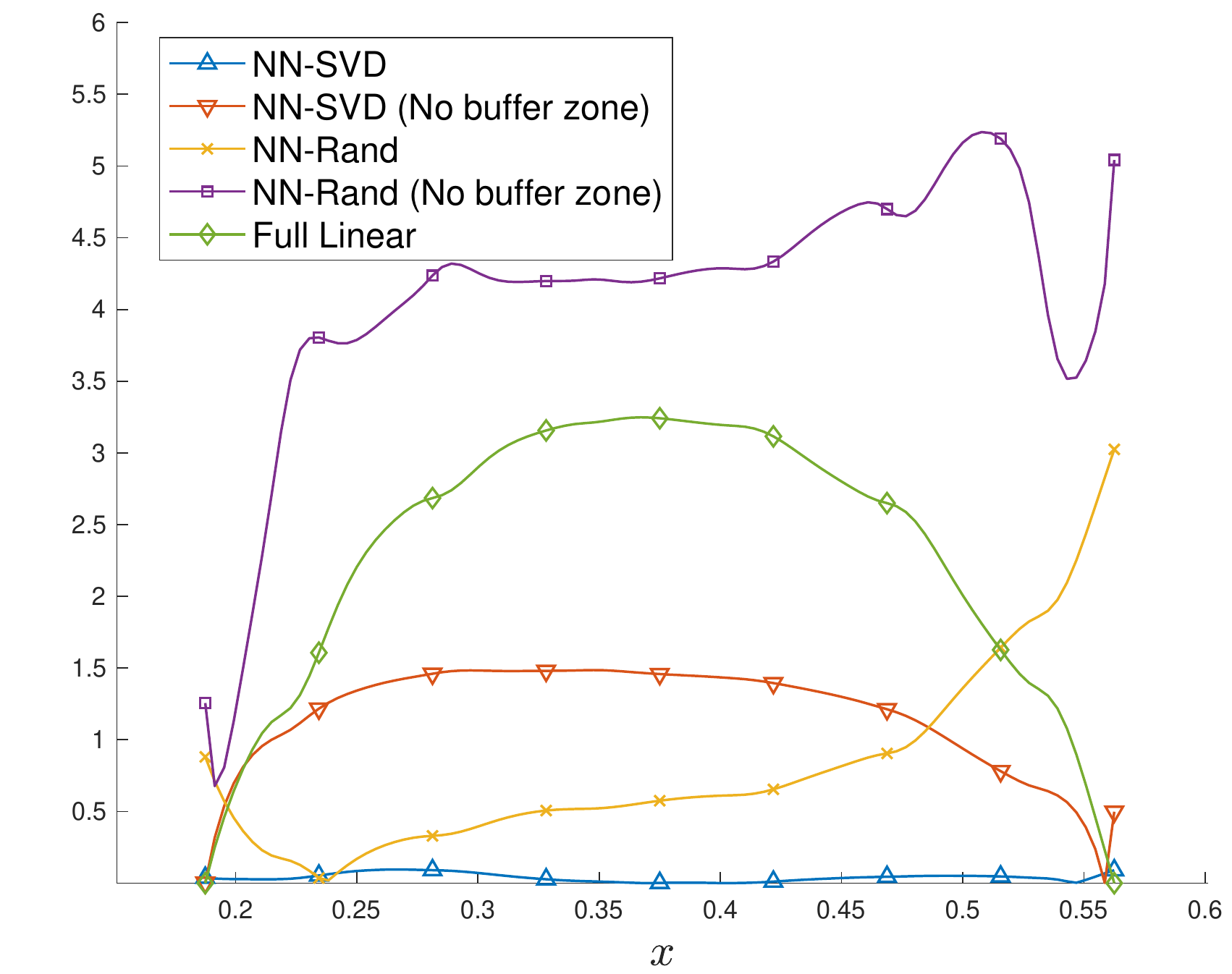}
  \includegraphics[width=0.35\textwidth]{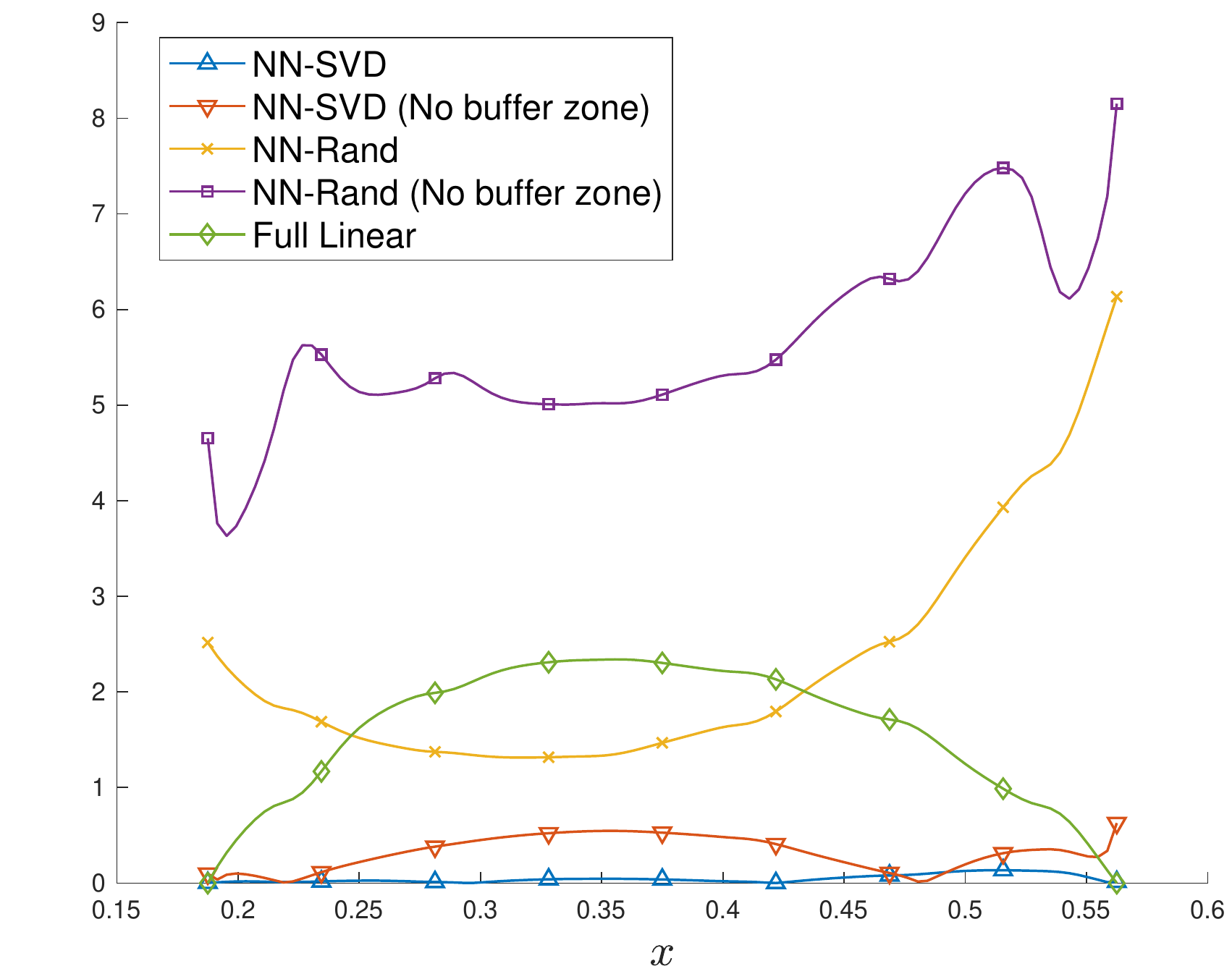}
  \caption{The top row shows the ground truths $\psi_{l,m}$ ($m = (2,2)$, $l = (2,1)$) of two samples in the test set. The bottom row shows the error $|\psi_{l,m}-\widetilde{\psi}_{l,m}|$, where $\widetilde{\psi}_{l,m}$ are computed by the low-rank SVD initialized $\mathcal{Q}_{m}^{\rmNN}$ (with and without buffer-zone), randomly initialized $\mathcal{Q}_{m}^{\rmNN}$ (with and without buffer-zone), and the linear operator $\mathcal{Q}_{m}^{\rmL}$.}
  \label{fig:psi_examples}
\end{figure}

To demonstrate generalization performance, we plot the predicted outputs for two typical examples in the test set in Figure~\ref{fig:psi_examples}. For comparison, we also plot the outputs produced by randomly initialized neural network and the linear operator $\mathcal{Q}_m^{\rmL}$. It can be seen that the low-rank SVD-initialized neural network has the best performance among all the initialization methods.

\begin{figure}[b!]
  \centering
  \includegraphics[width=0.35\textwidth]{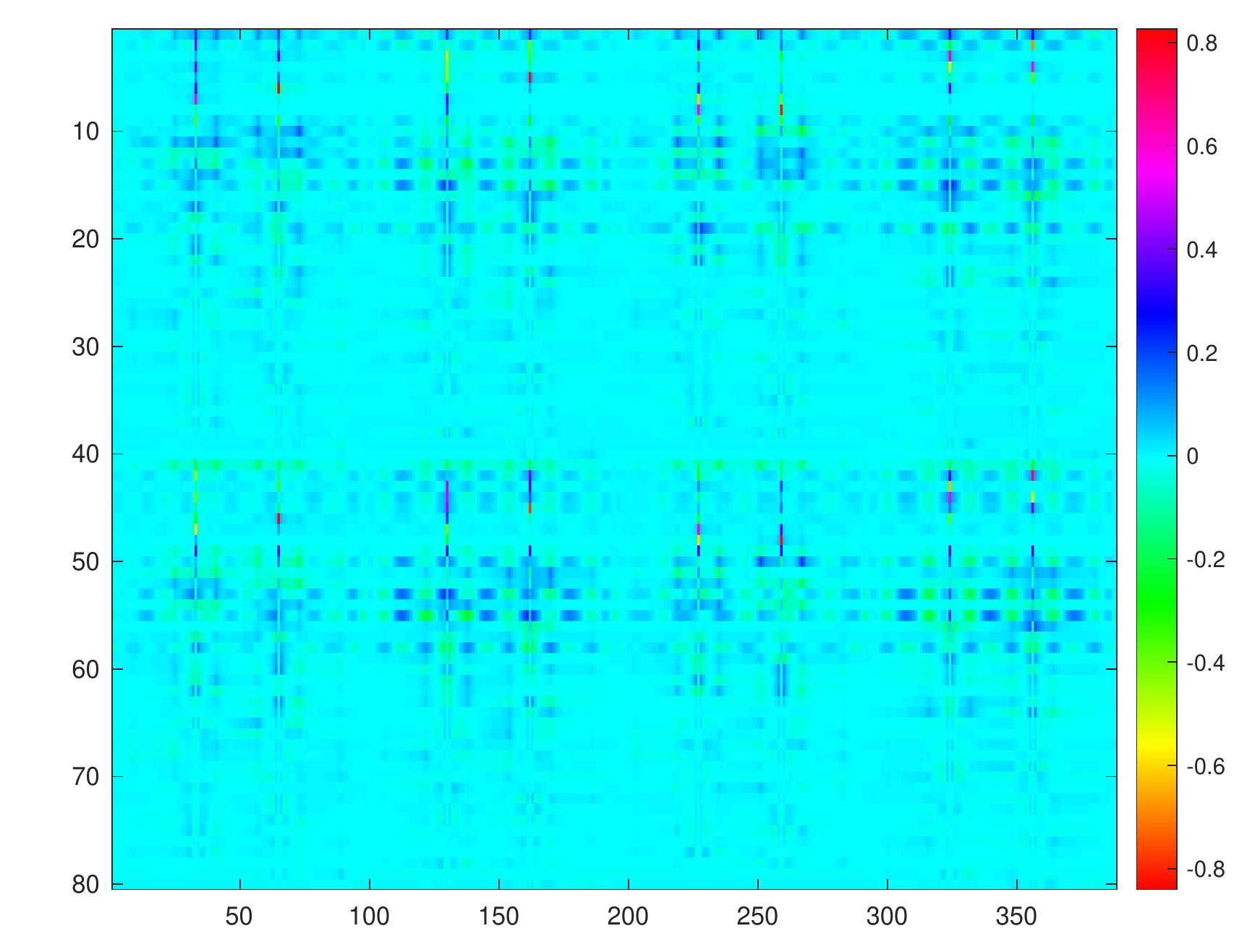}
  \includegraphics[width=0.35\textwidth]{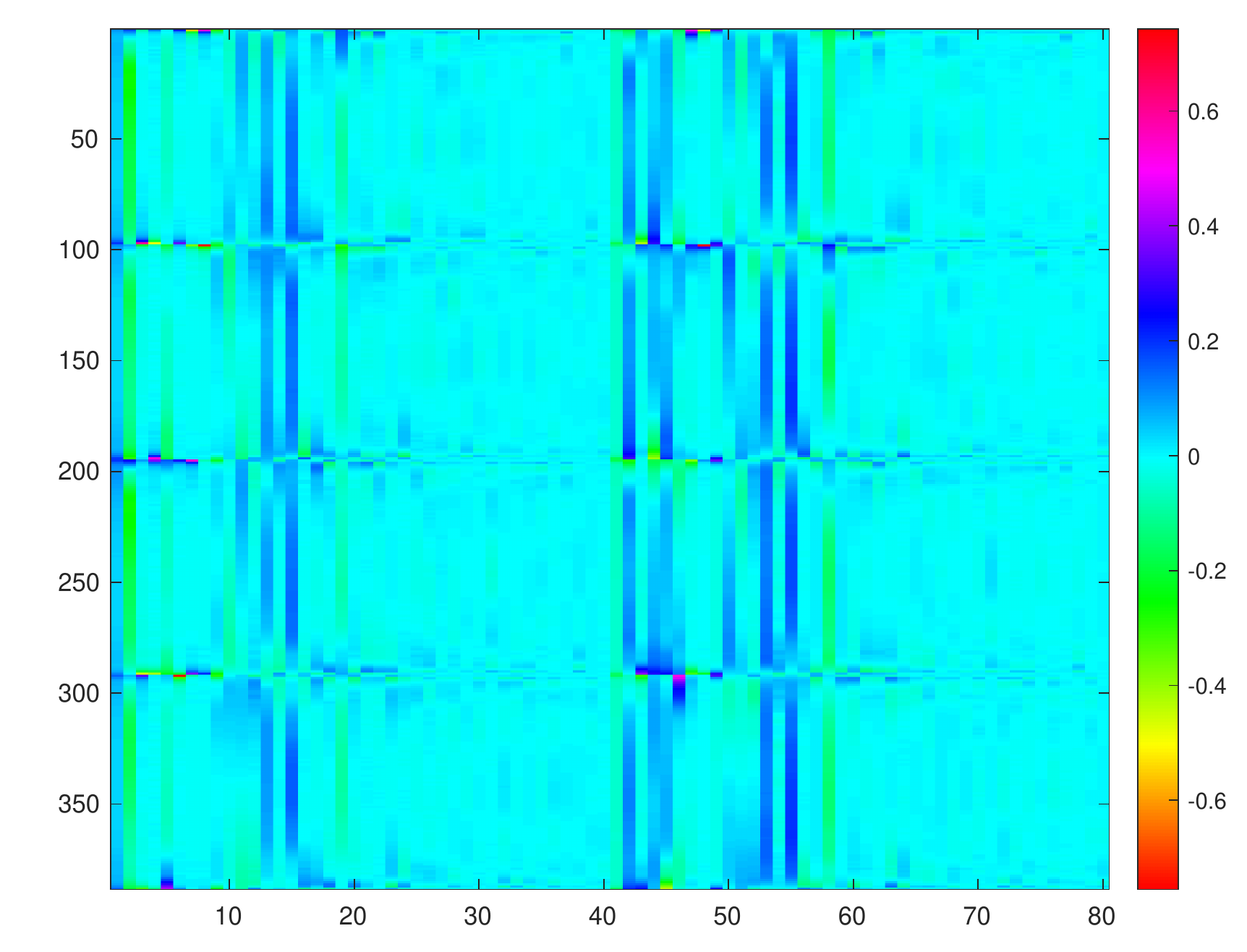}
  \\
  \includegraphics[width=0.35\textwidth]{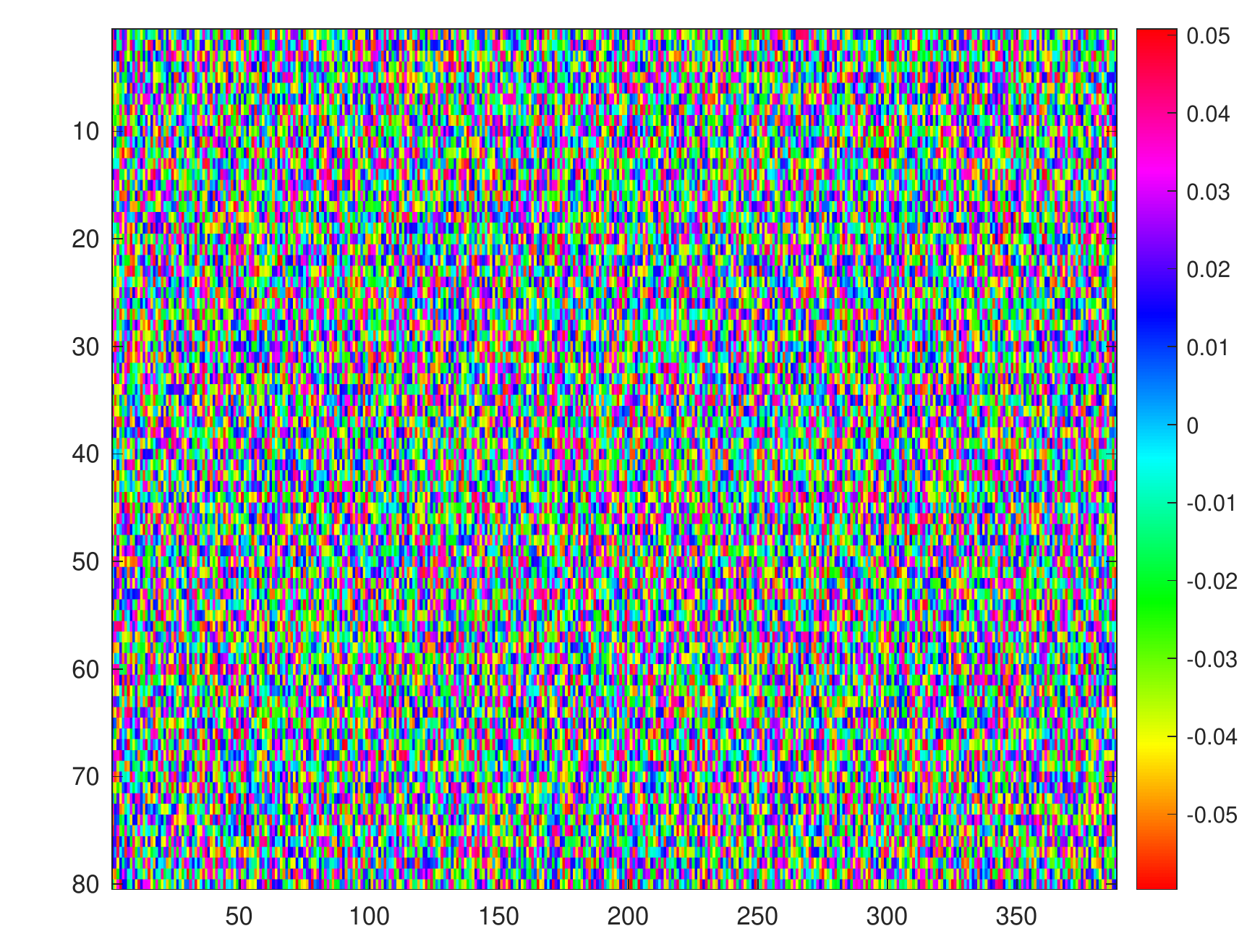}
  \includegraphics[width=0.35\textwidth]{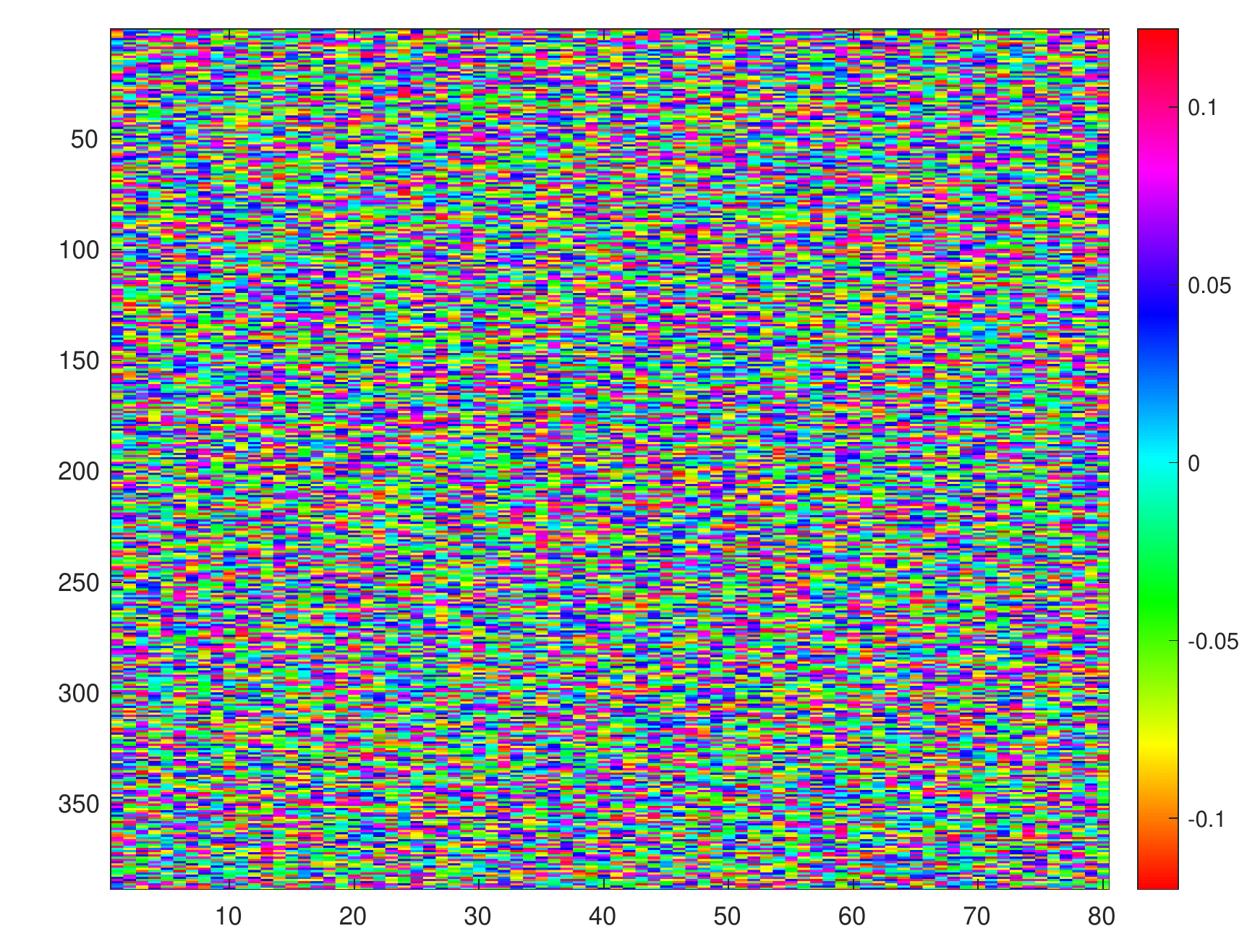}
  \caption{The first row shows the final weight matrices $W_1$ (left), $W_2$ (right) obtained the for SVD-initialized model on patch $m = (2,2)$. The second row shows the final weight matrices $W_1$ (left), $W_2$ (right) for randomly initialized model on patch $m = (2,2)$. In both cases, training  data is obtained by enlarging the patch.}
  \label{fig:parameter_semi}
\end{figure}

We note too that the neural network models initialized by the SVD of linear PDEs tend to be more interpretable.
Figure~\ref{fig:parameter_semi} shows the final weight matrices for models initialized by different methods. It can be seen that SVD-initialized model yields weight matrices with recognizable structure: the parameters for higher modes are near zero, and only the top $~25$ modes in the positive and negative halves are nontrivial.
By comparison, the trained weight matrices using randomly initialized parameters do not show any pattern or structure.

\subsubsection{Online phase: Schwarz iteration}

\begin{figure}
  \centering
  \includegraphics[width=0.3\textwidth]{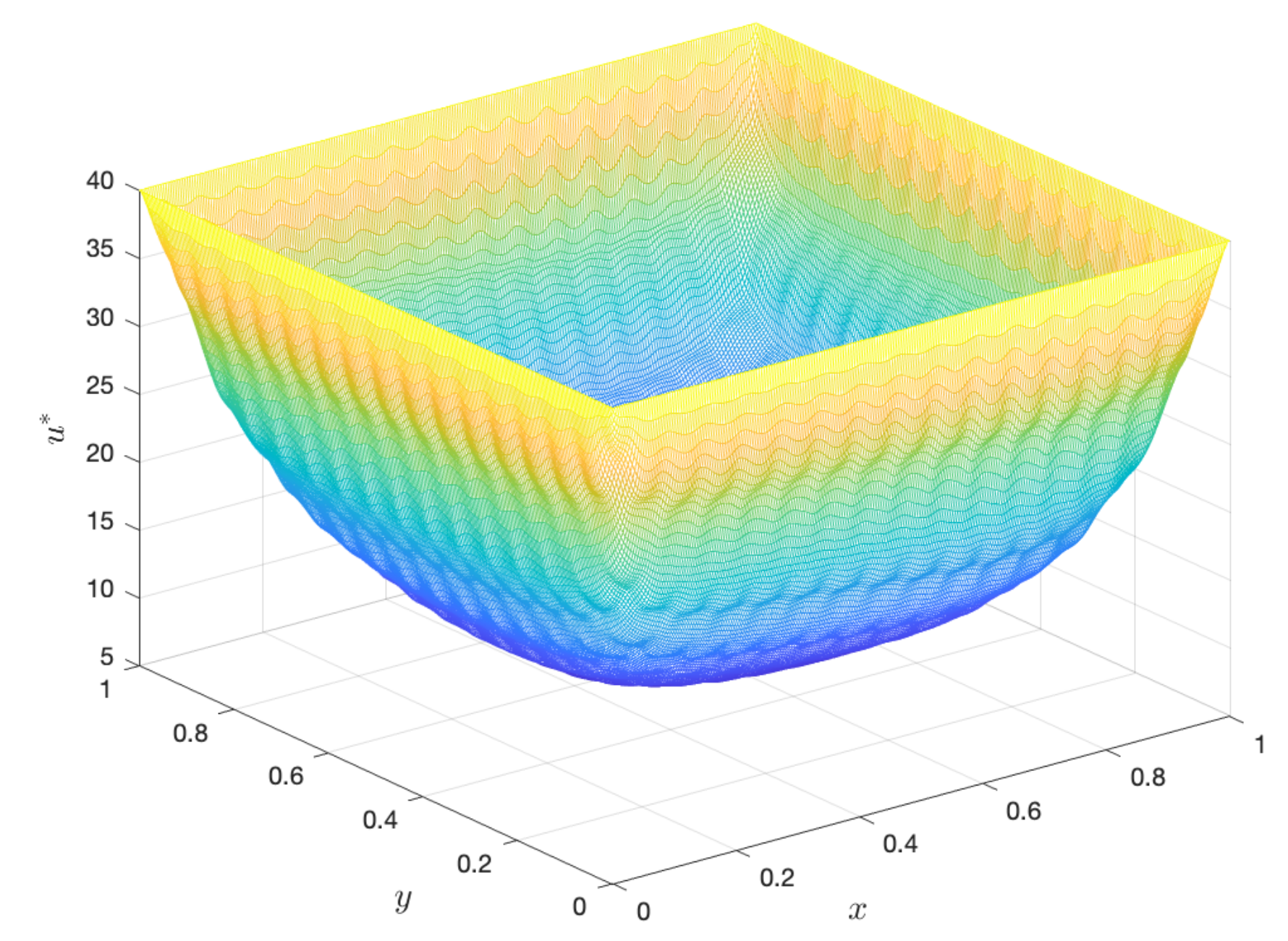}
  \includegraphics[width=0.3\textwidth]{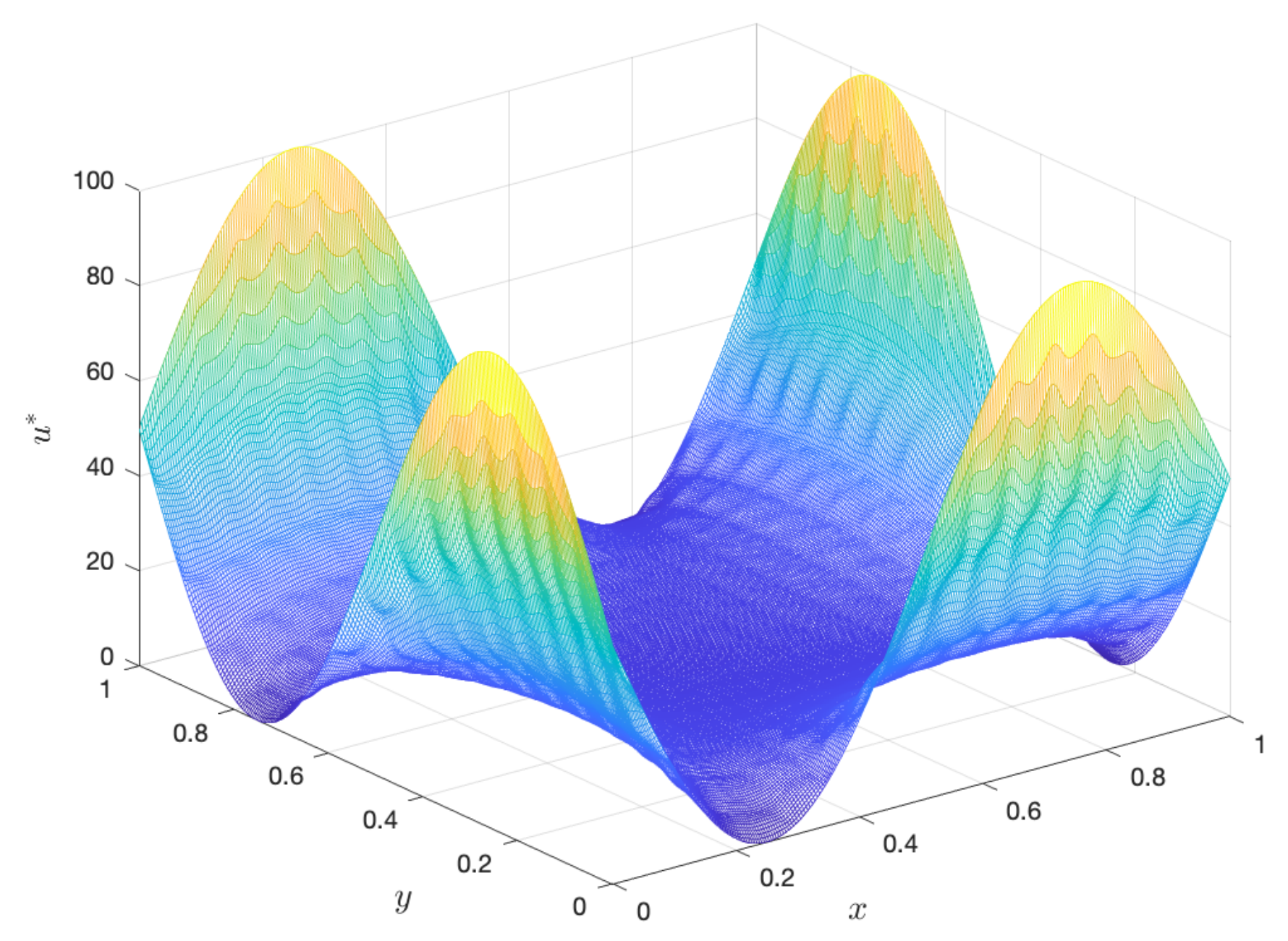}
  \includegraphics[width=0.3\textwidth]{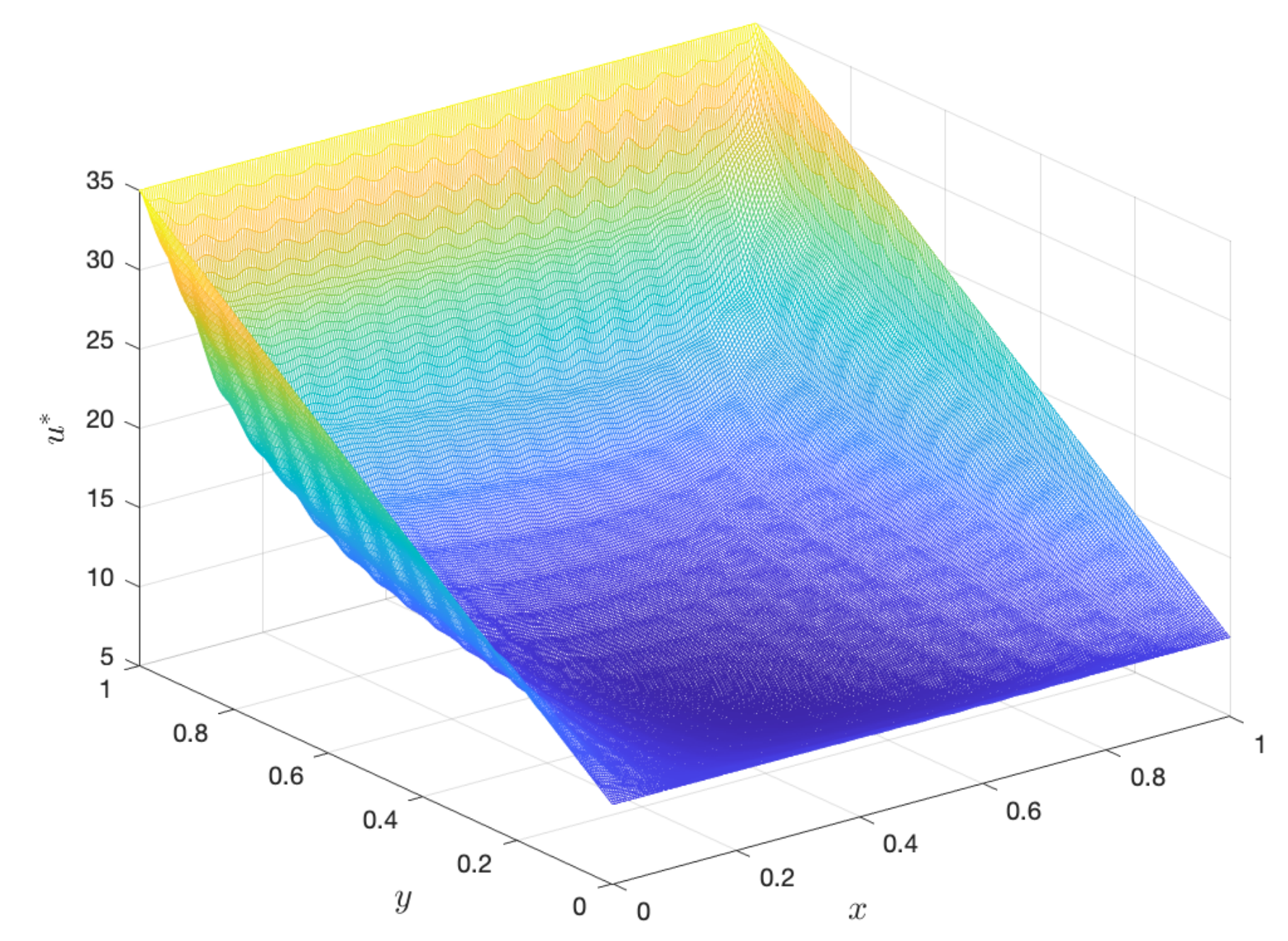}
  \\
  \includegraphics[width=0.3\textwidth]{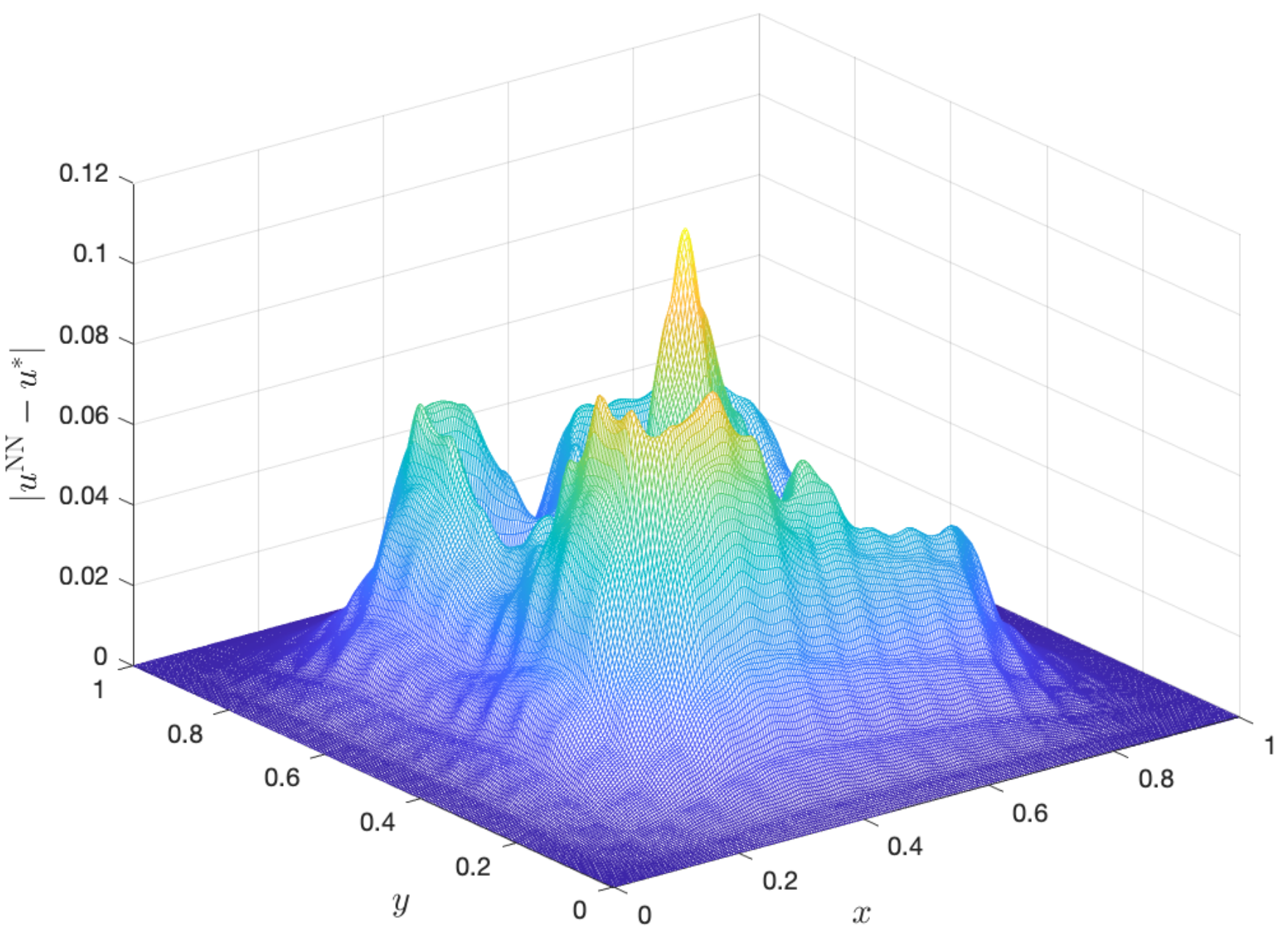}
  \includegraphics[width=0.3\textwidth]{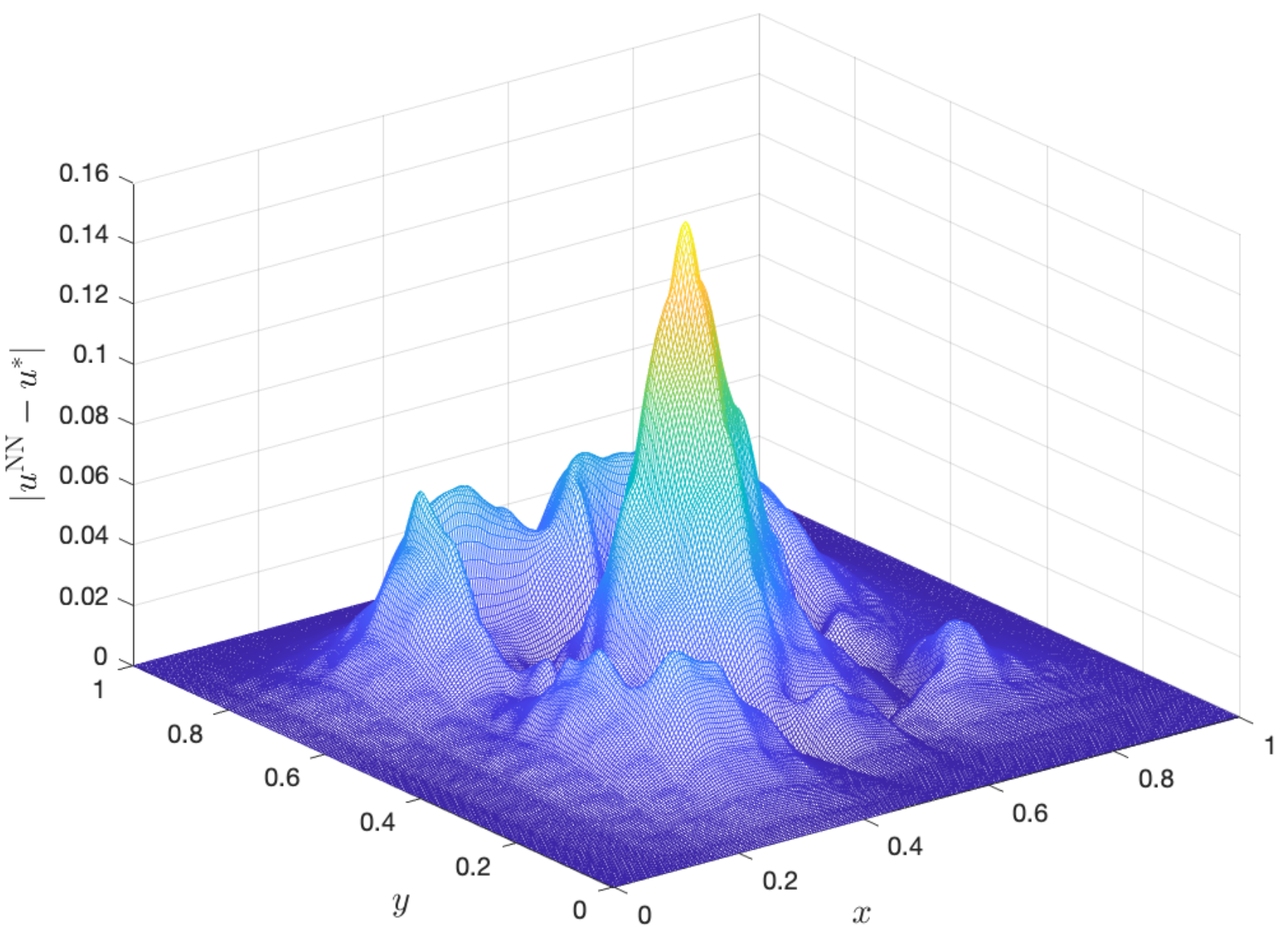}
  \includegraphics[width=0.3\textwidth]{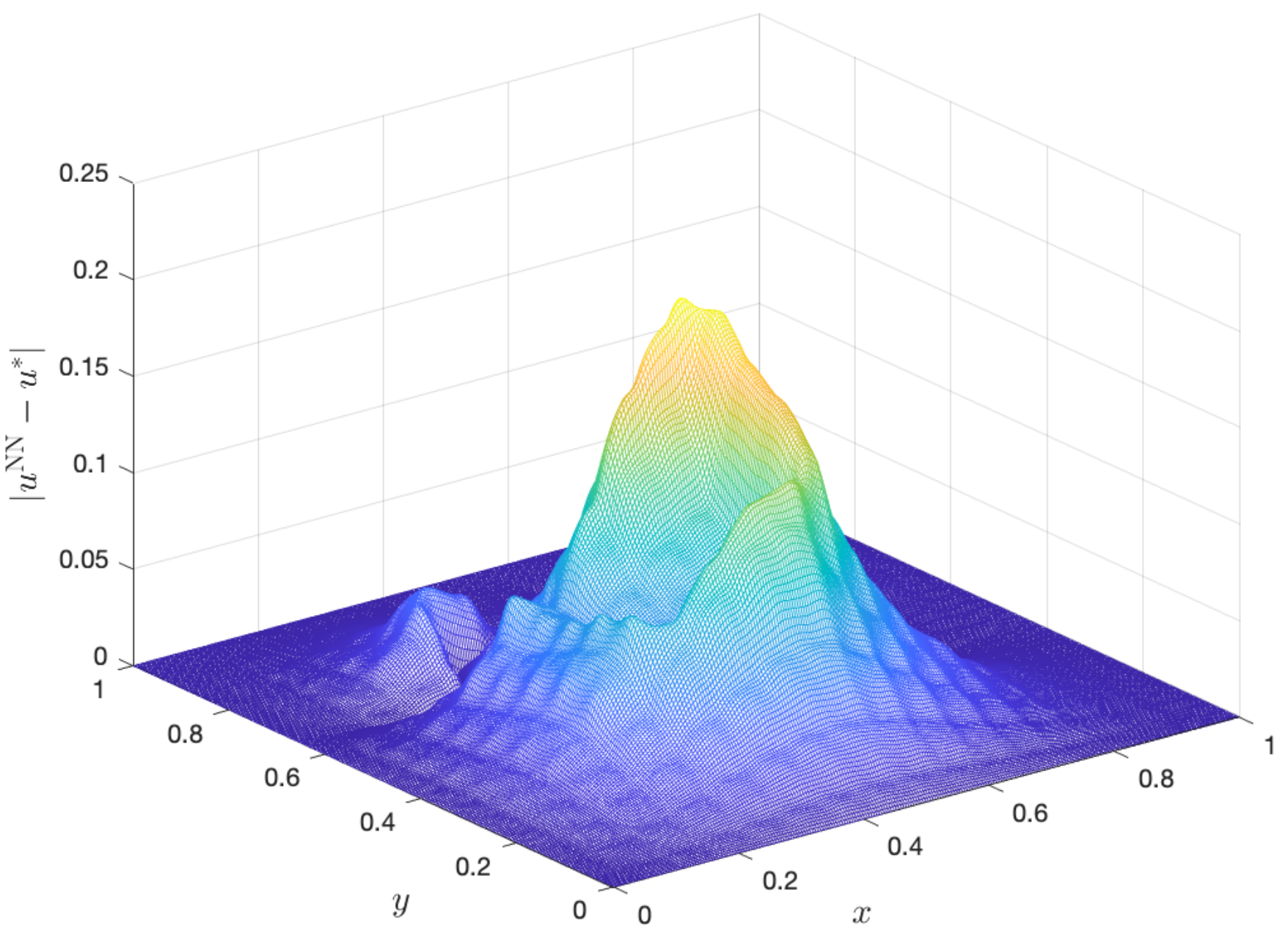}
  \caption{The first row shows the ground truth solutions $u^\ast$ for boundary conditions 1 to 3 from left to right. The second row shows the absolute error $|u^\rmNN-u^\ast|$ for boundary conditions 1 to 3 from left to right. (Note the much smaller vertical scale used in the second row.)}
  \label{fig:u_err_u1}
\end{figure}

We show results obtained by using the NN approximation $\mathcal{Q}_m^{\rmNN}(\theta_m)$ of the boundary-to-boundary map inside the Schwarz iteration.
Table~\ref{tbl:BCs} shows the boundary conditions used for the three problems we tested.
(The same medium \eqref{eqn:medium} is used in all cases.)
We use $\delta_0=10^{-4}$ for the tolerance in Algorithm~\ref{alg:NN_online}, and use the full accuracy local solvers as in the generation of training data set. In Figure~\ref{fig:u_err_u1}, we plot the ground truth solutions for different boundary conditions and the absolute error of $u^\rmNN$ obtained by neural network-based Schwarz iteration.
(Note that the scaling of the $y$-axis in the latter is different from the former.)
The relative errors obtained for the four variants of NN approximation along with the linear approximation $\mathcal{Q}_m^\rmL$ to the boundary-to-boundary map can be found in Tables~\ref{tbl:Err1} and \ref{tbl:Err3}.
Note that the smallest errors are attained by the variant that uses the SVD initialization and buffered patches.
To demonstrate the efficiency of our method, we compare the CPU time of neural network based-Schwarz method and the classical Schwarz method, using the same tolerance $\delta_0=10^{-4}$ for the latter.
The NNs we used for the test is trained by SVD initialization, and its training data is generated with buffer zone.
Since NN-produced local boundary-to-boundary map is only an approximation to the ground truth, for a fair comparison, we also run the reference local solution with a relaxed accuracy requirement.
The CPU time, number of iteration and error comparison can be found in Table~\ref{tbl:Runtime}.
In all three test cases, the NN approximate executes faster than the conventional local solution technique as a means of implementing the boundary-to-boundary map, while producing $H^1$ errors of the same order.
\begin{table}[t!]
	\centering
	\begin{tabular}{ c | c }
		\hline \hline
			No.	& Boundary condition \\
		\hline
		1  	& \tabincell {c}{$\phi(x,0) = 40$\,, $\phi(x,1) = 40$ \\ $\phi(0,y) = 40$\,, $\phi(1,y) = 40$} \\
        \hline
		2	& \tabincell {c}{$\phi(x,0) = 50-50\sin(2\pi x)$\,, $\phi(x,1) = 50+50\sin(2\pi x)$ \\ $\phi(0,y) = 50+50\sin(2\pi y)$\,, $\phi(1,y) = 50-50\sin(2\pi y)$} \\
        \hline
		3	& \tabincell {c}{$\phi(x,0) = 10$\,, $\phi(x,1) = 35$ \\ $\phi(0,y) = 10+25y$\,, $\phi(1,y) = 10+25y$} \\
		\hline\hline
	\end{tabular}
    \caption{Boundary conditions used in the global test.}
  \label{tbl:BCs}
\end{table}

\begin{table}[b!]
	\centering
	\begin{tabular}{ c |ccc | ccc }
		\hline \hline
			Problem Number	& \multicolumn{3}{c|}{1} & \multicolumn{3}{c}{2}\\
		\hline
		Relative Error  & $L^2$ & $H^1$ & $L^\infty$ & $L^2$ & $H^1$ & $L^\infty$\\
        \hline
		SVD-NN & \textbf{0.0013} & \textbf{0.0028} & \textbf{0.0029} & \textbf{0.0010} & \textbf{0.0010} & \textbf{0.0016} \\
        SVD-NN (No buffer zone) & 0.0042 & 0.0091 & 0.0078 & 0.0029 & 0.0030 & 0.0039 \\
		Rand-NN & 0.0425 & 0.0445 & 0.0907 & 0.0379 & 0.0188 & 0.0370 \\
        Rand-NN (No buffer zone) & 0.0882 & 0.0965 & 0.1555 & 0.0773 & 0.0400 & 0.0629 \\
		Linear & 0.0606 & 0.0644 & 0.1066 & 0.0505 & 0.0252 & 0.0415 \\
		\hline\hline
	\end{tabular}
    \caption{Relative error for global solutions by different methods.}
  \label{tbl:Err1}
\end{table}

\begin{table}[t!]
	\centering
	\begin{tabular}{ c | ccc }
		\hline \hline
			Problem Number &\multicolumn{3}{c}{3} \\
		\hline
		Relative Error &  $L^2$ & $H^1$ & $L^\infty$ \\
        \hline
		SVD-NN &  \textbf{0.0035} & \textbf{0.0059} & \textbf{0.0058} \\
        SVD-NN (No buffer zone) & 0.0235 & 0.0341 & 0.0346 \\
		Rand-NN &  0.1029 & 0.1293 & 0.1333 \\
        Rand-NN (No buffer zone) &  0.1739 & 0.2277 & 0.2078 \\
		Linear & 0.0614 & 0.0729 & 0.0776 \\
		\hline\hline
	\end{tabular}
    \caption{Relative error for global solutions by different methods. (Continued)}
  \label{tbl:Err3}
\end{table}

\begin{table}[t!]
	\centering
     \begin{tabular}{ c | cc | cc | cc}
		\hline \hline
		  Problem Number & \multicolumn{2}{c|}{1} & \multicolumn{2}{c|}{2} & \multicolumn{2}{c}{3}\\
        \hline
        Method &  NN & Classical &  NN & Classical &  NN & Classical\\
         \hline
        CPU time & \textbf{12.4} & 17.2 & \textbf{13.9} & 18.5 & \textbf{13.4} & 19.5\\
		Iteration & 30 & 29 & 30 & 29 & 34 & 35\\
		$H^1$ Error & 0.0035 & 0.0024 & 0.0012 & 0.0007 & 0.0062 & 0.0022\\
		\hline\hline
	\end{tabular}
    \caption{CPU time (s), number of iterations and the $H^1$ error of the classical Schwarz iteration and the neural network accelerated Schwarz iteration.}
  \label{tbl:Runtime}
\end{table}

\subsection{$p$-Laplace equations} \label{sec:numerics_pPoi}
\begin{figure}
  \centering
  \subfloat[$\kappa$]{
  \includegraphics[width=0.35\textwidth]{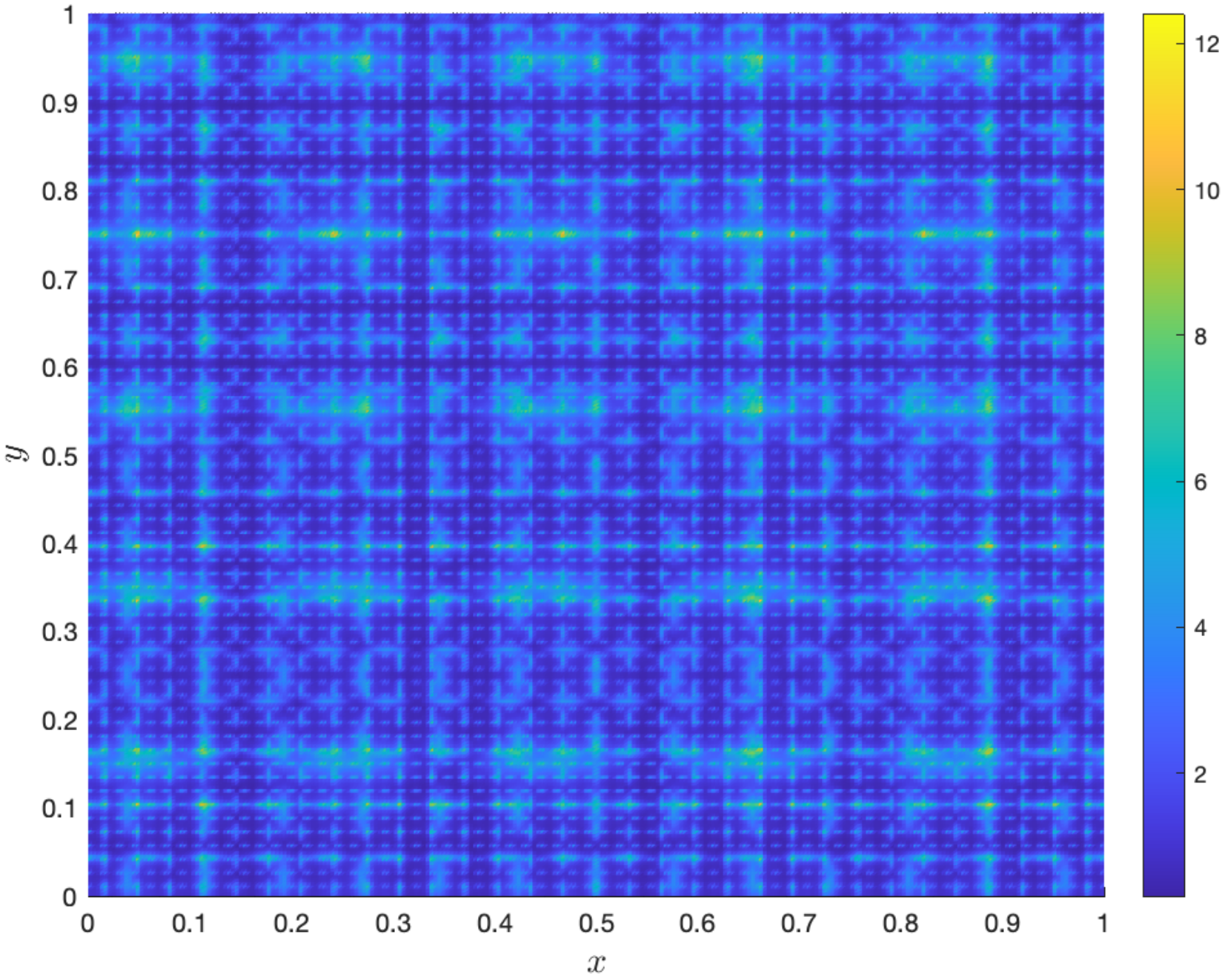}
  }
  \subfloat[$\log_{10}(\kappa|\nabla u^\eps|^{p-2})$]{
  \includegraphics[width=0.35\textwidth]{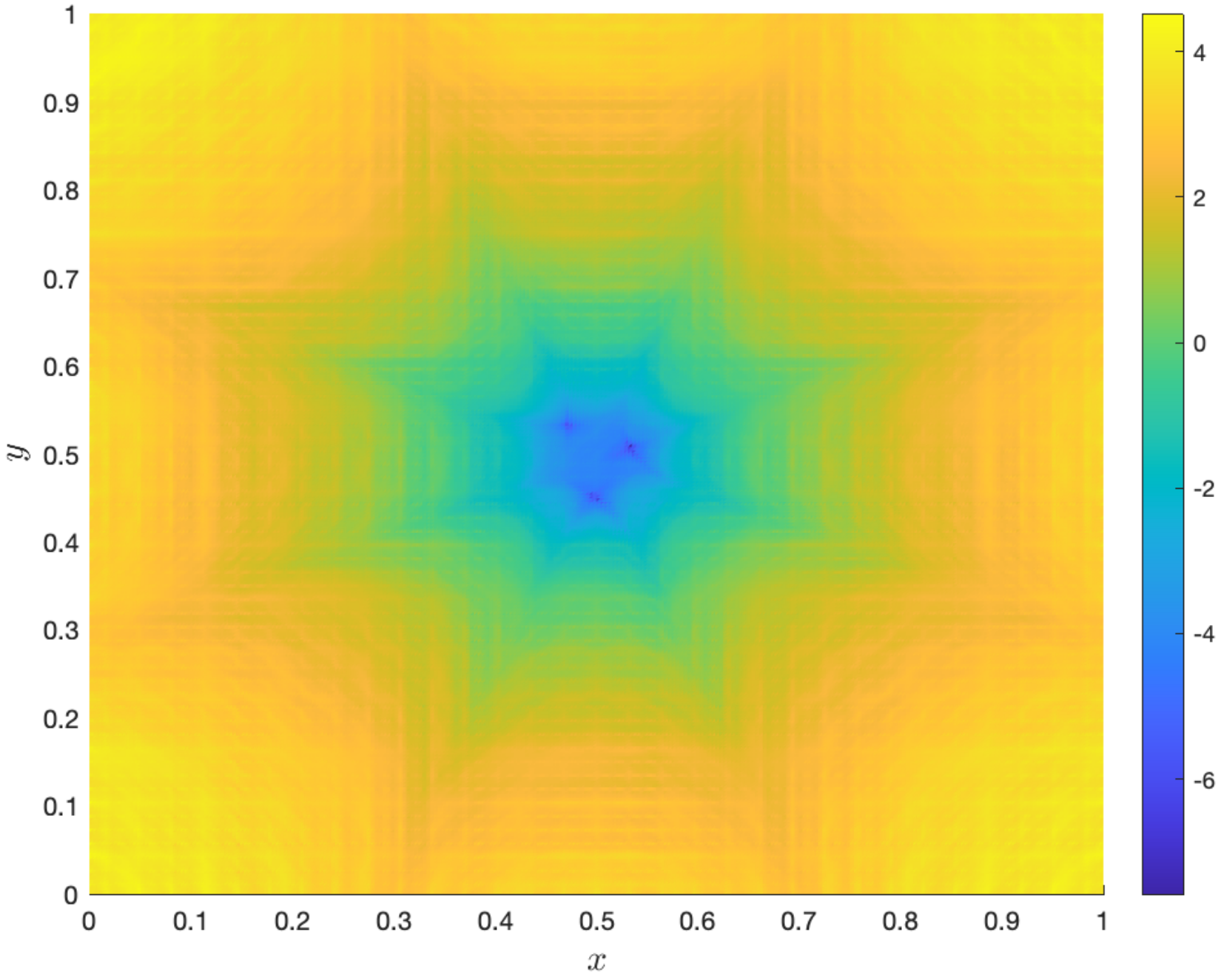}
  }
  \caption{Medium $\kappa$ and $\kappa|\nabla u^\eps|^{p-2}$ for $p$-Laplace equation. The solution $u^\eps$ is computed by boundary condition 1 (See Table \ref{tbl:BCs_pPoi}).}
  \label{fig:pPoi_medium}
\end{figure}

The second example concerns the multiscale $p$-Laplace elliptic equation~\cite{AsMa:1974principles,GeKoMaYv:2017homogenization,PrYe:2015reduced,ChEfShYe:2017multiscale,LiChZh:2021iterated} defined as follows:
\begin{equation}\label{eqn:pPoi_elliptic}
\begin{cases}
-\nabla \cdot (\kappa^\eps(x) |\nabla u^\eps|^{p-2} \nabla u^\eps(x)) = 0,\quad& x\in\Omega\,, \\
u^\eps(x) = \phi(x),\quad& x\in\partial\Omega\,,
\end{cases}
\end{equation}
where we use $p=6$ in this section, and the oscillatory medium is
\begin{equation}
\begin{aligned}
\kappa^\eps(x,y) = \frac{1}{6}
\bigg( &\frac{1.1+\sin(2\pi x/\eps_1)}{1.1+\sin(2\pi y/\eps_1)}
       + \frac{1.1+\sin(2\pi y/\eps_2)}{1.1+\cos(2\pi x/\eps_2)}
       + \frac{1.1+\cos(2\pi x/\eps_3)}{1.1+\sin(2\pi y/\eps_3)} \\
       &+ \frac{1.1+\sin(2\pi y/\eps_4)}{1.1+\cos(2\pi x/\eps_4)}
       + \frac{1.1+\cos(2\pi x/\eps_5)}{1.1+\sin(2\pi y/\eps_5)}
       + \sin(4x^2y^2) + 1 \bigg)\,,
\end{aligned}
\end{equation}
with $\eps_1 = 1/5, \eps_2 = 1/13, \eps_3 = 1/17, \eps_4 = 1/31, \eps_5 = 1/65$.
See Figure~\ref{fig:pPoi_medium} for an illustration of the medium.

\begin{figure}
  \centering
  \includegraphics[width=0.9\textwidth]{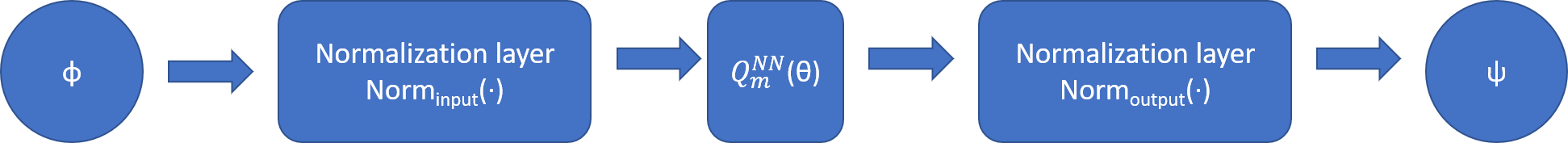}
  \caption{Neural network architecture for the boundary-to-boundary map in the p-Laplace equation.}
  \label{fig:NN_architecture}
\end{figure}

Noting that the differential equation in \eqref{eqn:pPoi_elliptic} is invariant when a constant is added or when multiplied by a constant, we use a normalization layer to improve the accuracy and robustness of the two-layer neural network.
Given a input boundary condition $\phi_m\in\Rb^{d_m}$, we define the input normalization layer by
\begin{equation}
\text{Norm}_{\text{input}}(\phi_m) = \frac{\phi_m - \widetilde{\phi}_m}{\max\{\|\phi_m\|_2,\epsilon_1\}}
\end{equation}
where $\widetilde{\phi}_m := 1/d_m \sum_{i=1}^{d_m} (\phi_m)_i$ is the mean, and the norm  is defined by $\|\phi_m\|_2^2 = \Delta x \sum_{i=1}^{d_m} (\phi_m)_i^2$.
We use $\epsilon_1 = 10^{-8}$ for the regularization constant.
The output normalization layer is defined by
\begin{equation}
\text{Norm}_{\text{output}}(\psi_m) = \max\{\|\phi_m\|_2,\epsilon_1\} \psi_m + \widetilde{\phi}_m \,.
\end{equation}
The overall architecture is illustrated in Figure~\ref{fig:NN_architecture}.

To compute both the reference solution and the patchwise solutions~\eqref{eqn:general_local_buffer}, we formulate the discretization using the standard piecewise linear finite element with uniform triangular grid, and solve with a preconditioned gradient descent method \cite{HuLiLi:2007:preconditioned}, where the line search parameter is computed by the Matlab function \underline{fminunc}.
The mesh size is  $\Delta x = 2^{-8} = {1}/{256}$.

For the domain decomposition, we set $M = 8$ with $\Delta x_{\rmo} = .03125$ to form $\Omega_m$.
The resulting input dimension is $d_m = 196$ and the output dimension is $p_m = 196$.
Training data is produced on enlarged patches $\ol{\Omega}_m$ with $\Delta x_{\text{b}} = .09375$.
On each patch, $1,000$ samples are generated with random distribution parameters $R_m = 10$ and $D = 3$.
To initialize the neural networks, we take $\mathcal{Q}_m^\rmL$ to be the boundary-to-boundary operator of the linear elliptic equation $-\nabla \cdot (\kappa^\eps(x) \nabla u(x)) = 0$.
We truncate the rank presentation of $\mathcal{Q}_m^\rmL$ at rank $r_m=36$, to preserve all singular values greater than a tolerance $\delta_1=10^{-2}$ so that the width of the hidden layer is $h_m = 72$.

\begin{figure}[b!]
  \centering
  \includegraphics[width=0.6\textwidth]{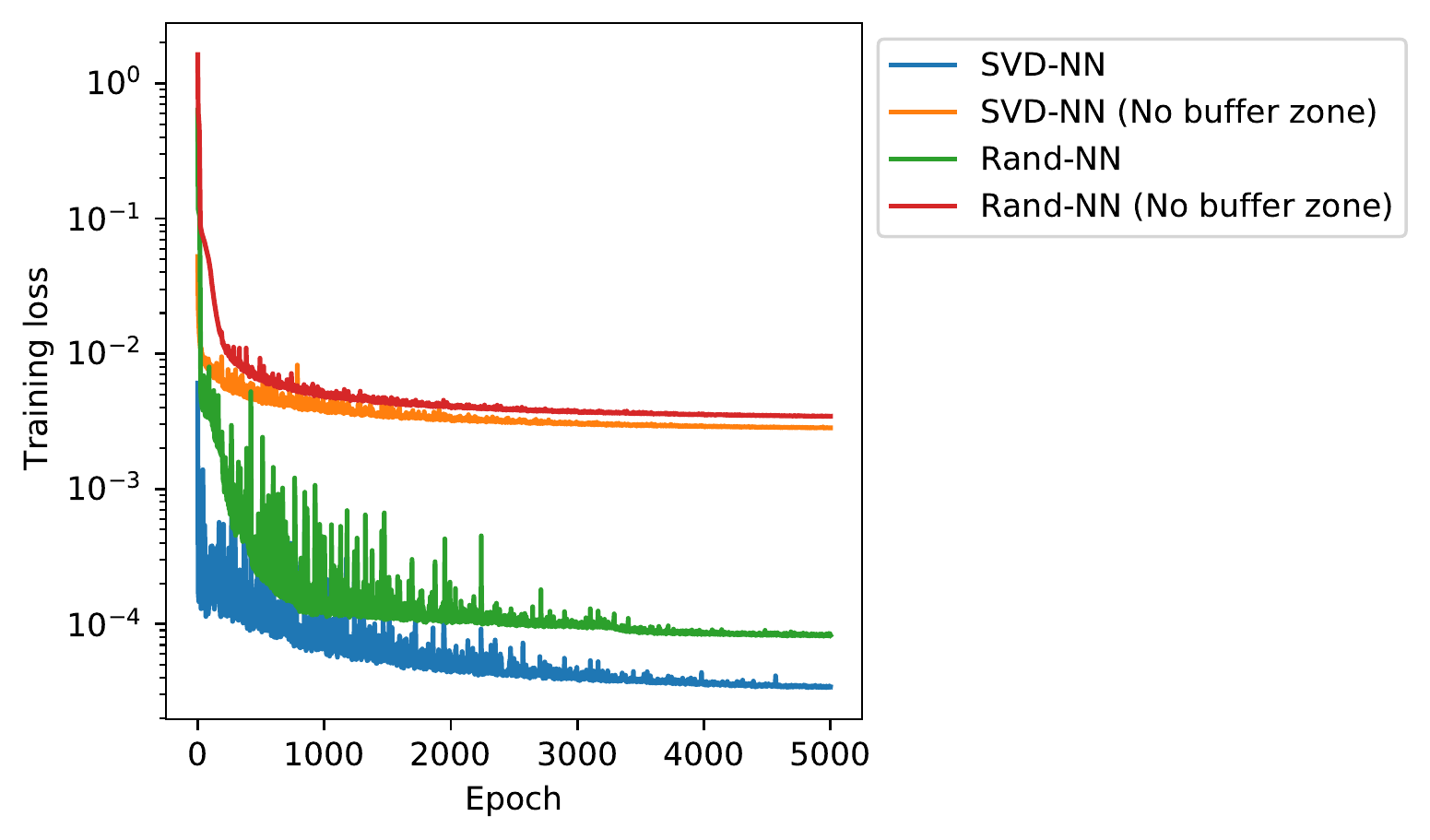}
  \caption{Training loss using loss function $\Lc$ \eqref{eqn:loss_detail} for patch $m=(2,2)$. We use the default random initialization method in PyTorch, which generate the weights and biases in each layer uniformly from $(-\sqrt{d_\rminput},\sqrt{d_\rminput})$ with $d_\rminput$ being the input dimension of the layer.}
  \label{fig:Lc_learning_curve_pPoi}
\end{figure}

\subsubsection{Offline training}

Here we show the improvements in the training process of $\mathcal{Q}_m^\rmNN$ by using the sampling and initialization strategies in Subsection~\ref{sec:implement}.
Figure~\ref{fig:Lc_learning_curve_pPoi} shows the training loss vs epochs for learning  $\mathcal{Q}_{m}^{\rmNN}$ for the patch $m=(2,2)$ using 1,000 samples.
The four variants are the same as in  Figure~\ref{fig:Lc_learning_curve}.
As for the previous example, the most effective training loss is for the variant in which samples are computed from buffered patches, using a reduced SVD initialization based on the linear approximate operator.
\begin{figure}[t!]
  \centering
  \includegraphics[width=0.45\textwidth]{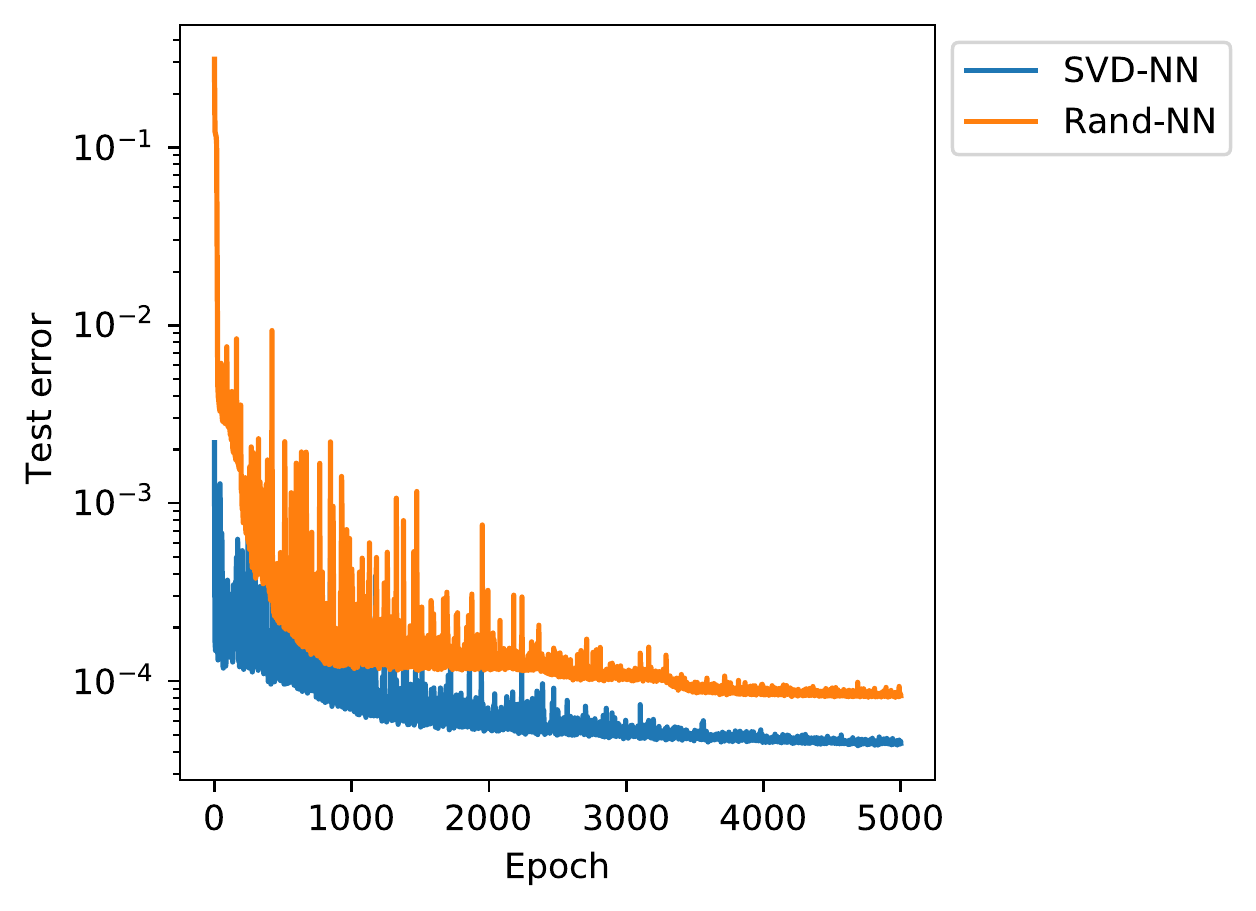}
  \caption{Testing error during the training for patch (2,2).}
  \label{fig:Lc_learning_curve_test_pPoi}
\end{figure}
\begin{figure}[b!]
  \centering
  \includegraphics[width=0.35\textwidth]{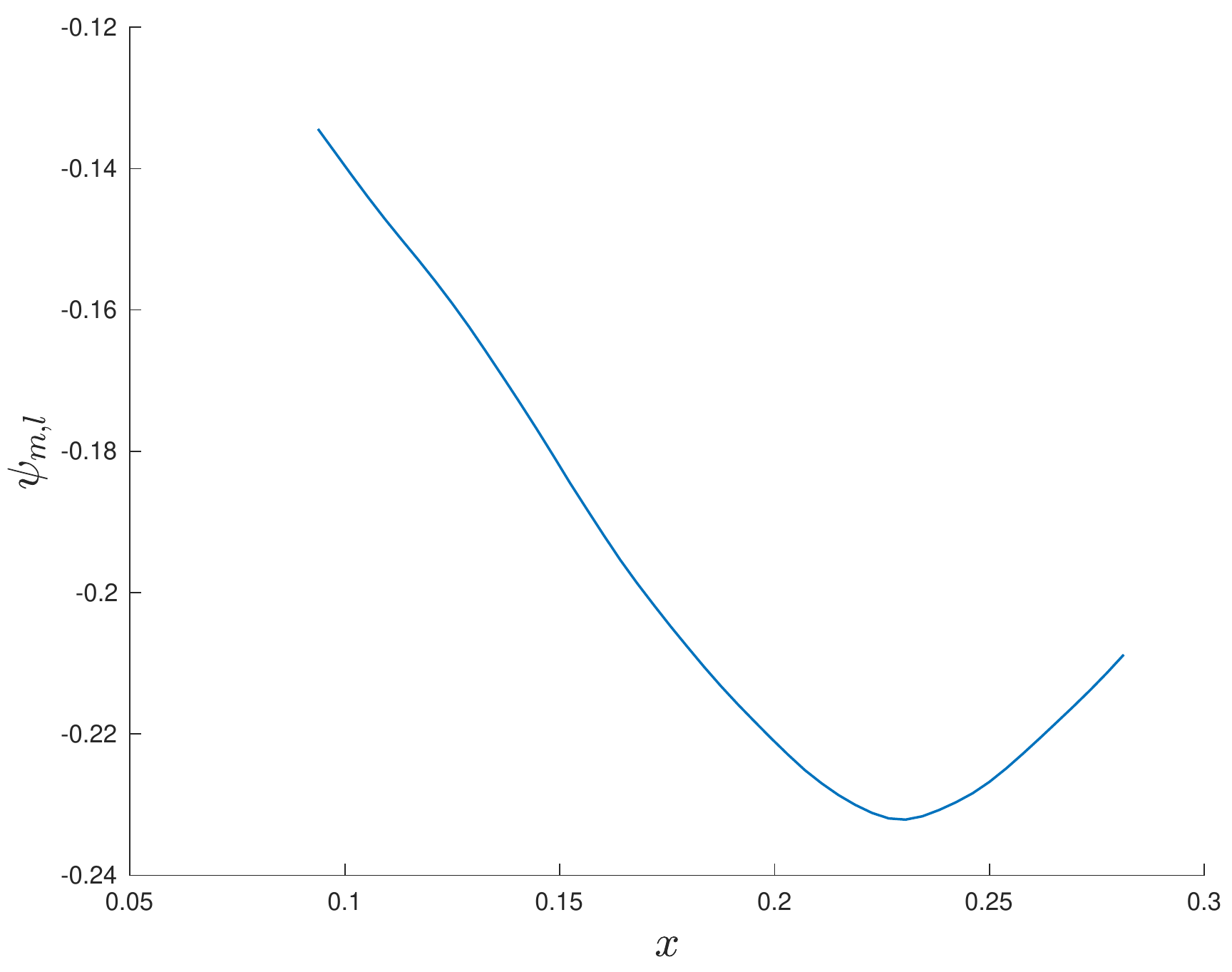}
  \includegraphics[width=0.35\textwidth]{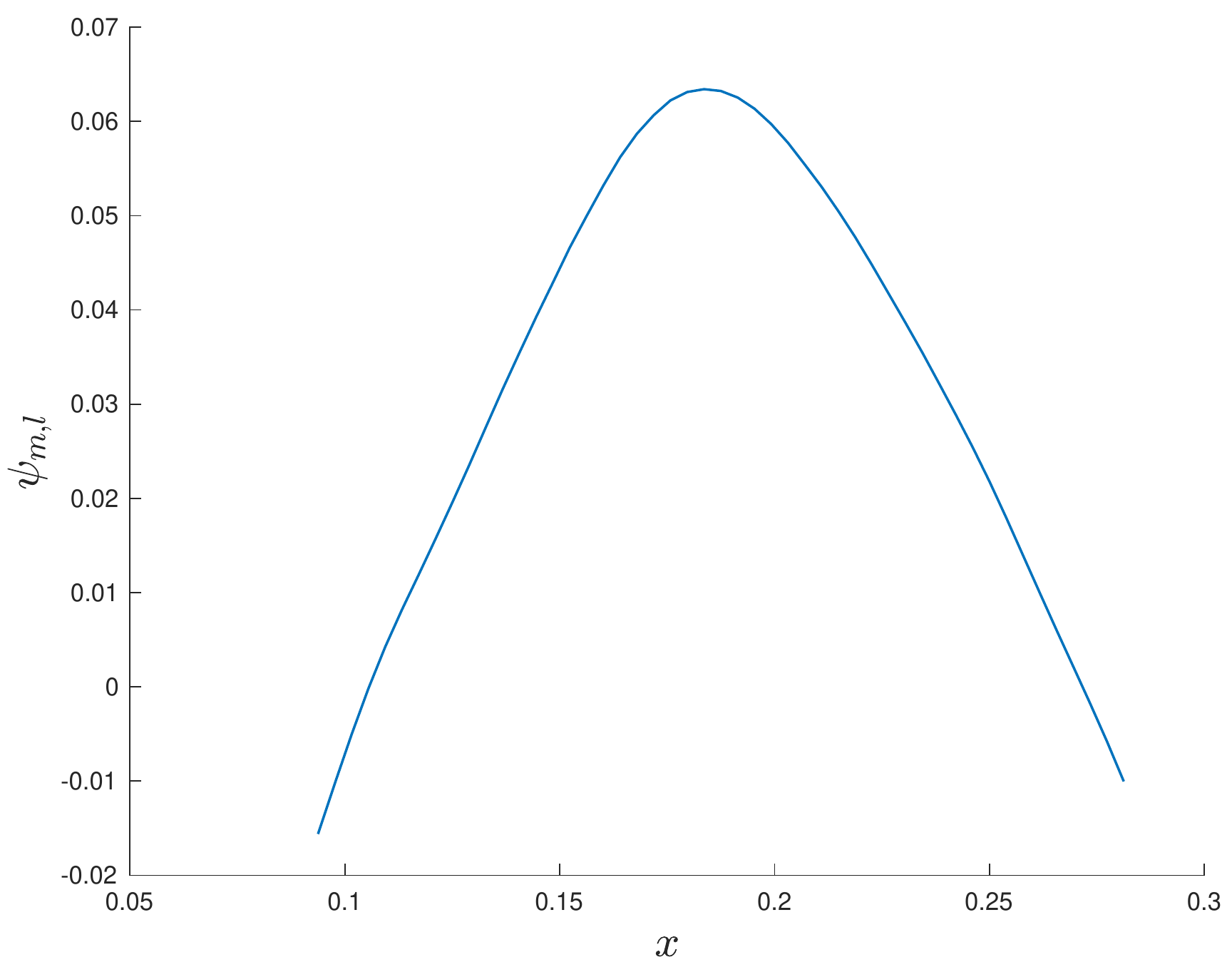}
  \includegraphics[width=0.35\textwidth]{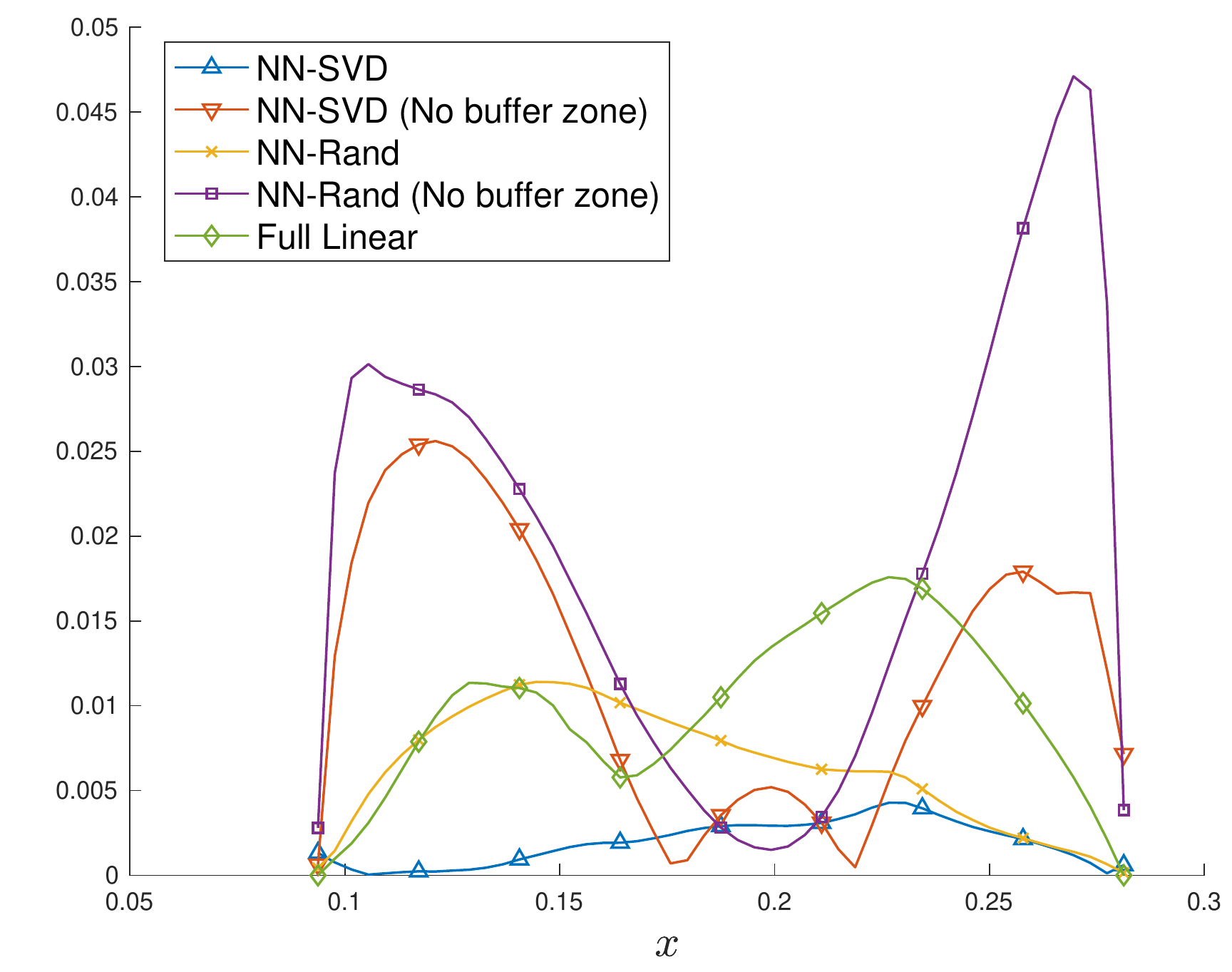}
  \includegraphics[width=0.35\textwidth]{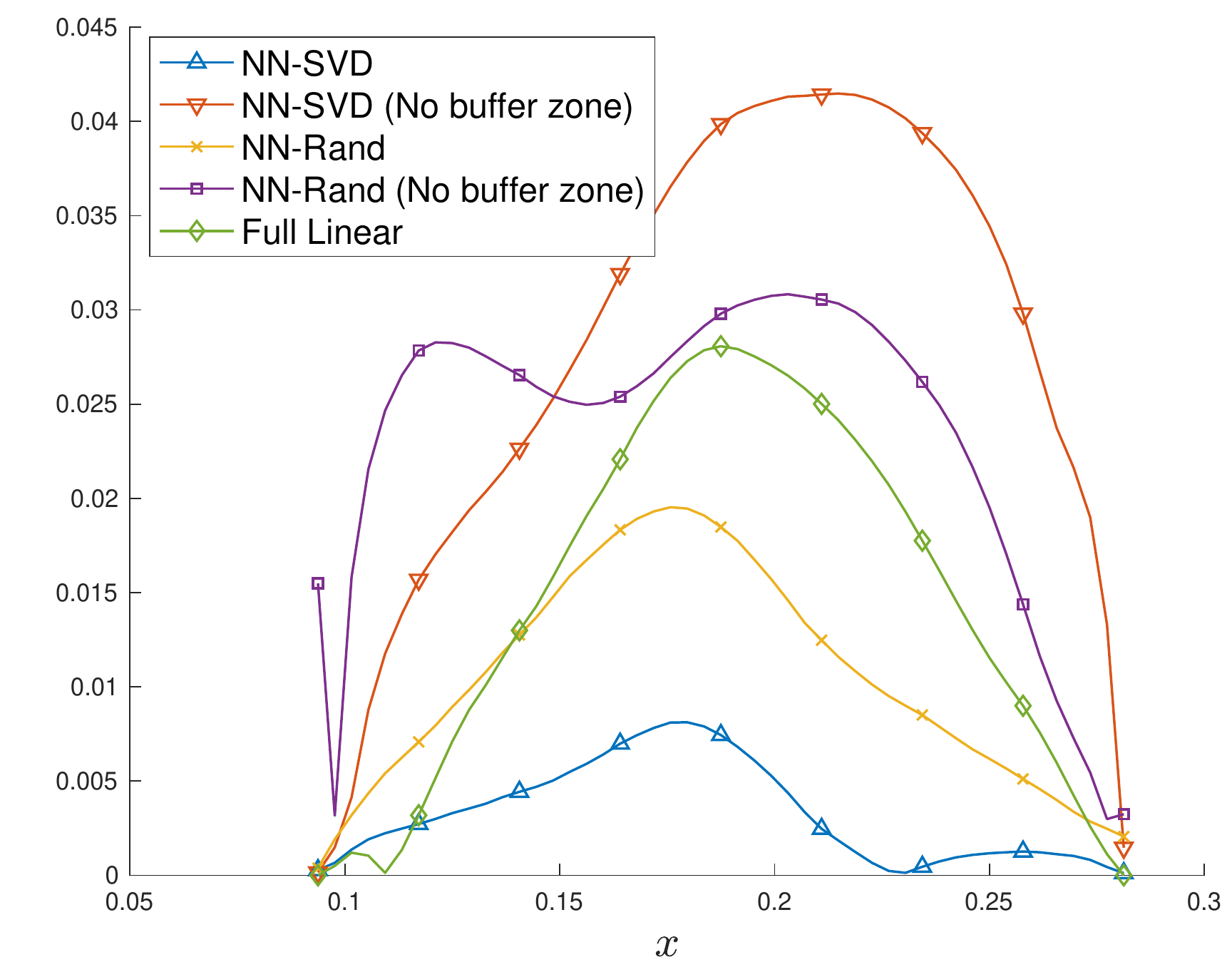}
  \caption{The top row shows the ground truths $\psi_{l,m}$ ($m = (2,2)$, $l = (2,1)$) of two samples in the test set. The bottom row shows the error $|\psi_{l,m}-\widetilde{\psi}_{l,m}|$, where $\widetilde{\psi}_{l,m}$ are computed by the low-rank SVD initialized $\mathcal{Q}_{m}^{\rmNN}$ (with and without buffer zone), randomly initialized $\mathcal{Q}_{m}^{\rmNN}$ (with and without buffer zone), and the linear operator $\mathcal{Q}_{m}^{\rmL}$.}
  \label{fig:psi_examples_pPoi}
\end{figure}

We generate a test data set from the same distribution as the buffered training data set with 100 samples for patch $m=(2,2)$. The test errors~\eqref{eqn:loss_detail} in the training process for different models are plotted in Figure~\ref{fig:Lc_learning_curve_test_pPoi}. As for the training loss, the variant with buffered patches and SVD initialization gives the best results.

To demonstrate generalization performance on this example, we plot the predicted outputs for two typical examples in the test set in Figure~\ref{fig:psi_examples_pPoi}.
For comparison, we also plot the outputs produced by randomly initialized neural network and the linear operator $\mathcal{Q}_m^{\rmL}$.
The low-rank SVD-initialized neural network shows best reconstruction performance.

\begin{figure}[t!]
  \centering
  \includegraphics[width=0.35\textwidth]{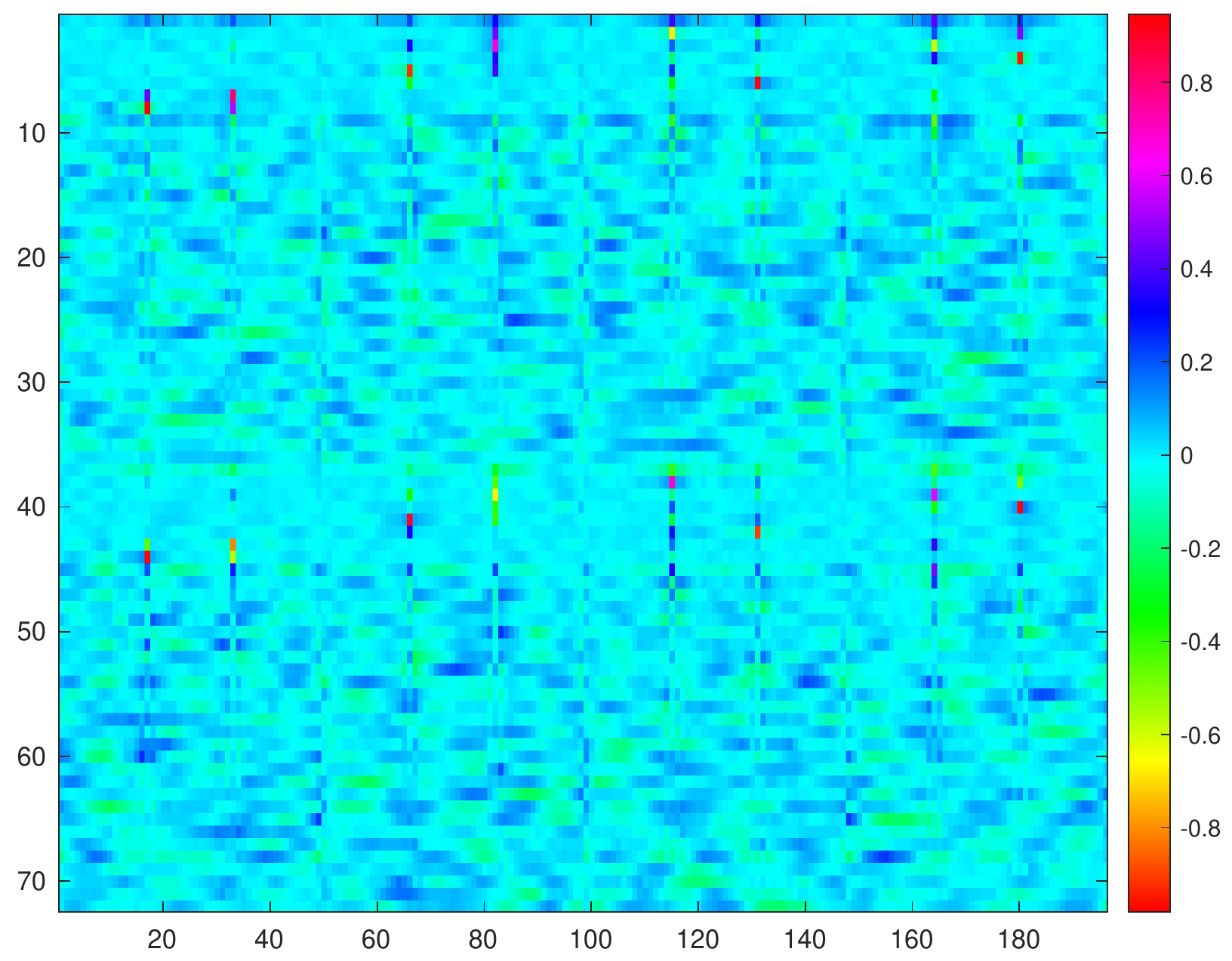}
  \includegraphics[width=0.35\textwidth]{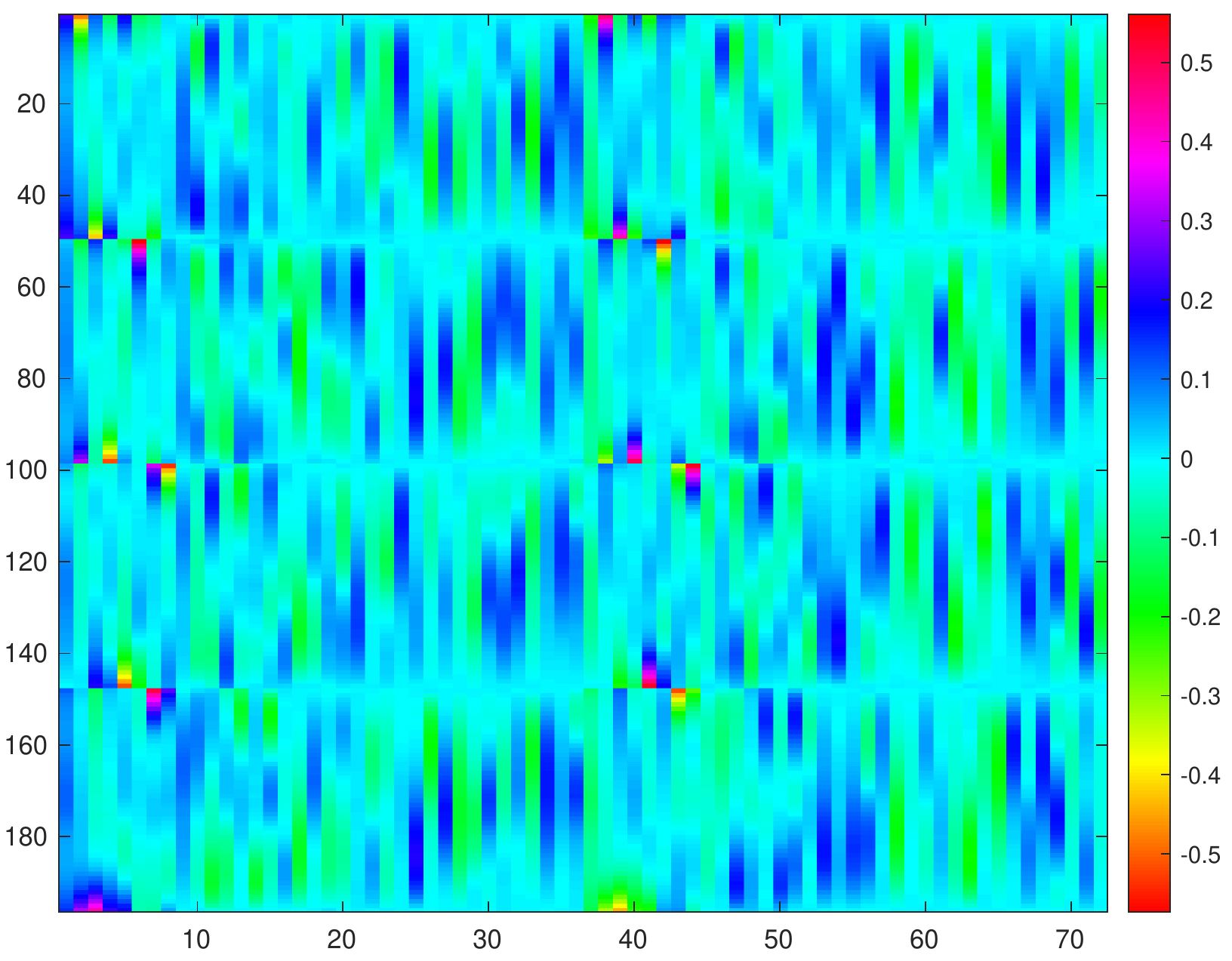}
  \\
  \includegraphics[width=0.35\textwidth]{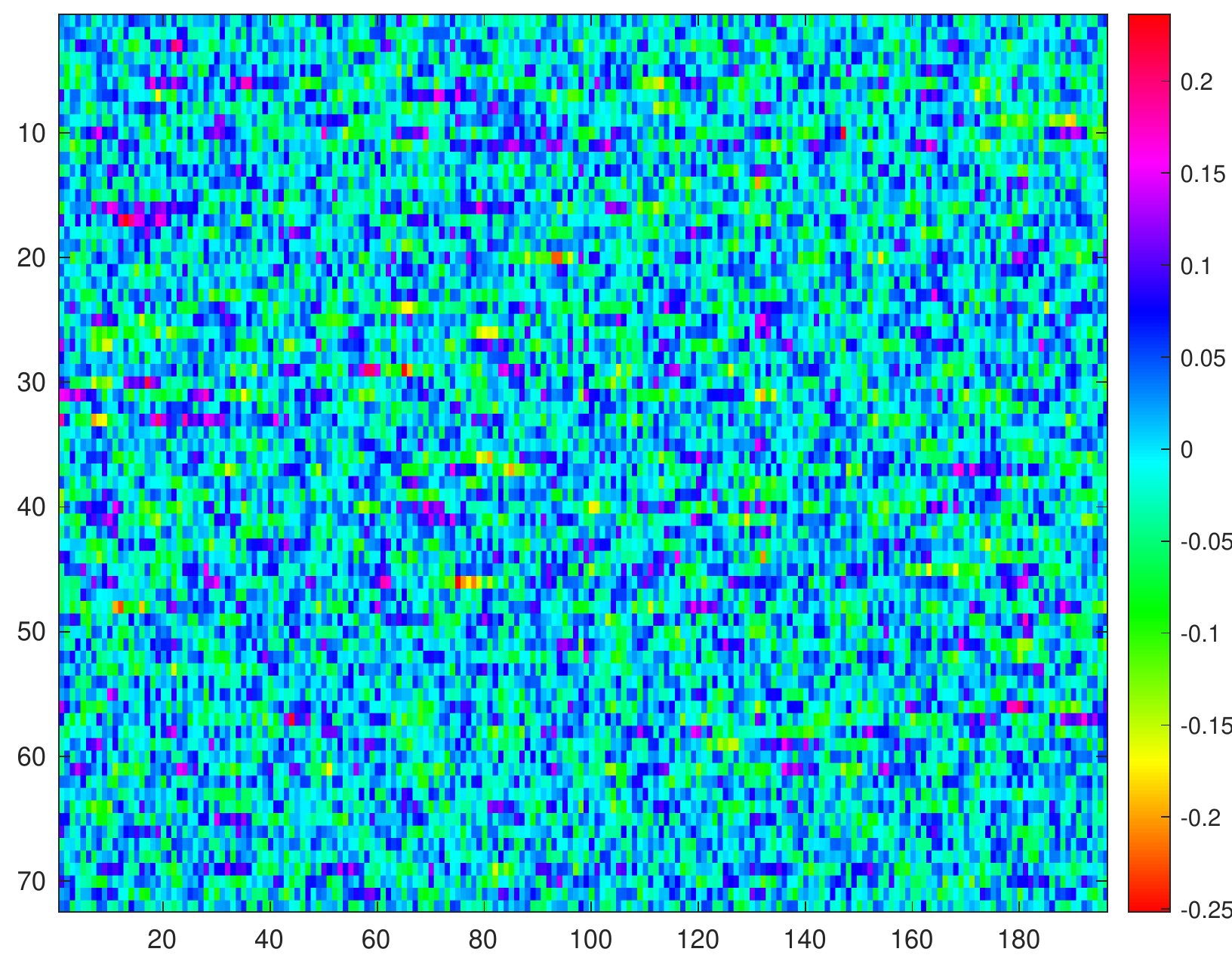}
  \includegraphics[width=0.35\textwidth]{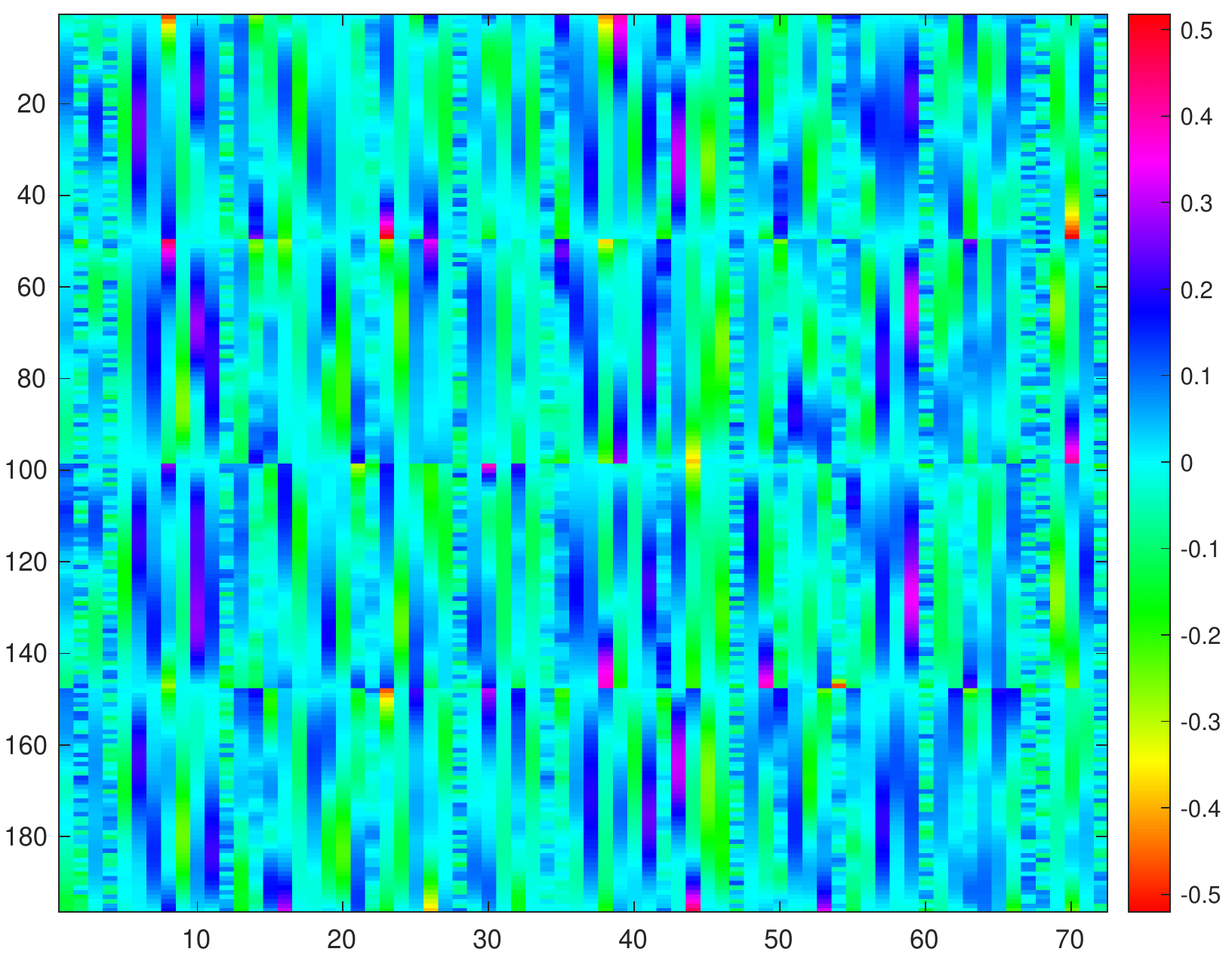}
  \caption{The first row shows the final weight matrices $W_1$ (left), $W_2$ (right) for SVD-initialized model on patch $m = (2,2)$. The second row shows the weight matrices $W_1$ (left), $W_2$ (right) for randomly initialized model on patch $m = (2,2)$. In both cases, training  data is obtained by enlarging the patch.}\label{fig:parameter_pPoi}
\end{figure}

Figure~\ref{fig:parameter_pPoi} show the final weight matrices for models initialized by different methods. All weights have non-trivial values, suggesting that the NN has appropriate dimensions for approximating $\mathcal{Q}^\eps_m$.
Although the structure of the $W_2$ matrix looks roughly similar for each case, the $W_1$ matrices are quite different in character, with the randomly initialized version at bottom left having no obvious structure.

\subsubsection{Schwarz iteration: Online solutions}

\begin{table}[b!]
	\centering
	\begin{tabular}{ c | c }
		\hline \hline
			No.	& Boundary condition \\
		\hline
		1  	& \tabincell {c}{$\phi(x,0) = -\sin(2\pi x)$\,, $\phi(x,1) = \sin(2\pi x)$ \\ $\phi(0,y) = \sin(2\pi y)$\,, $\phi(1,y) = -\sin(2\pi y)$}\\
        \hline
		2	& \tabincell {c}{$\phi(x,0) = -\sin(4\pi x)$\,, $\phi(x,1) = \sin(4\pi x)$ \\ $\phi(0,y) = \sin(4\pi y)$\,, $\phi(1,y) = -\sin(4\pi y)$}\\
        \hline
		3   & \tabincell {c}{$\phi(x,0) = -1$\,, $\phi(x,1) = 1$ \\ $\phi(0,y) = 2y^2-1$\,, $\phi(1,y) = 2y^2-1$} \\
		\hline\hline
	\end{tabular}
    \caption{Boundary conditions for $p$-Laplace equation~\eqref{eqn:pPoi_elliptic} used in the global test.}
  \label{tbl:BCs_pPoi}
\end{table}

\begin{figure}[t!]
  \centering
  \includegraphics[width=0.3\textwidth]{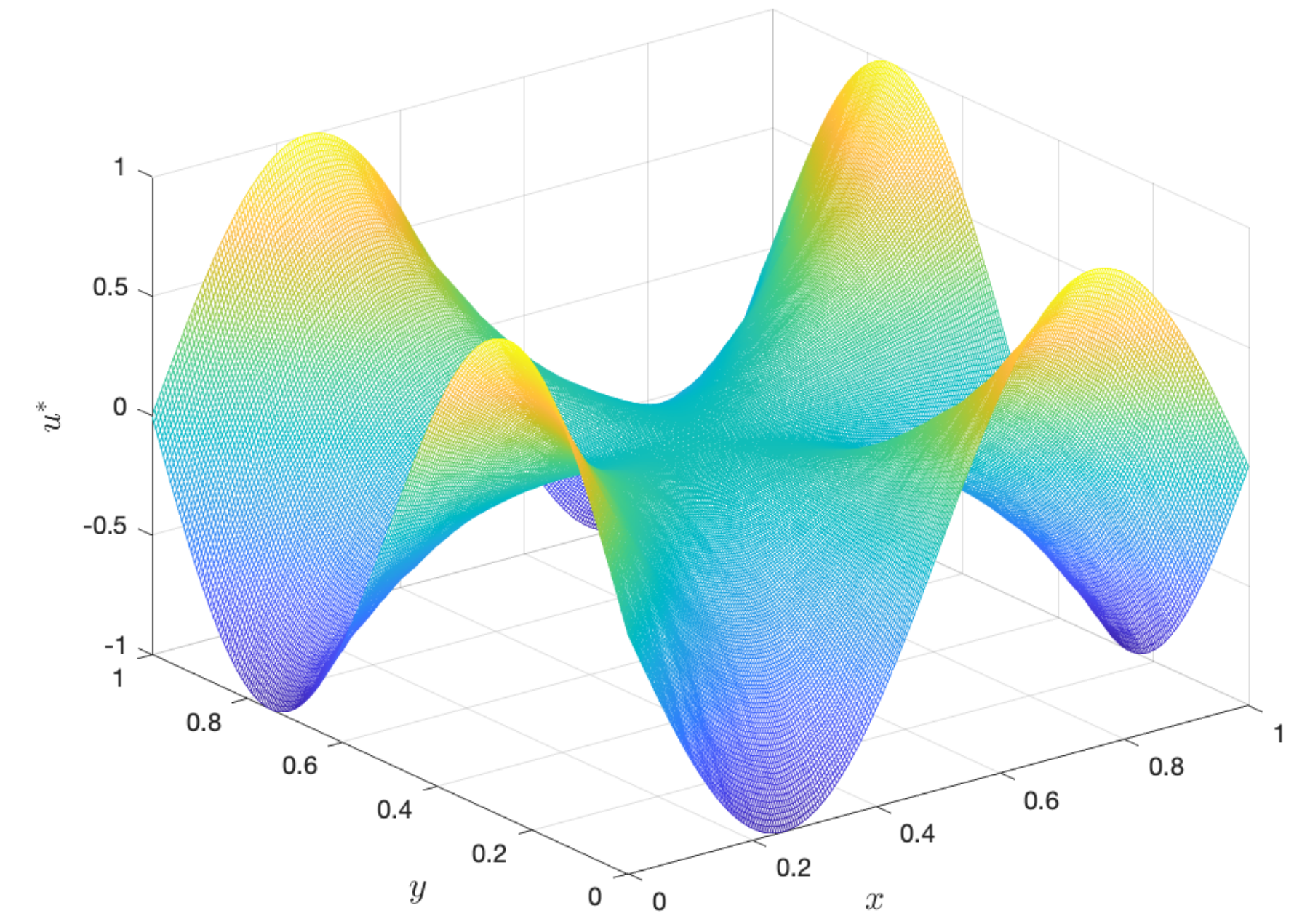}
  \includegraphics[width=0.3\textwidth]{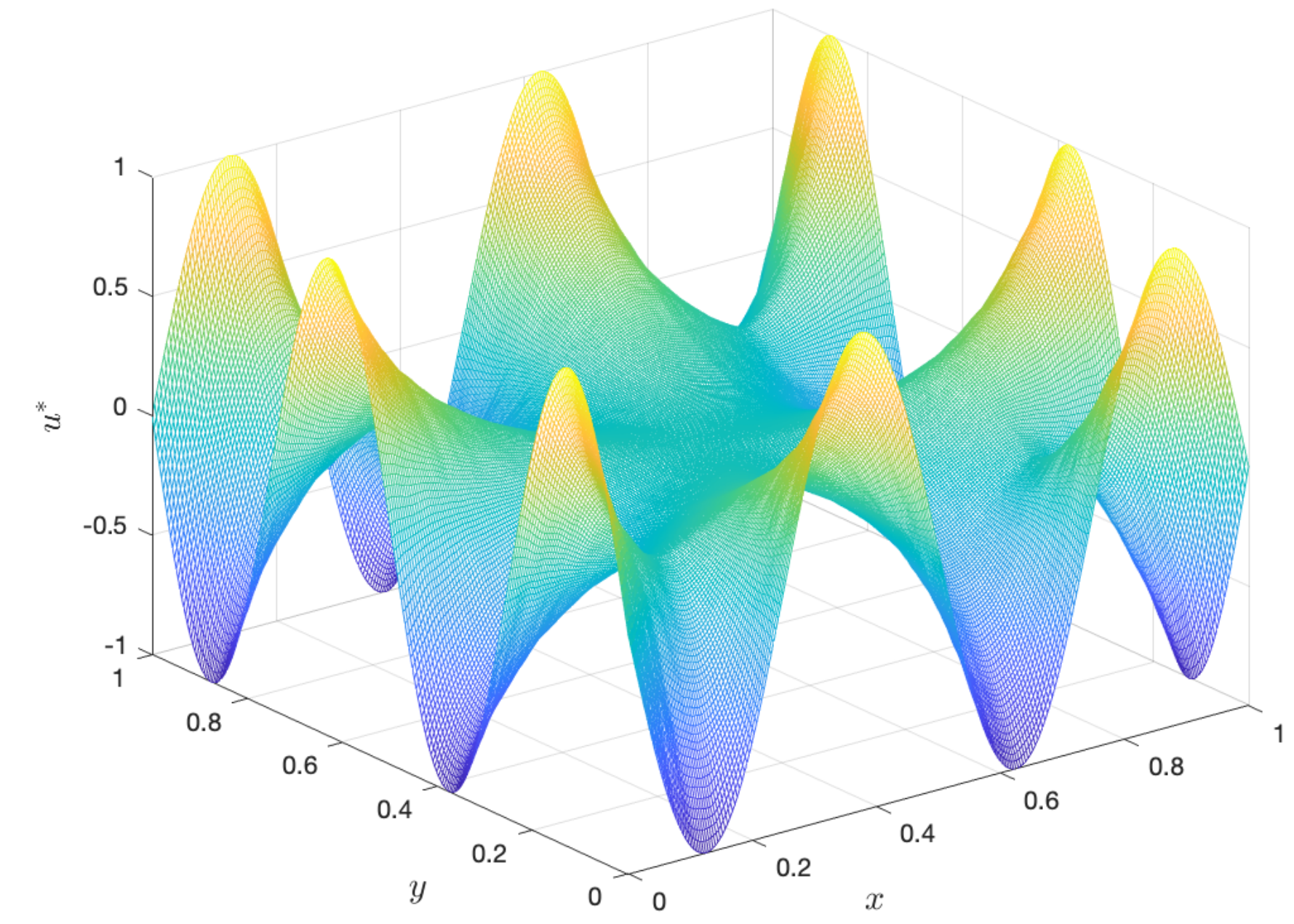}
  \includegraphics[width=0.3\textwidth]{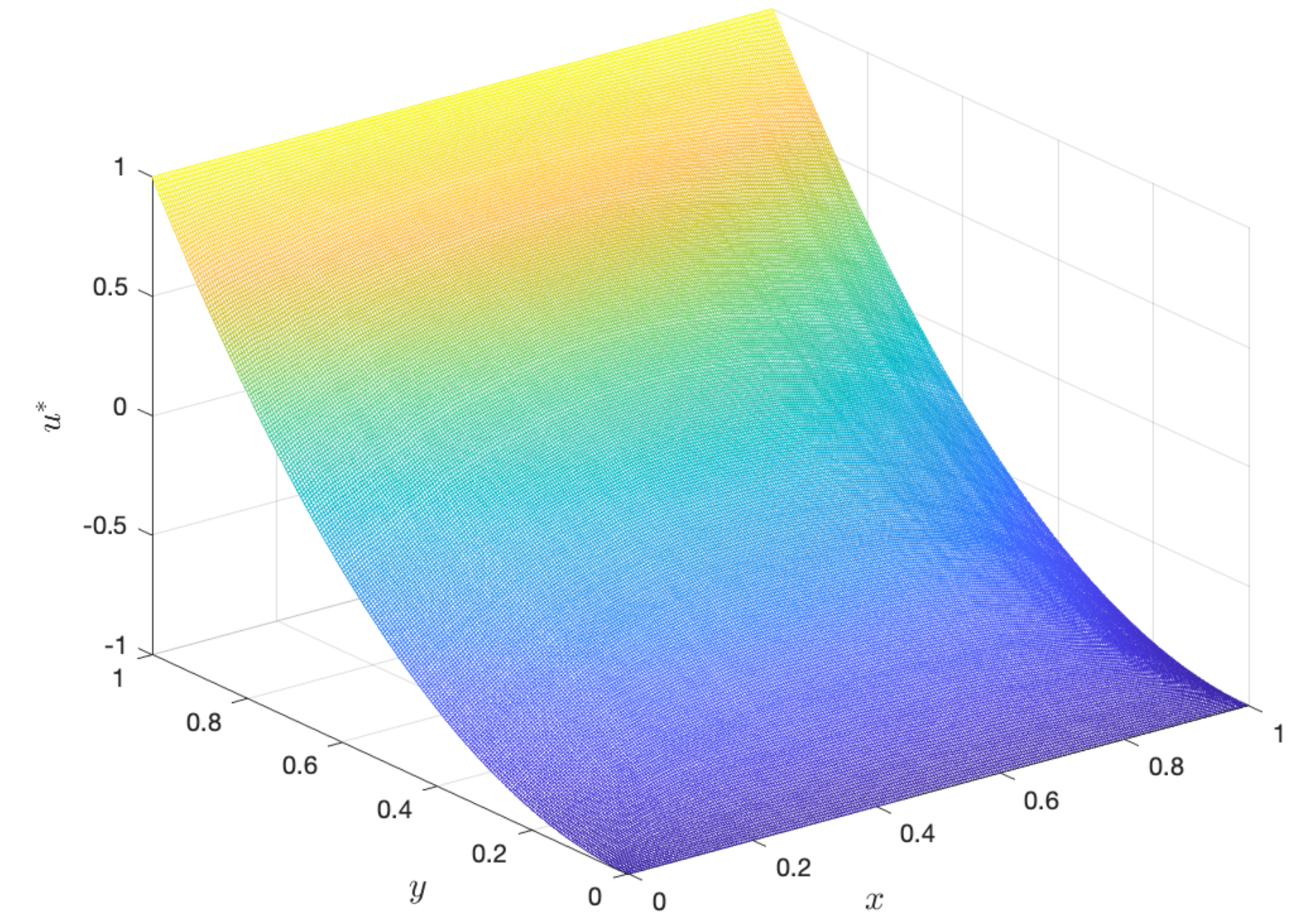}
  \\
  \includegraphics[width=0.3\textwidth]{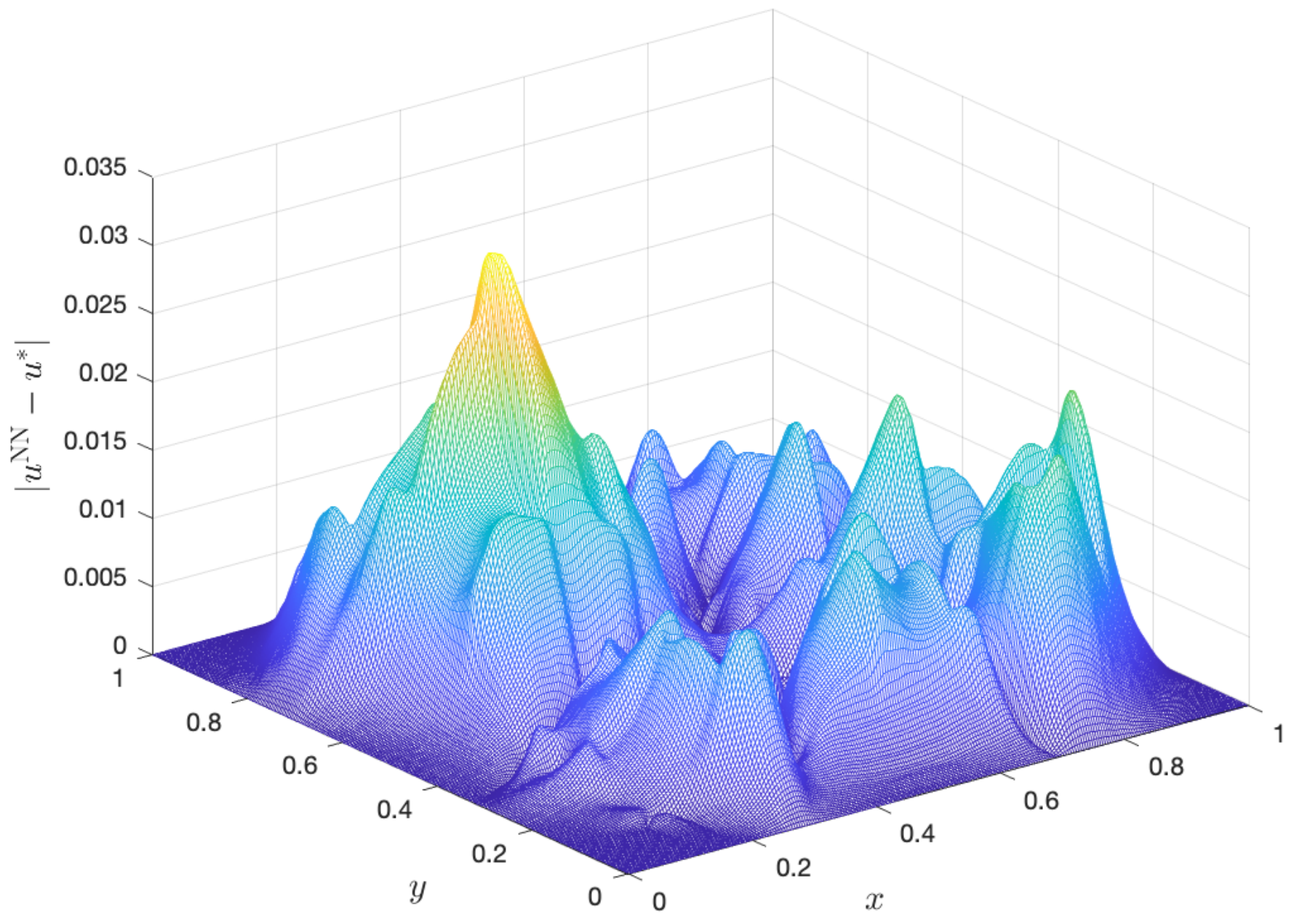}
  \includegraphics[width=0.3\textwidth]{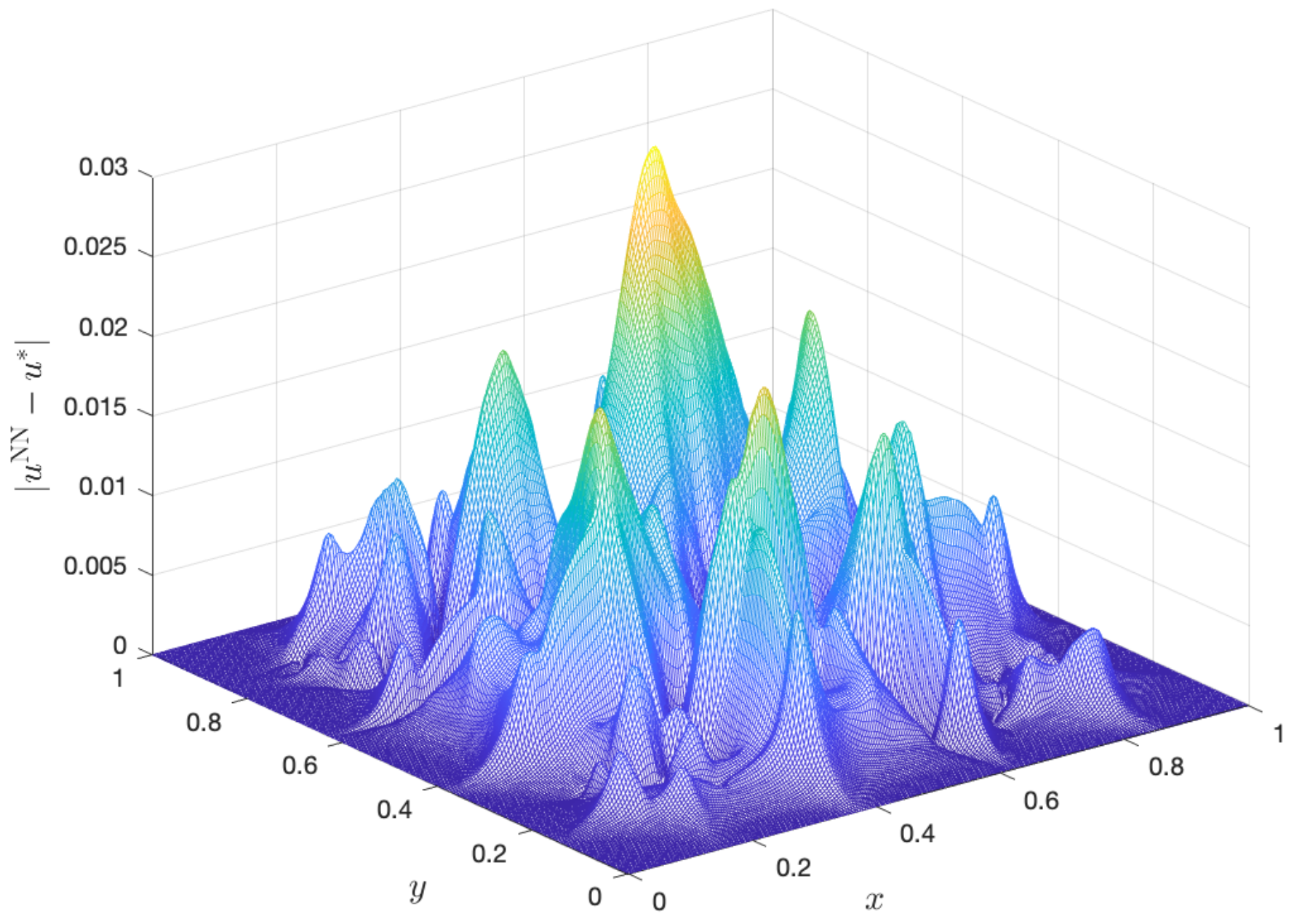}
  \includegraphics[width=0.3\textwidth]{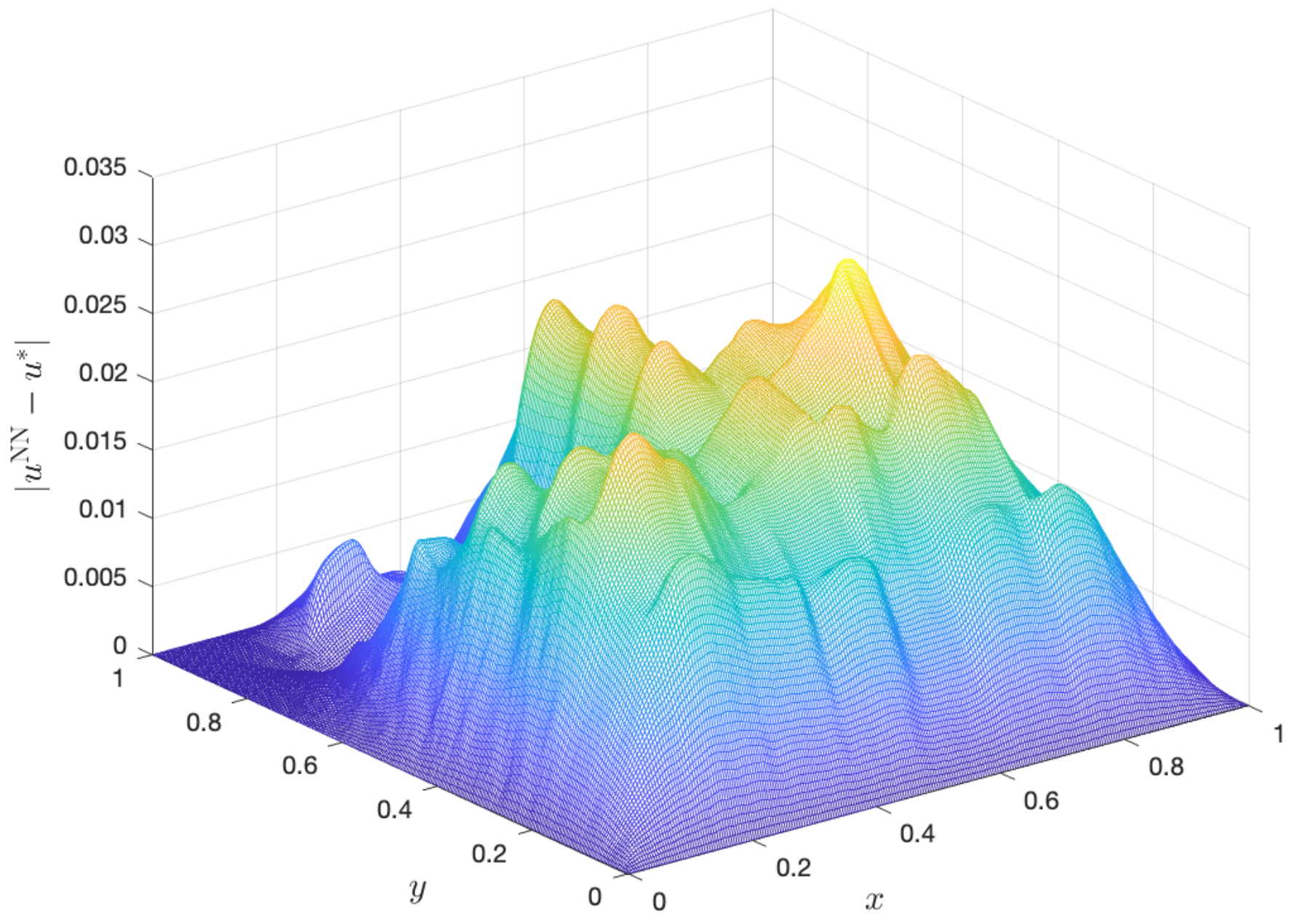}
  \caption{The first row shows the ground truth solution $u^\ast$ for $p$-Laplace equation~\eqref{eqn:pPoi_elliptic} for boundary condition 1 to 3 from left to right. The second row shows the absolute error $|u^\rmNN-u^\ast|$ for boundary condition 1 to 3 from left to right.}
  \label{fig:u_err_pPoi}
\end{figure}

Next, we apply the neural networks to the Schwarz iteration and show the global test performance. In Table~\ref{tbl:BCs_pPoi} we list the boundary conditions for three problems used in the test.
We use tolerance $\delta_0=10^{-4}$ in Algorithm \ref{alg:NN_online} and use the full accuracy  local solvers as in the generation of training  data set. Figure~\ref{fig:u_err_pPoi} shows ground-truth solutions for different boundary conditions and the absolute error of $u^\rmNN$ obtained by neural network-based Schwarz iteration (plotted on a different scale).
Error norms for the different methods can be found in Table~\ref{tbl:Err1_pPoi}.

To demonstrate the efficiency of our method, we compare the CPU time of neural network based-Schwarz method and the classical Schwarz method with tolerance $\delta_0=10^{-4}$ in Algorithm \ref{alg:general}.
The NNs are trained using SVD initialization, with training data generated with buffer zones on the patches.
The local solvers in the reference solution are chosen so that the local accuracy is at the same level as the NN-approximation, making for a fair comparison.
The CPU time, number of iteration and error comparison can be found in Table~\ref{tbl:Runtime_pPoi}.
Compared with the classical Schwarz iteration, the reduced method updates local iterations much faster, while producing $H^1$ errors of the same order.

\begin{table}[t!]
	\centering
	\begin{tabular}{ c |ccc | ccc }
		\hline \hline
			No. BC	& \multicolumn{3}{c|}{1} & \multicolumn{3}{c}{2}\\
		\hline
		Relative Error  & $L^2$ & $H^1$ & $L^\infty$ & $L^2$ & $H^1$ & $L^\infty$\\
        \hline
		SVD-NN & \textbf{0.0199} & \textbf{0.0314} & \textbf{0.0324} & \textbf{0.0171} & \textbf{0.0250} & \textbf{0.0290} \\
        SVD-NN (No buffer zone) & 0.0935 & 0.1793 & 0.1398 & 0.0874 & 0.1052 & 0.1346 \\
		Rand-NN & 0.0280 & 0.0400 & 0.0480 & 0.0260 & 0.0331 & 0.0367 \\
        Rand-NN (No buffer zone) & 0.1062 & 0.1793 & 0.1412 & 0.0696 & 0.1023 & 0.1119 \\
		Linear & 0.0623 & 0.1178 & 0.0909 & 0.0606 & 0.0990 & 0.0751 \\
		\hline\hline
	\end{tabular}
    \caption{Relative error for $p$-Laplace equation~\eqref{eqn:pPoi_elliptic} by different methods.}
  \label{tbl:Err1_pPoi}
\end{table}

\begin{table}[b!]
	\centering
	\begin{tabular}{ c | ccc }
		\hline \hline
			No. BC &\multicolumn{3}{c}{3} \\
		\hline
		Relative Error &  $L^2$ & $H^1$ & $L^\infty$ \\
        \hline
		SVD-NN &  \textbf{0.0215} & \textbf{0.0443} & \textbf{0.0311} \\
        SVD-NN (No buffer zone) &  0.1204 & 0.3331 & 0.2173 \\
		Rand-NN &  0.0241 & 0.0578 & 0.0411 \\
        Rand-NN (No buffer zone) &  0.2748 & 0.4380 & 0.3947 \\
		Linear & 0.1390 & 0.1861 & 0.1624 \\
		\hline\hline
	\end{tabular}
    \caption{Relative error for $p$-Laplace equation~\eqref{eqn:pPoi_elliptic} by different methods. (Continued)}
  \label{tbl:Err2_pPoi}
\end{table}

\begin{table}
	\centering
    \begin{tabular}{ c | cc | cc | cc}
		\hline \hline
		  Problem Number & \multicolumn{2}{c|}{1} & \multicolumn{2}{c|}{2} & \multicolumn{2}{c}{3}\\
        \hline
        Method &  NN & Classical &  NN & Classical &  NN & Classical\\
         \hline
        CPU time & \textbf{35.0} & 87.8 & \textbf{27.8} & 68.3 & \textbf{117.7} & 302.2\\
		Iteration & 52 & 54 & 37 & 38 & 151 & 146\\
		$H^1$ Error & 0.0392 & 0.0231 & 0.0363 & 0.0256 & 0.0457 & 0.0124\\
		\hline\hline
	\end{tabular}
    \caption{CPU time (s), number of iterations and the $H^1$ error of the classical Schwarz iteration and the neural network accelerated Schwarz iteration for $p$-Laplace equation~\eqref{eqn:pPoi_elliptic}.}
  \label{tbl:Runtime_pPoi}
\end{table}

\section{Conclusion} \label{sec:conclusion}

We have presented a reduced-order neural network-based Schwarz method for multiscale nonlinear elliptic PDEs.
In each iteration, the Schwarz method requires evaluation of a boundary-to-boundary map for each of the subdomains (patches).
This map has high dimensional input and output spaces but is compressible due to the existence of a homogenization limit.
A neural network can approximate high-dimensional maps using a number of parameters relaxed  significantly from the dimension of data, and thus is a perfect fit to learn the boundary-to-boundary operator.
Our method trains two-layer neural networks (with many fewer parameters than the input and output dimensions) to learn the boundary-to-boundary operators in an offline stage.
In an online stage, the neural networks serve as surrogates of local solvers in the Schwarz iteration, leading to significant speedup over classical approaches.
Our approach is illustrated with two examples: a semilinear elliptic equation and a $p$-Laplace equation.

\bigskip
\noindent {\bf Acknowledgments.} The work of all authors supported in part by the National Science Foundation via grant DMS-2023239. The work of SW is further supported in part by National Science Foundation via grant 1934612, a DOE Subcontract 8F-30039 from Argonne National Laboratory, and an AFOSR subcontract UTA20-001224 from UT-Austin. The work of SC, ZD and QL is further supported in part by NSF-DMS-1750488 and ONR-N00014-21-1-2140. The authors thank two anonymous reviewers for their detailed, expert reviews and helpful suggestions.

%\clearpage
\bibliographystyle{siam}
\bibliography{ref}

\end{document}